\documentclass{amsart}
\usepackage{amssymb}

\usepackage{amscd}
\usepackage{thmdefs}
\usepackage[dvips]{epsfig}

\input{tcilatex}
\theoremstyle{definition}
\theoremstyle{remark}
\numberwithin{equation}{section}

\input tcilatex

\begin{document}
\title[Labeling Operators]{Applications of Automata and Graphs: Labeling-Operators in Hilbert Space 
\textrm{II}}
\author{Ilwoo Cho and Palle E. T. Jorgensen}
\address{Saint Ambrose Univ., Dep. of Math., 421 Ambrose Hall, 518 W. Locust St.,
Davenport, Iowa, 52803, U. S. A. / University of Iowa, Dep. of Math., 14
McLean Hall, Iowa City, Iowa, 52242-1419, U. S. A.}
\email{chowoo@sau.edu, jorgense@math.uiowa.edu}
\thanks{}
\date{Mar., 2008}
\dedicatory{}
\thanks{The second named author was supported by the U. S. National Science
Foundation.}
\maketitle

\begin{abstract}
We introduced a family of infinite graphs directly associated with a class
of von Neumann automaton model $\mathcal{A}_{G}.$ These are finite state
models used in symbolic dynamics: stimuli models and in control theory. In
the context of groupoid von Neumann algebras, and an associated fractal
group, we prove a classification theorem for representations of automata.
\end{abstract}

\strut

\textbf{Table of Contents}

\strut

$
\begin{array}{lll}
\text{{\small 1. Introduction}} &  & \,\,\text{{\small 2}} \\ 
\text{{\small 2. Background}} &  & \,\,\text{{\small 7}} \\ 
\text{{\small \quad 2.1. Amalgamated Free Probability}} &  & \,\,\text{%
{\small 7}} \\ 
\text{{\small \quad 2.2. Basic Concepts}} &  & \,\,\text{{\small 9}} \\ 
\text{{\small 3. Definitions}} &  & \text{{\small 12}} \\ 
\text{{\small \quad 3.1. Categorial Groupoids \& Groupoid Actions}} &  & 
\text{{\small 12}} \\ 
\text{{\small \quad 3.2. Automata and Fractal Groups}} &  & \text{{\small 14}%
} \\ 
\text{{\small \quad 3.3. Right Graph von Neumann Algebras}} &  & \text{%
{\small 18}} \\ 
\text{{\small \quad 3.4. }}M\text{{\small -Valued Right Graph }}W^{*}\text{%
{\small -Probability Spaces}} &  & \text{{\small 21}} \\ 
\text{{\small 4. Labeled Graph Groupoids and Graph Automata}} &  & \text{%
{\small 26}} \\ 
\text{{\small \quad 4.1. Labeled Graph Groupoids}} &  & \text{{\small 27}}
\\ 
\text{{\small \quad 4.2. The Operation }}\theta \text{{\small \ and }}\omega
_{+} &  & \text{{\small 34}} \\ 
\text{{\small \quad 4.3. Graph Automata}} &  & \text{{\small 35}} \\ 
\text{{\small \quad 4.4. Graph-Automata Trees}} &  & \text{{\small 38}} \\ 
\text{{\small 5. Labeling Operators of Graph Groupoids}} &  & \text{{\small %
43}} \\ 
\text{{\small \quad 5.1. The Labeling Operator of }}\Bbb{G}\text{{\small \
in }}M_{G} &  & \text{{\small 43}} \\ 
\text{{\small \quad 5.2. }}D_{G}\text{{\small -Valued Free Distributional
Data of Labeling Operators}} &  & \text{{\small 45}} \\ 
\text{{\small \quad 5.3. Graph-Theoretical Characterization of }}\mathcal{W}%
_{2n}^{(m)}\text{{\small \ and }}\mathcal{W}_{2n}^{(c)} &  & \text{{\small 52%
}} \\ 
\text{{\small 6. Examples}} &  & \text{{\small 54}} \\ 
\text{{\small References}} &  & \text{{\small 57}}
\end{array}
$

\strut

\strut

\strut \strut \strut \strut

\section{Introduction\strut}

\strut

This paper is the second of two papers studying representations of graphs by
operators in Hilbert space, and their applications, with the help of
Automata Theory. While graph theory is traditionally considered part of
discrete mathematics, in this paper we show that applications of tools from
automata and operators on Hilbert spaces yield global results for
representations of a class of infinite graphs, as well as spin-off
applications. We begin with an outline of the use of automata, and more
generally, of finite state models (FSMs) in the processing of numbers, or
more importantly in sampling and in quantization of digitized information
such as speech signals and digital images. In these models, the finite input
states of a particular FSM might be frequency-bands (for example a
prescribed pair of high-pass and low-pass digital filters), or a choice of
subdivision filters; where the subdivision refers to data sets with
self-similarity; such as is typically seen in fractals. Hence, these
applications make connections to discrete wavelet algorithms as used in
signal and image processing, as well as in science and engineering. If the
input-options for a particular FSM are chosen from a prescribed system of
low-pass and high-pass filters, the resulting discrete model can then be
realized by operators on Hilbert spaces. Similarly, images are digitized
into matrix shapes used in computer programs for compression of images.\strut

In a general operator theoretic setting, this paper introduces the relevant
representations of the generators of graphs and automata. Hence data from
FSMs and graphs are represented with the use of Hilbert space geometry.
Recall that digital images are typically given by matrices of pixels, and
that spectral analysis and matrix operations can be done with operators in
Hilbert space. However, the Hilbert spaces needed in a particular
application are typically not immediately apparent from the particular
engineering problem under discussion. This paper focused on making the link
between graphs and automata on the discrete side to Hilbert space operators
and representations on the spectral side.\strut

The list of the advantages deriving from the use of operator algebraic
methods in the analysis of graphs includes the following: starting with a
particular graph $G,$ a main question is that of identifying decompositions
and global invariants for $G.$ Motivated by ergodic theory, a further goal
is to analyze global properties of $G$ in terms of its local features, and
an operator algebraic decomposition is a powerful tool in the identification
of local vs global features and invariants of $G.$ For the purpose, one may
use that the theory of von Neumann algebras (rings of operators) allows
factor decompositions (See [66]). And the von Neumann factors then further
allow type classification, type $I,$ $II,$ and $III,$ with sub-types; See
[66] for the fundamentals in the theory of von Neumann algebras, and [64],
[65] for some relevant applications to infinite graphs of this idea. In
addition, the analysis of graphs uses infinite determinant, and the
Fuglede-Kadison determinant [63] (from von Neumann algebra theory) is a
powerful tool in graph computation because it passes to direct Hilbert space
decompositions.\strut \strut \strut \strut \strut \strut 

The potential-theoretic part of the subject may be understood as a
mathematical idealization of electrical networks on infinite graphs $G$ (See
[43], [44] and [45]). We address two issues from Analysis: Find
representations (by operators on Hilbert spaces) of the graph systems $G,$
and identify a class of operators whose spectral theory captures significant
information about $G.$ We focus on the graphs themselves, and our motivation
derives in part from work by Strichartz and Kigami (See [46], [47] and [48])
and others for restricted classes of fractals. The focus there is the
adaptation of a rescaling and energy-renormalization on graphs and an
adaptation to fractal models $X$; for example, to Brownian motion on $X,$ or
to a version of differential operators on spaces of functions on $X.$ Here,
our focus is on Operator Theory needed for analysis of graphs and automata
(finite-state machines: e.g., [54]) such as are used in for example signal
processing algorithm, e.g., those based on a discrete multiresolution (e.g.,
see [52]).

\strut

We begin our paper with an outline of three trends: (i) Hilbert space and
Spectral Analysis on the graphs, (ii) Analysis on associated fractals
arising from automata, and finally (iii) the interplay between (i) and (ii).
We further stress some of the differences between the two.

\strut

\strut The main purpose of this paper is to introduce a new algebraic
structures having certain fractal property. In [10], [11], [13], [14] and
[15], we introduced graph groupoids induced by countable directed graphs. A
graph groupoid is a categorial groupoid having as a base the set of all
vertices of the given graph. We know that every groupoid having only one
base element is a group. So, if $G$ is a finite directed graph with graph
groupoid $\Bbb{G},$ and if the vertex set $V(G)$ $\subset $ $\Bbb{G}$
consists of only one element, then the graph groupoid $\Bbb{G}$ of this
one-vertex-$n$-loop-edge graph $G$, for $n$ $=$ $\left| E(G)\right| $ $\in $ 
$\Bbb{N},$ is a group; futhermore, it is group-isomorphic to the free group $%
F_{n}$ with $n$-generators (See [10] and [11]). Here, $E(G)$ denotes the
edge set of $G.$

\strut

Then we can embed the edge set $E(G)$ into the edge set $E(\widehat{G}),$
where $\widehat{G}$ is the shadowed graph of $G.$ Put the lattice weights $%
l_{\pm 1},$ ..., $l_{\pm n}$, associate to the edges $e_{1}^{\pm 1},$ ..., $%
e_{n}^{\pm 1}$ in $E(\widehat{G}),$ where $e_{j}^{-1}$ $\in $ $E(\widehat{G}%
) $ means the shadow of $e_{j}$ $\in $ $E(G)$ (See Section 2), where $l_{\pm
j} $ $=$ $(0,$ $\pm e^{j})$ in $\Bbb{R}^{2}.$ Here, $e$ is the natural
number in $\Bbb{R}$. More precisely, we define the labeling set $\mathcal{X}%
_{0}$ by $\pm X$ $=$ $\{l_{\pm 1},$ ..., $l_{\pm n}\}.$ Then define the
labeling map

$\strut $

\begin{center}
$\Bbb{\varphi }$ $:$ $\mathcal{X}_{0}$ $\times $ $E(\widehat{G})$ $%
\rightarrow $ $\mathcal{X}_{0}$
\end{center}

by

\begin{center}
$\Bbb{\varphi }\left( l,\text{ }e\right) $ $=$ the weight $l_{e}$ of $e.$
\end{center}

\strut

Also, define the shifting map

$\strut $

\begin{center}
$\Bbb{\psi }$ $:$ $\mathcal{X}_{0}$ $\times $ $E(\widehat{G})$ $\rightarrow $
$E(\widehat{G})$
\end{center}

by

\begin{center}
$\Bbb{\psi }\left( l,\text{ }e\right) $ $=$ $e.$
\end{center}

Then we have the automaton $\mathcal{A}_{G}$ $=$ $<\pm X_{0},$ $E(\widehat{G}%
),$ $\Bbb{\varphi },$ $\Bbb{\psi }>.$ Then this automaton $\mathcal{A}_{G}$ $%
=$ $<\pm X_{0},$ $E(\widehat{G}),$ $\varphi ,$ $\psi >$ generates a fractal
group $\Gamma $ $\overset{\text{Group}}{=}$ $F_{n}$. We further address the
general case where the vertex set contains more than one element? The answer
of this question was provided in [17]. In this paper, we will consider more
general case than this.

\strut

Our paper is organized as follows. The remaining of Introduction contains an
overview, including some definitions; meant only as sketch. Full definitions
and statements of theorems will be in the respective main sections of the
paper.

\strut

Section 2 introduces those fundamentals from amalgamated free probability
theory which will be needed later. Here, we use the term ``free
probability'' in the sense of D. Voiculescu (e.g., See [5]), a framework
where random variables are noncommutative operators in Hilbert space, as
opposed to commuting measurable functions on a probability space. While
Voiculescu's notion free independence is modeled on the classical concept of
independence of families of random variables, there are important
differences. They are reflected in the other notions we need from the
noncommutative setting: operator-valued $*$-moments and $*$-cumulants. In
addition, we will need fundamentals from the theory of partitions. This
material allows applications of operator algebra tools to our present graph
analysis of automata. Moreover, within this framework, we will need an
extension of the notion of group crossed products from operator algebra
theory: we introduce a (goupoid) crossed product of a von Neumann algebra by
a groupoid (See Definition 2.1).

\strut

Section 3 focuses on graphs and graph groupoids: categorial groupoids and
groupoid actions. Here, actions refers to actions by transformations on von
Neumann algebras.

\strut

A main application will be to von Neumann's automata (i.e., certain
self-reproducing dynamical structures) introduced in Section 3.2. From each
automaton $\mathcal{A}$, we define an automata group $\Gamma (\mathcal{A}),$
and we identify a subclass of $\mathcal{A}$ for which the corresponding
group is a fractal (Definition 3.1). The rest of Section 3 deals with
actions on von Neumann algebras which are induced by countable directed
graphs. Starting with the graph, in Theorem 3.2, we show that a groupoid
crossed product von Neumann algebra $\Bbb{M}_{G}$ has an amalgamated free
structures over its diagonal $W^{*}$-subalgebra $\Bbb{D}_{G}.$ We then use
this in Theorem 3.4 and 3.5 in breaking up a given $W^{*}$-probability space
into basic constituents. Specially, we show that every groupoid crossed
product algebra $\Bbb{M}_{G},$ as a $W^{*}$-probability space, decomposes in
a specific way into a direct decomposition in such a manner that the
constituents in the decomposition are so-called $\Bbb{D}_{G}$-free blocks.
Moreover, Theorem 3.6 offers a classification of the $\Bbb{D}_{G}$-free
blocks which may occur. Applications, Example 3.5 and 3.6 are given to edge
graphs, while Section 4 covers graph-automata (Definition 4.7). In our
further analysis of automata-labeled graph groupoids and their von Neumann
algebras, this is brought to bear. In Section 4.4, we specialize to
graph-automata trees, spelling out more detailed results.

\strut

The theorems from Section 3 and 5 are strengthened in the remaining two
sections of our paper with the use of labeling operators. Theorem 5.4 makes
use of\ a condition for labeling operators in deciding $\Bbb{D}_{G}$%
-freeness, while Theorem 5.6 and 5.8 yield explicit formulas for $\Bbb{D}_{G}
$-valued joint moments, and joint cumulants, respectively. This is then
further refined, and applied in Section 6, the fractal case in Example 6.1.

\strut

We realize that, by giving weights on the edges of a given graph $G$, which
are determined by the out-degrees of vertices, we can get the corresponding
automaton $\mathcal{A}_{G}$. By observing the properties of $\mathcal{A}%
_{G}, $ we can determine the groupoidal version of a fractal property,
relative to a fractal group. We will say that a graph groupoid $\Bbb{G}$ is
a fractaloid, if $\mathcal{A}_{G}$ acts ``fully'' on a $\left| V(\widehat{G}%
)\right| $-copies of regular growing directed trees. For example, if a
directed graph $G$ is a one-flow circulant graph, then the graph groupoid $%
\Bbb{G}$ is a fractaloid. How about the given graph groupoids are not
fractaloids? This is our motivation of this work.

\strut \strut \strut

Similar to [10] and like [17], we can have a right graph von Neumann algebra 
$\Bbb{M}_{G}$ $=$ $\overline{\Bbb{C}[\beta (\Bbb{G})]}^{w},$ as a $W^{*}$%
-subalgebra of the operator algebra $B(H_{G}),$ consisting of all bounded
linear operators on $H_{G},$ where $H_{G}$ is a generalized Fock space,
called the graph Hilbert space induced by $\Bbb{G}.$ We are interested in
the case where $\Bbb{G}$ is a fractaloid. Then, similar to the classical
case, we can define the Hecke-type operator $\tau $ $\in $ $\Bbb{M}_{G}.$
Instead of considering the pure operator-theoretical data of $\tau ,$ we
will observe the amalgamated free distributional data of $\tau .$ Since $%
\tau $ is self-adjoint, the amalgamated free moments of it contains the
operator-valued spectral measure theoretical information. This means that
the free moments of $\tau $ will contain the operator-theoretical data of $%
\tau .$ In [17], we observed the special cases where the given graph
groupoids are fractaloids.

\strut \strut

From the theory of algebras of operators in Hilbert space, we will need von
Neumann algebra constructions (e.g., [59]), free probability, in particular,
amalgamated free products (e.g., [5] and [23]), groupoids and groupoid
actions (e.g., [19] and [60]), and Hecke operators (e.g., [27]). In
technical discussions, we will use ``von Neumann algebras'' and ``$W^{*}$%
-algebras'' synonymously.. If $H$ is a given Hilbert space, the algebra of
all bounded operators on $H$ is denoted by $B(H).$ From Graph Theory, we use
such notations as ``sets of vertices'', ``sets of edges'', ``loops'' and
``degrees'', etc. If $G$ is a given directed graph, we introduce a
``shadow'' construction from a reversal of edges, denoted by $G^{-1}$; see
details below. From symbolic dynamics, we shall use fundamentals of
``automata'', as well as free constructions, such as the free group with
multi-generators. Hence, our next section will contain a number of
definitions that will needed later.

\strut

Since our paper is interdisciplinary and directed at several audiences, we
have included details from one area of mathematics which might not be
familiar to readers from another.

\strut \strut \strut \strut \strut \strut \strut \strut \strut

A \emph{graph} is a set of objects called \emph{vertices} (or points or
nodes) connected by links called \emph{edges} (or lines). In a \emph{%
directed graph}, the two directions are counted as being distinct directed
edges (or arcs). A graph is depicted in a diagrammatic form as a set of dots
(for vertices), jointed by curves (for edges). Similarly, a directed graph
is depicted in a diagrammatic form as a set of dots jointed by arrowed
curves, where the arrows point the direction of the directed edges.

\strut

Throughout this paper, every graph is a locally finite countably directed
graph. Recall that we say that a countably directed graph $G$ is \emph{%
locally finite} if each vertex of $G$ has only finitely many incident edges.
Equivalently, the degree of $v$ is finitely determined. Also, recall that
the degree $\deg (v)$ of a vertex $v$ is defined to be the sum of the
out-degree $\deg _{out}(v)$ and the in-degree $\deg _{in}(v),$ where

\strut

\begin{center}
$\deg _{out}(v)$ $\overset{def}{=}$ $\left| \{e\in E(G):e\text{ has its
initial vertex }v\}\right| $
\end{center}

and

\begin{center}
$\deg _{in}(v)$ $\overset{def}{=}$ $\left| \{e\in E(G):e\text{ has its
terminal vertex }v\}\right| .$
\end{center}

\strut

Let $\widehat{G}$ be the shadowed graph of $G$, in the sense of [10]. Then
we can consider the degree of each vertex of $\widehat{G},$ too, since $%
\widehat{G}$ is also a locally finite countable directed graph. Assume that

$\strut $

\begin{center}
$N$ $=$ $\max $ $\{\deg _{out}(v)$ $:$ $v$ $\in $ $V(\widehat{G})$ $=$ $%
V(G)\}.$
\end{center}

\strut

Then we can define the labeling set $X$ $=$ $\{1,$ $2,$ ..., $N\}.$ We
assign the weights $\{1,$ ..., $N\}$ to all elements in the edge set $E(%
\widehat{G})$ of the shadowed graph $\widehat{G}$ of $G.$ (This weighting
provides the weights of all elements in $\Bbb{G}$, which are the sequences
contained in $X_{0}^{\,\infty },$ where $X_{0}$ $=$ $\{0\}$ $\cup $ $X.$) We
will call this process placing the weights onto all elements of $\Bbb{G},$
the labeling process. This labeling process lets us construct the automaton $%
\mathcal{A}_{G}$ $=$ $<X_{0},$ $\Bbb{G},$ $\varphi ,$ $\psi $ $>$ induced by
the graph groupoid $\Bbb{G}.$ If the automaton $\mathcal{A}_{G}$ satisfies
certain \emph{fractal} property; we will call the graph groupoid $\Bbb{G}$ a
fractaloid. Clearly, the word ``fractaloid'' hints at ``fractal (graph)
groupoid''.

\strut \strut \strut

As in [10] and [11], we will fix a representation $(H_{G},$ $\beta )$ of a
graph groupoid $\Bbb{G},$ where $H_{G}$ is the graph Hilbert space and $%
\beta $ is a certain groupoid action of $\Bbb{G}.$ Let $X$ $=$ $\{x_{1},$ $%
x_{2},$ ..., $x_{N}\}$ be the labeling set. Then we can define the operator $%
\tau _{j}$ $\in $ $B(H_{G})$ by

\strut

\begin{center}
$\tau _{j}(\xi _{w})$ $=$ $\xi _{w}$ $\xi _{e}$ $=$ $\xi _{w\,e},$ for all $%
\xi _{w}$ $\in $ $\mathcal{B}_{H_{G}},$
\end{center}

\strut

for each $j$ $\in $ $X,$ where $\mathcal{B}_{H_{G}}$ is the Hilbert basis of 
$H_{G},$ whenever an edge $e$ has its weight $x_{j},$ and $w$ $e$ $\neq $ $%
\emptyset $ in $\Bbb{G}.$ Then we can have the operator $\tau $ $\in $ $%
B(H_{G})$ defined by

\strut

\begin{center}
$\tau $ $=$ $\sum_{j=1}^{N}$ $\tau _{j}.$
\end{center}

\strut

This Hecke-type operator on $H_{G}$ is said to be the labeling operator of $%
\Bbb{G}$ on $H_{G}.$ We will consider the free distributional data of this
operator $\tau $ on $H_{G}.$ If $\Bbb{G}$ is a fractaloid, then $\tau $ is
self-adjoint. So, the free moments of it contain the spectral theoretical
properties of $\tau .$ In particular, we can show that the amalgamated free
moments $\left( E\left( \tau ^{n}\right) \right) _{n=1}^{\infty }$ of $\tau $
are determined by the cardinalities $\left( \eta _{n}\right) _{n=1}^{\infty
} $ of certain subsets $\mathcal{W}_{n}^{(m)}$ of the collection of all
finite paths in the shadowed graph $\widehat{G}$ of the given graph $G,$
determined recursively for $n$ $\in $ $\Bbb{N}$. Also, the equivalent free
distributional data is provided.

\strut

\strut

\strut

\section{Background}

\strut \strut \strut \strut \strut \strut \strut \strut \strut \strut \strut
\strut \strut

\strut Recently, the countable directed graphs have been studied in Pure and
Applied Mathematics, because not only that they are involved by a certain
noncommutative structures but also that they visualize such structures.
Futhermore, the visualization has nice matricial expressions, (sometimes,
the operator-valued matricial expressions dependent on) adjacency matrices
or incidence matrices of the given graph. In particular, partial isometries
in an operator algebra can be expressed and visualized by directed graphs:
in [10], [11], [13] and [14], we showed that each edge (resp. each vertex)
of a graph matches to a partial isometry (resp. a projection) in a (graph)
Hilbert space, and in [15], we showed that any partial isometry (resp.
projections induced by this partial isometry) on an arbitrary separable
infinite Hilbert space matches to an edge (resp. vertices) of certain
(corresponding) directed graph.

\strut

The main purpose of this paper is to introduce algebraic and
operator-algebraic structures induced by countable directed graphs.
Basically, we will follow the same settings used in [13]. For convenience,
as we assumed at the beginning of this paper, we will restrict our interests
to the case where the countable directed graphs are locally finite. In [10]
and [11], starting with a countable directed graph, we assign certain
algebraic elements gotten from the admissibility, and then we assign partial
isometries to those elements. From these operators, we generated a von
Neumann algebra and then we considered free probabilistic properties of
them. In conclusion, we found the nice amalgamated free block structures of
such von Neumann algebras and this provides the full characterization of the
von Neumann algebras generated by graph groupoids.

\strut \strut

\strut

\subsection{Amalgamated Free Probability}

\strut \strut \strut \strut

\strut

Let $B$ $\subset $ $A$ be von Neumann algebras with $1_{B}$ $=$ $1_{A}$ and
assume that there is a \emph{conditional expectation} $E_{B}$ $:$ $A$ $%
\rightarrow $ $B$ satisfying that (i) $E_{B}$ is a ($\Bbb{C}$-)linear map,
(ii) $E_{B}(b)$ $=$ $b,$ for all $b$ $\in $ $B$, (iii) $E_{B}(b_{1}$ $a$ $%
b_{2})$ $=$ $b_{1}$ $E_{B}(a)$ $b_{2},$ for all $b_{1},$ $b_{2}$ $\in $ $B$
and $a$ $\in $ $A,$ (iv) $E_{B}$ is continuous under the given topologies
for $A$ (and $B$), and (v) $E_{B}(a^{*})$ $=$ $E_{B}(a)^{*}$ in $B,$ for all 
$a$ $\in $ $A.$ The algebraic pair $(A,$ $E_{B})$ is said to be a $B$\emph{%
-valued }$W^{*}$\emph{-probability space}. Every operator in $(A,$ $E_{B})$
is called a $B$\emph{-valued (free) random variable}. Any $B$-valued random
variables have their $B$-valued free distributional data: $B$-valued $*$%
-moments and $B$-valued $*$-cumulants of them. Suppose $a_{1},$ ..., $a_{s}$
are $B$-valued random variables in $(A,$ $E_{B}),$ where $s$ $\in $ $\Bbb{N}%
. $ The $(i_{1},$ ..., $i_{n})$-th joint $B$-valued $*$-\emph{moments} of $%
a_{1},$ ..., $a_{s}$ are defined by

\strut

\begin{center}
$E_{B}\left(
(b_{1}a_{i_{1}}^{r_{i_{1}}})(b_{2}\,a_{i_{2}}^{r_{i_{2}}})\,...(b_{n}%
\,a_{i_{n}}^{r_{i_{n}}})\right) ,$
\end{center}

\strut

and the $(j_{1},$ ..., $j_{k})$-th joint $B$-valued $*$\emph{-cumulants} of $%
a_{1},$ ..., $a_{s}$ are defined by

\strut

\begin{center}
$k_{k}^{B}\left( (b_{1}a_{j_{1}}^{r_{i_{1}}}),\text{ ..., }%
(b_{k}a_{j_{k}}^{r_{j_{k}}})\right) =\underset{\pi \in NC(k)}{\sum }%
E_{B:\,\pi }\left( b_{1}a_{j_{1}}^{r_{j_{1}}}\text{ },...,\text{ }%
b_{k}a_{j_{k}}^{r_{j_{k}}}\right) \mu (\pi ,$ $1_{k}),$
\end{center}

\strut

for all $(i_{1},$ ..., $i_{n})$ $\in $ $\{1,$ ..., $s\}^{n}$ and for all $%
(j_{1},$ ..., $j_{k})$ $\in $ $\{1,$ ..., $s\}^{k},$ for $n,$ $k$ $\in $ $%
\Bbb{N},$ where $b_{j}$ $\in $ $B$ are arbitrary and $r_{i_{1}},$ ..., $%
r_{i_{n}},$ $r_{j_{1}},$ ..., $r_{j_{k}}$ $\in $ $\{1,$ $*\}$ and $NC(k)$ is
the lattice of all \emph{noncrossing partitions} with its minimal element $%
0_{k}$ $=$ $\{(1),$ $(2),$ ..., $(k)\}$ and its maximal element $1_{k}$ $=$ $%
\{(1,$ $2,$ ..., $k)\},$ for all $k$ $\in $ $\Bbb{N},$ and $\mu $ is the 
\emph{Moebius functional} in the incidence algebra $\mathcal{I}.$ Here, $%
E_{B:\pi }(...)$ is the \emph{partition-depending }$B$\emph{-valued moment}.
For example, if $\pi $ $=$ $\{(1,$ $4),$ $(2,$ $3),$ $(5)\}$ in $NC(5),$ then

\strut

\begin{center}
$E_{B:\pi }\left( a_{1},a_{2},a_{3},a_{4},a_{5}\right) =E_{B}\left(
a_{1}E_{B}(a_{2}a_{3})a_{4}\right) E_{B}(a_{5}).$
\end{center}

\strut \strut \strut \strut

Recall that the lattice $NC(n)$ of all noncrossing partitions over $\{1,$
..., $n\}$ has its partial ordering ``$\leq $'',

\strut

\begin{center}
$\pi \leq \theta $ $\overset{def}{\Longleftrightarrow }$ for each block $V$
in $\pi ,$ $\exists $ a block $B$ in $\theta $ s.t.,. $V$ $\subseteq $ $B,$
\end{center}

\strut

for $\pi ,\theta $ $\in $ $NC(n),$ where ``$\subseteq $'' means the usual
set inclusion, for all $n$ $\in $ $\Bbb{N}.$ Also recall that the \emph{%
incidence algebra} $\mathcal{I}$ is the collection of all functionals

\strut

\begin{center}
$\xi $ $:$ $\cup _{n=1}^{\infty }\left( NC(n)\times NC(n)\right) $ $%
\rightarrow $ $\Bbb{C}$
\end{center}

\strut

satisfying that $\xi (\pi ,$ $\theta )$ $=$ $0,$ whenever $\pi $ $>$ $\theta 
$, with the usual function addition $(+)$ and the \emph{convolution} $(*)$
defined by

\strut

\begin{center}
$\left( \xi _{1}*\xi _{2}\right) $ $(\pi ,$ $\theta )$ $\overset{def}{=}$ $%
\underset{\pi \leq \sigma \leq \theta }{\sum }$ $\xi _{1}(\pi ,\sigma )\xi
_{2}(\sigma ,\theta ),$
\end{center}

\strut

for all $\xi _{1},$ $\xi _{2}$ $\in $ $\mathcal{I}.$ If we define the\emph{\
zeta functional} $\zeta $ $\in $ $\mathcal{I}$ by

\strut

\begin{center}
$\zeta (\pi ,$ $\theta )$ $=$ $1,$ \ for all \ $\pi $ $\leq $ $\theta $ in $%
NC(n),$ for all $n$ $\in $ $\Bbb{N},$
\end{center}

\strut

then it is the identity element of $\mathcal{I}$ having its convolution
inverse, the Moebius functional $\mu $. Thus the Moebius functional $\mu $
satisfies that

\strut

\begin{center}
$\mu (0_{n},$ $1_{n})$ $=$ $(-1)^{n-1}$ $c_{n-1}$ \ and \ $\underset{\pi \in
NC(n)}{\sum }$ $\mu (\pi ,$ $1_{n})$ $=$ $0,$
\end{center}

\strut

where $c_{k}$ $\overset{def}{=}$ $\frac{1}{k+1}$ $\left( 
\begin{array}{l}
2k \\ 
\,\,k
\end{array}
\right) $ is the $k$-th Catalan number, for all $k$ $\in $ $\Bbb{N}$ (See
[21]).

\strut \strut

The $B$-valued freeness on $(A,$ $E_{B})$ is characterized by the $B$-valued 
$*$-cumulants. Let $A_{1}$ and $A_{2}$ be $W^{*}$-subalgebras of $A$ having
their common $W^{*}$-subalgebra $B.$ We say that $A_{1}$ and $A_{2}$ are 
\emph{free over} $B$ in $(A,$ $E_{B})$ if all mixed $B$-valued $*$-cumulants
of $A_{1}$ and $A_{2}$ vanish. The subsets $X_{1}$ and $X_{2}$ of $A$ are
said to be free over $B$ in $(A,$ $E_{B})$ if the $W^{*}$-subalgebras $%
vN(X_{1},$ $B)$ and $vN(X_{2},$ $B)$ of $A$ are free over $B$ in $(A,$ $%
E_{B}),$ where $vN(S_{1},$ $S_{2})$ means the von Neumann algebra generated
by sets $S_{1}$ and $S_{2}.$ Similarly, we say that the $B$-valued random
variables $x$ and $y$ are free over $B$ in $(A,$ $E_{B})$ if the subsets $%
\{x\}$ and $\{y\}$ are free over $B$ in $(A,$ $E_{B})$ (Also See [21]).

\strut

Suppose two $W^{*}$-subalgebras $A_{1}$ and $A_{2}$ of $A,$ containing their
common $W^{*}$-subalgebra $B,$ are free over $B$ in $(A,$ $E_{B}).$ Then we
can create the $W^{*}$-subalgebra $vN(A_{1},$ $A_{2})$ of $A$ generated by $%
A_{1}$ and $A_{2}.$ It is denoted by $A_{1}$ $*_{B}$ $A_{2},$ which is
called the $B$-valued free product (sub)algebra of $A_{1}$ and $A_{2}$ (over 
$B$). Assume that the given von Neumann algebra $A$ is generated by its $%
W^{*}$-subalgebras $A_{i}$, containing their common $W^{*}$-subalgebra $B,$
for $i$ $\in $ $\Lambda .$ Also, assume that all $A_{i}$'s are free over $B$
from each other in $(A,$ $E_{B}).$ Then $A$ is denoted by $\underset{i\in
\Lambda }{\,*_{B}}$ $A_{i}.$ i.e., the von Neumann algebra $A$ is said to be
a $B$-\emph{valued free product algebra} of $A_{i}$'s (over $B$).

\strut

Suppose a von Neumann algebra $A$ is a $B$-free product algebra $\underset{%
i\in \Lambda }{\,*_{B}}$ $A_{i}.$ Then $A$ is Banach-space isomorphic to the
Banach space

\strut

\begin{center}
$B$ $\oplus $ $\left( \oplus _{n=1}^{\infty }\left( \underset{i_{1}\neq
i_{2},\,i_{2}\neq i_{3},\,...,\,i_{n-1}\neq i_{n}}{\oplus }\left(
A_{i_{1}}^{o}\text{ }\otimes _{B}\text{ ... }\otimes _{B}\text{ }%
A_{i_{n}}^{o}\right) \right) \right) $
\end{center}

with

\begin{center}
$A_{i_{j}}^{o}$ $\overset{def}{=}$ $A_{i_{j}}$ $\ominus $ $B,$ for all $j$ $%
= $ $1,$ ..., $n,$
\end{center}

\strut

where $\otimes _{B}$ is the $B$-valued tensor product.\strut

\strut \strut \strut

\strut \strut

\subsection{Basic Concepts}

\strut \strut \strut

\strut

For a given countable directed graph $G,$ we can define the oppositely
directed graph $G^{-1},$ with $V(G^{-1})$ $=$ $V(G)$ and $E(G^{-1})$ $=$ $%
\{e^{-1}$ $:$ $e$ $\in $ $E(G)\},$ where each element $e^{-1}$ satisfies
that $e$ $=$ $v_{1}$ $e$ $v_{2}$ in $E(G)$, with $v_{1},$ $v_{2}$ $\in $ $%
V(G),$ if and only if $e^{-1}$ $=$ $v_{2}$ $e^{-1}$ $v_{1},$ in $E(G^{-1}).$
This opposite directed edge $e^{-1}$ $\in $ $E(G^{-1})$ of $e$ $\in $ $E(G)$
is called the shadow of $e.$ Also, this new graph $G^{-1}$, induced by $G,$
is said to be the shadow of $G.$ It is clear that $(G^{-1})^{-1}$ $=$ $G.$

\strut

Define the \emph{shadowed graph} $\widehat{G}$ $=$ $G$ $\cup $ $G^{-1}$ of $%
G $ by a directed graph with its vertex set $V(\widehat{G})$ $=$ $V(G)$ $=$ $%
V(G^{-1})$ and its edge set $E(\widehat{G})$ $=$ $E(G)$ $\cup $ $E(G^{-1})$,
where $G^{-1}$ is the \emph{shadow} of $G$. Then we can construct the \emph{%
free semigroupoid} $\Bbb{F}^{+}(\widehat{G})$ of the shadowed graph $%
\widehat{G},$ as a collection of all vertices and finite paths of $\widehat{G%
}$ with its binary operation called the \emph{admissibility}, where the
admissibility is nothing but the direction-depending connectedness of
elements in $\Bbb{F}^{+}(\widehat{G})$. Notice that all finite paths in $%
\Bbb{F}^{+}(\widehat{G})$ are the words in $V(\widehat{G})$ $\cup $ $E(%
\widehat{G}).$ By defining the \emph{reduction }(RR) on $\Bbb{F}^{+}(%
\widehat{G}),$ we can construct the graph groupoid $\Bbb{G}.$ i.e., the 
\emph{graph groupoid }$\Bbb{G}$ is a set of all ``reduced'' words in $E(%
\widehat{G}),$ with the inherited admissibility on $\Bbb{F}^{+}(\widehat{G}%
), $ where the reduction (RR) on $\Bbb{G}$ is

\strut

(RR)\qquad $\qquad \qquad w$ $w^{-1}$ $=$ $v$ and $w^{-1}w$ $=$ $v^{\prime
}, $

\strut

for all $w$ $=$ $v$ $w$ $v^{\prime }$ $\in $ $\Bbb{G},$ with $v,$ $v^{\prime
}$ $\in $ $V(\widehat{G}).$ In fact, this graph groupoid $\Bbb{G}$ is indeed
a categorial groupoid with its base $V(\widehat{G})$.

\strut \strut

For an arbitrary fixed von Neumann algebra $M$ in an operator algebra $B(K)$
of all bounded linear operators on a Hilbert space $K,$ we define a crossed
product algebra $\Bbb{M}_{G}$ $=$ $M$ $\times _{\beta }$ $\Bbb{G}$ of $M$
and $\Bbb{G}$ via a groupoid action $\beta $ $:$ $\Bbb{G}$ $\rightarrow $ $%
B(K$ $\otimes $ $H_{G}),$ where

$\strut $

\begin{center}
$H_{G}$ $\overset{def}{=}$ $\underset{w\in FP_{r}(\widehat{G})}{\oplus }$ $(%
\Bbb{C}$ $\cdot $ $\xi _{w})$
\end{center}

\strut

is the Hilbert space with its Hilbert basis $\{\xi _{w}$ $:$ $w$ $\in $ $%
FP_{r}(\widehat{G})\},$ where

\strut

\begin{center}
$FP_{r}(\widehat{G})$ $\overset{def}{=}$ $\Bbb{G}$ $\setminus $ $\left( V(%
\widehat{G})\text{ }\cup \text{ }\{\emptyset \}\right) .$
\end{center}

\strut

We will call $H_{G}$ the graph Hilbert space. A groupoid action $\beta ,$
called the graph-representation, is a nonunital intertwined partial
representation, determining the bounded operators $\beta _{w}$ on $K$ $%
\otimes $ $H_{G}$ satisfying that

\strut

\begin{center}
$\beta _{w}\left( m\right) $ $R_{w}$ $R_{w}^{*}=R_{w}^{*}$ $m$ $R_{w},$
\end{center}

and

\begin{center}
$\beta _{v}(m)$ $=$ $m,$
\end{center}

\strut \strut

for all \ $m$ $\in $ $M$,\ $w$ $\in $ $FP_{r}(\widehat{G})$ and $v$ $\in $ $%
V(\widehat{G}),$ where $R_{w}$ are the \emph{right} multiplication operators
on $H_{G}$ with their symbols $\xi _{w}$, for all $w$ $\in $ $\Bbb{G}.$ The
adjoint $R_{w}^{*}$ of $R_{w}$ are defined to be $R_{w^{-1}},$ for all $w$ $%
\in $ $\Bbb{G}.$

\strut

\begin{definition}
The crossed product algebra $\Bbb{M}_{G}$ $=$ $M$ $\times _{\beta }$ $\Bbb{G}
$ is the von Neumann algebra generated by $M$ and $\{R_{w}$ $:$ $w$ $\in $ $%
\Bbb{G}\},$ satisfying the above $\beta $-relations. We call $\Bbb{M}_{G}$ a
right graph von Neumann algebra induced by $G$ over $M$, via $\beta .$ A
right graph von Neumann algebra $\Bbb{M}_{G}$ has its canonical $W^{*}$%
-subalgebra $\Bbb{D}_{G}$ $\overset{def}{=}$ $\underset{v\in V(\widehat{G})}{%
\oplus }$ $(M$ $\cdot $ $R_{v}),$ called the $M$-diagonal subalgebra of $%
\Bbb{M}_{G}.$
\end{definition}

\strut

\begin{remark}
Recall that, in [10] and [11], we observed the (left) multiplication
operators $L_{w}$'s, for all $w$ $\in $ $\Bbb{G},$ instead of using right
multiplication operators $R_{w}$'s.
\end{remark}

\strut

Similar to [10], we construct an amalgamated $W^{*}$-probability space $(%
\Bbb{M}_{G},$ $E)$ over $\Bbb{D}_{G},$ where $E$ $:$ $\Bbb{M}_{G}$ $%
\rightarrow $ $\Bbb{D}_{G}$ is the canonical conditional expectation. Under
this setting, we realized that $\Bbb{M}_{G}$ is $*$-isomorphic to a $\Bbb{D}%
_{G}$-valued reduced free product $\underset{e\in E(G)}{*_{\Bbb{D}_{G}}^{r}}$
$\Bbb{M}_{e}$ of the $\Bbb{D}_{G}$-free blocks $\Bbb{M}_{e}$ indexed by
directed edges $e$ $\in $ $E(G),$ where

\strut

\begin{center}
$\Bbb{M}_{e}$ $\overset{def}{=}$ $vN\left( M\text{ }\times _{\beta }\text{ }%
\Bbb{G}_{e},\text{ }\Bbb{D}_{G}\right) ,$ for all $e$ $\in $ $E(G).$
\end{center}

\strut

Here, $\Bbb{G}_{e}$ is a substructure (or a subgroupoid) of $\Bbb{G},$
consisting of all reduced words in $\{e,$ $e^{-1}\},$ with the inherited
admissibility on $\Bbb{G},$ for all $e$ $\in $ $E(G).$ Notice that the
reduction of the free product ``$*_{\Bbb{D}_{G}}^{r}$'' is completely
dependent upon the admissibility on the graph groupoid $\Bbb{G}$.

\strut

Also, similar to [11], we can characterize the $\Bbb{D}_{G}$-free blocks $%
\Bbb{M}_{e}$'s of our right graph von Neumann algebra $\Bbb{M}_{G}.$ Because
of the setting, in fact, the results are almost same. Especially, we show
that each $\Bbb{D}_{G}$-free block $\Bbb{M}_{e}$ of $\Bbb{M}_{G}$ is $*$%
-isomorphic to a certain von Neumann algebra contained in either $\mathcal{M}%
_{lp}^{\beta }$ or $\mathcal{M}_{non-lp}^{\beta },$ where

$\strut $

\begin{center}
$\mathcal{M}_{lp}^{\beta }$ $=$ $\{vN(M$ $\times _{\lambda _{e}}$ $\Bbb{Z},$ 
$\Bbb{D}_{G})$ $:$ $\lambda _{e}$ $=$ $\beta \mid _{\Bbb{S}_{e}},$ $e$ is a
loop edge$\}$
\end{center}

and

\begin{center}
$\mathcal{M}_{non-lp}^{\beta }$ $=$ $\{vN(M_{2}^{\beta _{e}}(M),$ $\Bbb{D}%
_{G})$ $:$ $e$ is a non-loop edge$\}$,
\end{center}

\strut

where $M$ $\times _{\lambda ^{(e)}}$ $\Bbb{Z}$ is a classical group crossed
products of $M$ and the infinite cyclic abelian group $\Bbb{Z}$ via a group
action $\lambda ^{(e)}$ satisfying that $\lambda _{(e)}$ $=$ $\beta $ $\mid
_{\Bbb{G}_{e}},$ and where $M_{2}^{\beta _{e}}(M)$ is a $W^{*}$-subalgebra
of $M_{2}(M)$ $=$ $M$ $\otimes _{\Bbb{C}}$ $M_{2}(\Bbb{C})$ satisfying the $%
\beta _{e}$-relation: $\beta _{e}(m)$ $R_{e}$ $R_{e}^{*}$ $=$ $R_{e}^{*}$ $m$
$R_{e}$, for all $m$ $\in $ $M$ and for $e$ $\in $ $E(G),$ where $M_{2}(\Bbb{%
C})$ is the matricial algebra generated by all $(2$ $\times $ $2)$-matrices.
In particular, if $M$ $=$ $\Bbb{C},$ we can conclude that each $\Bbb{D}_{G}$%
-free block $\Bbb{M}_{e}$ of $\Bbb{M}_{G}$ is $*$-isomorphic to either $%
L^{\infty }(\Bbb{T})$ or $M_{2}(\Bbb{C}),$ where $\Bbb{T}$ is the unit
circle in $\Bbb{C}$, for all $e$ $\in $ $E(G).$ This characterization says
that the study of right graph von Neumann algebras is the investigation of
graph groupoids and the above two types of von Neumann algebras.\strut

\strut

\textbf{Assumption} In this paper, we only consider the case where the fixed
von Neumann algebra $M$ is $\Bbb{C}.$ Then every right graph von Neumann
algebras $\Bbb{M}_{G}$ $=$ $\Bbb{C}$ $\times _{\beta }$ $\Bbb{G}$ are all $*$%
-isomorphic to $M_{G}$ $=$ $\overline{\Bbb{C}[\Bbb{G}]}^{w},$ as $W^{*}$%
-subalgebras of $B(H_{G}),$ for any choice of $\beta ,$ by the linearity of $%
\beta $ on $\Bbb{C}$. We will say that the von Neumann algebra $M_{G}$ is 
\emph{the} right graph von Neumann algebra of $G.$ $\square $

\strut \strut

We will consider the labeling operator $\tau $ of $\Bbb{G}$, which is an
analogue of classical Heck-type operators (of groups). As an element of the
right graph von Neumann algebra $M_{G}$, we can verify that $\tau $ has its
decomposition $\tau $ $=$ $\oplus _{j=1}^{N}$ $\tau _{j},$ where

\strut

\begin{center}
$N$ $\overset{def}{=}$ $\max \{\deg _{out}(v)$ $:$ $v$ $\in $ $V(\widehat{G}%
)\}$
\end{center}

\strut and

\begin{center}
$\tau _{j}$ $(\xi _{w})$ $\overset{def}{=}$ $\left\{ 
\begin{array}{ll}
\xi _{w}\text{ }\xi _{e}=\xi _{w\,e} & 
\begin{array}{l}
\,\,\text{if }w\,e\neq \emptyset \text{ and} \\ 
e\text{ has its weight }j
\end{array}
\\ 
&  \\ 
\xi _{w}\text{ }\xi _{\emptyset }=\xi _{\emptyset }=0 & 
\begin{array}{l}
\text{if there is no edge }e\text{ such that} \\ 
e\text{ has its weight }j\text{ and }w\,e\neq \emptyset .
\end{array}
\end{array}
\right. $
\end{center}

\strut

The study of the labeling operator $\tau ,$ itself, is interesting in
Operator Theory and Quantum Physics, but we will concentrate on observing
its free distributional data of the operator $\tau $. The data will show how
the graph groupoid $\Bbb{G}$ acts inside $B(H_{G}).$ Futhermore, if we
consider the general case where we have $B(K$ $\otimes $ $H_{G}),$ for any
Hilbert spaces $K,$ we could show how $\Bbb{G}$ can acts on $K$ $\otimes $ $%
H_{G},$ too.

\strut

\strut

\strut

\section{Definitions}

\strut

\strut

In this chapter, we will introduce the definitions and concepts which we
will use in the rest of this paper. We will review groupoids, groupoid
actions, automata. Also, we observe the basic properties of right graph von
Neumann algebras.\strut

\strut

\strut

\subsection{Categorial Groupoids and Groupoid Actions\strut}

\strut

\strut

We say an algebraic structure $(\mathcal{X},$ $\mathcal{Y},$ $s,$ $r)$ is a 
\emph{(categorial) groupoid} if it satisfies that (i) $\mathcal{Y}$ $\subset 
$ $\mathcal{X},$ (ii) for all $x_{1},$ $x_{2}$ $\in $ $\mathcal{X},$ there
exists a partially-defined binary operation $(x_{1},$ $x_{2})$ $\mapsto $ $%
x_{1}$ $x_{2},$ for all $x_{1},$ $x_{2}$ $\in $ $\mathcal{X},$ depending on
the source map $s$ and the range map $r$ satisfying that:

\strut

(ii-1) $x_{1}$ $x_{2}$ is well-determined, whenever $r(x_{1})$ $=$ $s(x_{2})$
and in this case, $s(x_{1}$ $x_{2})$ $=$ $s(x_{1})$ and $r(x_{1}$ $x_{2})$ $%
= $ $r(x_{2}),$ for $x_{1},$ $x_{2}$ $\in $ $\mathcal{X},$

\strut

(ii-2) $(x_{1}$ $x_{2})$ $x_{3}$ $=$ $x_{1}$ $(x_{2}$ $x_{3})$, if they are
well-determined in the sense of (ii-1), for $x_{1},$ $x_{2},$ $x_{3}$ $\in $ 
$\mathcal{X},$

\strut

(ii-3) if $x$ $\in $ $\mathcal{X},$ then there exist $y,$ $y^{\prime }$ $\in 
$ $\mathcal{Y}$ such that $s(x)$ $=$ $y$ and $r(x)$ $=$ $y^{\prime },$
satisfying $x$ $=$ $y$ $x$ $y^{\prime }$ (Here, the elements $y$ and $%
y^{\prime }$ are not necessarily distinct),

\strut \strut

(ii-4) if $x$ $\in $ $\mathcal{X},$ then there exists a unique element $%
x^{-1}$ for $x$ satisfying $x$ $x^{-1}$ $=$ $s(x)$ and $x^{-1}$ $x$ $=$ $%
r(x).$

\strut

Thus, every group is a groupoid $(\mathcal{X},$ $\mathcal{Y},$ $s,$ $r)$
with $\left| \mathcal{Y}\right| $ $=$ $1$ (and hence $s$ $=$ $r$ on $%
\mathcal{X}$). This subset $\mathcal{Y}$ of $\mathcal{X}$ is said to be the 
\emph{base of} $\mathcal{X}$. Remark that we can naturally assume that there
exists the \emph{empty element} $\emptyset $ in a groupoid $\mathcal{X}.$
The empty element $\emptyset $ means the products $x_{1}$ $x_{2}$ are not
well-defined, for some $x_{1},$ $x_{2}$ $\in $ $\mathcal{X}.$ Notice that if 
$\left| \mathcal{Y}\right| $ $=$ $1$ (equivalently, if $\mathcal{X}$ is a
group), then the empty word $\emptyset $ is not contained in the groupoid $%
\mathcal{X}.$ However, in general, whenever $\left| \mathcal{Y}\right| $ $%
\geq $ $2,$ a groupoid $\mathcal{X}$ always contain the empty word. So, if
there is no confusion, we will naturally assume that the empty element $%
\emptyset $ is contained in $\mathcal{X}.$

\strut

It is easily checked that our graph groupoid $\Bbb{G}$ of a countable
directed graph $G$ is indeed a groupoid with its base $V(\widehat{G}).$
i.e., every graph groupoid $\Bbb{G}$ of a countable directed graph $G$ is a
groupoid $(\Bbb{G},$ $V(\widehat{G}),$ $s$, $r)$, where $s(w)$ $=$ $s(v$ $w)$
$=$ $v$ and $r(w)$ $=$ $r(w$ $v^{\prime })$ $=$ $v^{\prime },$ for all $w$ $%
= $ $v$ $w$ $v^{\prime }$ $\in $ $\Bbb{G}$ with $v,$ $v^{\prime }$ $\in $ $V(%
\widehat{G}).$ i.e., the vertex set $V(\widehat{G})$ $=$ $V(G)$ is a base of 
$\Bbb{G}.$

\strut

Let $\mathcal{X}_{k}$ $=$ $(\mathcal{X}_{k},$ $\mathcal{Y}_{k},$ $s_{k},$ $%
r_{k})$ be groupoids, for $k$ $=$ $1,$ $2.$ We say that a map $f$ $:$ $%
\mathcal{X}_{1}$ $\rightarrow $ $\mathcal{X}_{2}$ is a \emph{groupoid
morphism} if (i) $f$ is a function, (ii) $f(\mathcal{Y}_{1})$ $\subseteq $ $%
\mathcal{Y}_{2},$ (iii) $s_{2}\left( f(x)\right) $ $=$ $f\left(
s_{1}(x)\right) $ in $\mathcal{X}_{2},$ for all $x$ $\in $ $\mathcal{X}_{1}$%
, and (iv) $r_{2}\left( f(x)\right) $ $=$ $f\left( r_{1}(x)\right) $ in $%
\mathcal{X}_{2},$ for all $x$ $\in $ $\mathcal{X}_{1}.$ If a groupoid
morphism $f$ is bijective, then we say that $f$ is a \emph{%
groupoid-isomorphism}, and the groupoids $\mathcal{X}_{1}$ and $\mathcal{X}%
_{2}$ are said to be \emph{groupoid-isomorphic}.

\strut

Notice that, if two countable directed graphs $G_{1}$ and $G_{2}$ are
graph-isomorphic, via a graph-isomorphism $g$ $:$ $G_{1}$ $\rightarrow $ $%
G_{2},$ in the sense that (i) $g$ is bijective from $V(G_{1})$ onto $%
V(G_{2}),$ (ii) $g$ is bijective from $E(G_{1})$ onto $E(G_{2}),$ (iii) $%
g(e) $ $=$ $g(v_{1}$ $e$ $v_{2})$ $=$ $g(v_{1})$ $g(e)$ $g(v_{2})$ in $%
E(G_{2}),$ for all $e$ $=$ $v_{1}$ $e$ $v_{2}$ $\in $ $E(G_{1}),$ with $%
v_{1},$ $v_{2}$ $\in $ $V(G_{1}),$ then the graph groupoids $\Bbb{G}_{1}$
and $\Bbb{G}_{2}$ are groupoid-isomorphic. More generally, if two graphs $%
G_{1}$ and $G_{2}$ have graph-isomorphic shadowed graphs $\widehat{G_{1}}$
and $\widehat{G_{2}}, $ then $\Bbb{G}_{1}$ and $\Bbb{G}_{2}$ are
groupoid-isomorphic.

\strut \strut \strut

Let $\mathcal{X}$ $=$ $(\mathcal{X},$ $\mathcal{Y},$ $s,$ $r)$ be a
groupoid. We say that this groupoid $\mathcal{X}$ \emph{acts on a set }$Y$
if there exists a groupoid action $\pi $ of $\mathcal{X}$ such that $\pi (x)$
$:$ $Y$ $\rightarrow $ $Y$ is a well-determined function, for all $x$ $\in $ 
$\mathcal{X}.$ Sometimes, we call the set $Y,$ a $\mathcal{X}$\emph{-set}.
We are interested in the case where a $\mathcal{X}$-set $Y$ is a Hilbert
space. The nicest example of a groupoid action acting on a Hilbert space is
a graph-representation in the sense of [10] and [11]. (Also, see [14] and
[15] and the action $\beta $ in Section 2.5.)

\strut

Let $\mathcal{X}_{1}$ $\subset $ $\mathcal{X}_{2}$ be a subset, where $%
\mathcal{X}_{2}$ $=$ $(\mathcal{X}_{2},$ $\mathcal{Y}_{2},$ $s,$ $r)$ is a
groupoid, and assume that $\mathcal{X}_{1}$ $=$ $(\mathcal{X}_{1},$ $%
\mathcal{Y}_{1},$ $s,$ $r),$ itself, is a groupoid, where $\mathcal{Y}_{1}$ $%
=$ $\mathcal{X}_{2}$ $\cap $ $\mathcal{Y}_{2}.$ Then we say that the
groupoid $\mathcal{X}_{1}$ is a \emph{subgroupoid} of $\mathcal{X}_{2}.$

\strut \strut

Recall that we say that a countable directed graph $G_{1}$ is a \emph{%
full-subgraph} of a countable directed graph $G_{2},$ if

\strut

\begin{center}
$E(G_{1})$ $\subseteq $ $E(G_{2})$
\end{center}

and

\begin{center}
$V(G_{1})$ $=$ $\{v$ $\in $ $V(G_{1})$ $:$ $e$ $=$ $v$ $e$ or $e$ $=$ $e$ $%
v, $ $\forall $ $e$ $\in $ $E(G_{1})\}.$
\end{center}

\strut

Remark the difference between full-subgraphs and subgraphs: We say that $%
G_{1}^{\prime }$ is a \emph{subgraph} of $G_{2},$ if

\strut

\begin{center}
$V(G_{1}^{\prime })$ $\subseteq $ $V(G_{2})$
\end{center}

and

\begin{center}
$E(G_{1}^{\prime })$ $=$ $\{e$ $\in $ $E(G_{2})$ $:$ $e$ $=$ $v_{1}$ $e$ $%
v_{2},$ for $v_{1},$ $v_{2}$ $\in $ $V(G_{1}^{\prime })\}.$
\end{center}

\strut

We can see that the graph groupoid $\Bbb{G}_{1}$ of $G_{1}$ is a subgroupoid
of the graph groupoid $\Bbb{G}_{2}$ of $G_{2},$ whenever $G_{1}$ is a
full-subgraph of $G_{2}.$

\strut \strut

\strut

\subsection{Automata and Fractal Groups}

\strut \strut

\strut \strut

Automata Theory is the study of abstract machines, and we are using it in
the formulation given by von Neumann. It is related to the theory of formal
languages. In fact, automata may be thought of as the class of formal
languages they are able to recognize. In von Neumann's version, an automaton
is a finite state machine (FSM). i.e., a machine with input of symbols,
transitions through a series of states according to a transition function
(often expressed as a table). The transition function tells the automata
which state to go to next, given a current state and a current symbol. The
input is read sequentially, symbol by symbol, for example as a tape with a
word written on it, registered by the head of the automaton; the head moves
forward over the tape one symbol at a time. Once the input is depleted, the
automaton stops. Depending on the state in which the automaton stops, it is
said that the automaton either accepts or rejects the input. The set of all
the words accepted by the automaton is called the language of the automaton.
For the benefit for the readers, we offer the following references for the
relevant part of Automata Theory: [1], [35], [36], [37], [50] and [51].

\strut

Let the quadruple $\mathcal{A}$ $=$ $<D,$ $Q,$ $\varphi ,$ $\psi >$ be
given, where $D$ and $Q$ are sets and

\strut

\begin{center}
$\varphi $ $:$ $D$ $\times $ $Q$ $\rightarrow $ $Q$ \ \ \ and \ \ $\psi $ $:$
$D$ $\times $ $Q$ $\rightarrow $ $D$
\end{center}

\strut

are maps. We say that $D$ and $Q$ are the (finite) \emph{alphabet} and the 
\emph{state set} of $\mathcal{A},$ respectively and we say that $\varphi $
and $\psi $ are the \emph{output function} and the s\emph{tate transition
function}, respectively. In this case, the quadruple $\mathcal{A}$ is called
an \emph{automaton}. If the map $\psi (\bullet ,$ $q)$ is bijective on $D,$
for any fixed $q$ $\in $ $Q,$ then we say that the automaton $A$ is \emph{%
invertible}. Similarly, if the map $\varphi (x,$ $\bullet )$ is bijective on 
$Q,$ for any fixed $x$ $\in $ $D,$ then we say that the automaton $\mathcal{A%
}$ is \emph{reversible}. If the automaton $\mathcal{A}$ is both invertible
and reversible, then $\mathcal{A}$ is said to be \emph{bi-reversible}.

\strut

To help visualize the use of automata, a few concrete examples may help.
With some oversimplification, they may be drawn from the analysis and
synthesis of input / output models in Engineering, often referred to as
black box diagram: excitation variables, response variables, and
intermediate variables (e.g., see [54] and [55]).

\strut

In our presentation above, the $D$ (the chosen finite alphabet) often takes
different forms on the side of input $D_{i}$ and output $D_{o}.$ In popular
automata that models stimuli of organisms, the three sets input $D_{i},$
output $D_{o}$, and the state set $Q,$ could be as in the following
prototypical three examples:

\strut

\begin{example}
Models stimuli of organisms:

\strut 

$\qquad D_{i}$ $=$ $\{$positive stimulus, negative stimulus$\},$

$\qquad D_{o}$ $=$ $\{$reaction, no reaction$\},$

and

$\qquad Q$ $=$ $\left\{ 
\begin{array}{c}
\text{reaction to last positive stimulus,} \\ 
\text{no reaction to last positive stimulus}
\end{array}
\right\} .$
\end{example}

\strut

\begin{example}
In a control model for say a steering mechanism in a vehicle:

\strut 

$\qquad D_{i}$ $=$ $\{$right, left$\},$

$\qquad D_{o}$ $=$ $\{$switch on, switch off$\}$ or $\{$lamp on, lamp off$\}$

and

$\qquad Q$ $=$ $\left\{ 
\begin{array}{c}
\text{right-turning direction signal} \\ 
\text{left-turning direction signal}
\end{array}
\right\} .$
\end{example}

\strut

\begin{example}
In a model for quantization in Signal Processing:

\strut 

\begin{center}
$D$ $=$ $D_{i}$ $=$ $D_{o}$ $=$ $\left\{ 
\begin{array}{c}
\text{assignments from a bit alphabet,} \\ 
\text{with the bits referring to the value} \\ 
\text{of pulses-in and pulses-out in a} \\ 
\text{signal processing algorithm}
\end{array}
\right\} ,$
\end{center}

\strut 

for example, on a discrete multiresolution (e.g., [50]), and

\strut 

\begin{center}
$Q$ $=$ $
\begin{array}{l}
\text{a subset of the Cartesian product of} \\ 
\text{copies of }D,\text{ fixing finite number of} \\ 
\text{times, i.e., }D\text{ }\times \text{ ... }\times \text{ }D
\end{array}
$
\end{center}
\end{example}

\strut \strut

Recently, various algebraists have studied automata and the corresponding
automata groups (Also, see [1], [22], [35] and [37]). We will consider a
certain special case, where $Q$ is a free semigroupoid of a shadowed graph.

\strut

Roughly speaking, a \emph{undirected tree} is a connected simplicial graph
without loop finite paths. Recall that a (undirected) graph is \emph{%
simplicial}, if the graph has neither loop-edges nor multi-edges. A\emph{\
directed tree} is a connected simplicial graph without loop finite paths,
with directed edges. In particular, we say that a directed tree $\mathcal{T}%
_{n}$ is a $n$-\emph{regular tree}, if $\mathcal{T}_{n}$ is rooted and
one-flowed infinite directed tree, having the same out-degrees $n,$ for all
vertices. For example, the $2$-regular tree $\mathcal{T}_{2}$ can be
depicted by

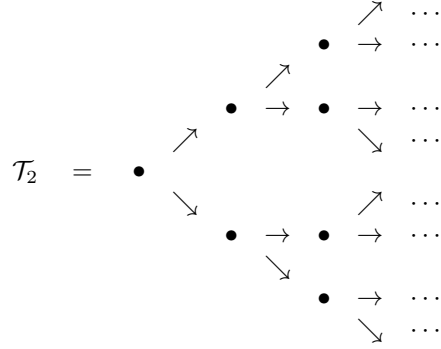
\strut

\begin{center}
$\mathcal{T}_{2}\quad =$\quad $
\begin{array}{lllllll}
&  &  &  &  & \nearrow & \cdots \\ 
&  &  &  & \bullet & \rightarrow & \cdots \\ 
&  &  & \nearrow &  &  &  \\ 
&  & \bullet & \rightarrow & \bullet & \rightarrow & \cdots \\ 
& \nearrow &  &  &  & \searrow & \cdots \\ 
\bullet &  &  &  &  &  &  \\ 
& \searrow &  &  &  & \nearrow & \cdots \\ 
&  & \bullet & \rightarrow & \bullet & \rightarrow & \cdots \\ 
&  &  & \searrow &  &  &  \\ 
&  &  &  & \bullet & \rightarrow & \cdots \\ 
&  &  &  &  & \searrow & \cdots
\end{array}
$
\end{center}

\strut

Let $\mathcal{A}$ $=$ $<D,$ $Q,$ $\varphi ,$ $\psi >$ be an automaton with $%
\left| D\right| $ $=$ $n.$ Then, we can construct \emph{automata action}s of 
$\mathcal{A}$ on $\mathcal{T}_{n}.$ Let's fix $q$ $\in $ $Q.$ Then the
action of $\mathcal{A}_{q}$ is defined on the finite words $D_{*}$ of $D$ by

\strut

\begin{center}
$\mathcal{A}_{q}\left( x\right) $ $\overset{def}{=}$ $\varphi (x,$ $q),$ for
all $x$ $\in $ $D,$
\end{center}

\strut

and recursively,

$\strut $

\begin{center}
$\mathcal{A}_{q}\left( (x_{1},\text{ }x_{2},\text{ ..., }x_{m})\right) $ $=$ 
$\varphi \left( x_{1},\text{ }\mathcal{A}_{q}(x_{2},...,x_{m})\right) ,$
\end{center}

\strut

for all $(x_{1},$ ..., $x_{m})$ $\in $ $D_{*},$ where

\strut

\begin{center}
$D_{*}$ $\overset{def}{=}$ $\cup _{m=1}^{\infty }$ $\left( \left\{ (x_{1},%
\text{ ..., }x_{m})\in D^{m}\left| 
\begin{array}{c}
\text{ }x_{k}\in D,\text{ for all} \\ 
k=1,...,m
\end{array}
\right. \right\} \right) .$
\end{center}

\strut

Then the automata actions $\mathcal{A}_{q}$'s are acting on the $n$-regular
tree $\mathcal{T}_{n}$. In other words, all images of automata actions are
regarded as an elements in the free semigroupoid $\Bbb{F}^{+}(\mathcal{T}%
_{n})$ of the $n$-regular tree. i.e.,

\strut

\begin{center}
$V(\mathcal{T}_{n})$ $\supseteq $ $D_{*}$
\end{center}

\strut

and its edge set

\strut

\begin{center}
$\strut E(\mathcal{T}_{n})$ $\supseteq $ $\{\mathcal{A}_{q}(x)$ $:$ $x$ $\in 
$ $D,$ $q$ $\in $ $Q\}.$
\end{center}

\strut

This makes us to illustrate how the automata actions work.

\strut

Let $\mathcal{C}$ $=$ $\{\mathcal{A}_{q}$ $:$ $q$ $\in $ $Q\}$ be the
collection of automata actions of the given automaton $\mathcal{A}$ $=$ $<D,$
$Q,$ $\varphi ,$ $\psi >$. Then we can create a group $G(\mathcal{A})$
generated by the collection $\mathcal{C}.$ This group $G(\mathcal{A})$ is
called the \emph{automata group generated by} $\mathcal{A}.$ The generator
set $\mathcal{C}$ of $G(\mathcal{A})$ acts ``fully'' on the $\left| D\right| 
$-regular tree $\mathcal{T}_{\left| D\right| },$ we say that this group $G(%
\mathcal{A})$ is a \emph{fractal group}. There are many ways to define
fractal groups, but we define them in the sense of automata groups. (See [1]
and [37]. In fact, in [37], Batholdi, Grigorchuk and Nekrashevych did not
define the term ``fractal'', but they provide the fractal properties.)

\strut

Now, we will define a fractal group more precisely (Also see [1]). Let $%
\mathcal{A}$ be an automaton and let $\Gamma $ $=$ $G(\mathcal{A})$ be the
automata group generated by the automata actions acting on the $n$-regular
tree $\mathcal{T}_{n},$ where $n$ is the cardinality of the alphabet of $%
\mathcal{A}.$ By $St_{\Gamma }(k),$ denote the subgroup of $\Gamma $ $=$ $G(%
\mathcal{A})$, consisting of those elements of $\Gamma ,$ acting trivially
on the $k$-th level of $\mathcal{T}_{n},$ for all $k$ $\in $ $\Bbb{N}$ $\cup 
$ $\{0\}.$

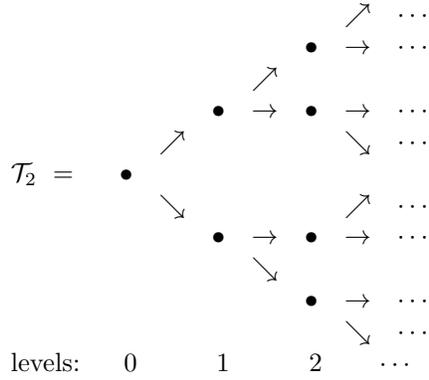
\strut

\begin{center}
$
\begin{array}{ll}
\mathcal{T}_{2}\text{ }= & 
\begin{array}{lllllll}
&  &  &  &  & \nearrow & \cdots \\ 
&  &  &  & \bullet & \rightarrow & \cdots \\ 
&  &  & \nearrow &  &  &  \\ 
&  & \bullet & \rightarrow & \bullet & \rightarrow & \cdots \\ 
& \nearrow &  &  &  & \searrow & \cdots \\ 
\bullet &  &  &  &  &  &  \\ 
& \searrow &  &  &  & \nearrow & \cdots \\ 
&  & \bullet & \rightarrow & \bullet & \rightarrow & \cdots \\ 
&  &  & \searrow &  &  &  \\ 
&  &  &  & \bullet & \rightarrow & \cdots \\ 
&  &  &  &  & \searrow & \cdots
\end{array}
\\ 
\text{levels:} & \,\,\,\,0\qquad \quad 1\qquad \quad 2\qquad \cdots
\end{array}
$
\end{center}

\strut

Analogously, for a vertex $u$ in $\mathcal{T}_{n},$ define $St_{\Gamma }(u)$
by the subgroup of $\Gamma ,$ consisting of those elements of $\Gamma ,$
acting trivially on $u.$ Then

$\strut $

\begin{center}
$St_{\Gamma }(k)$ $=$ $\underset{u\,:\,\text{vertices of the }k\text{-th
level of }\mathcal{T}_{n}}{\cap }$ $\left( St_{\Gamma }(u)\right) .$
\end{center}

\strut

For any vertex $u$ of $\mathcal{T}_{n},$ we can define the algebraic
projection $p_{u}$ $:$ $St_{\Gamma }(u)$ $\rightarrow $ $\Gamma .$

\strut

\begin{definition}
Let $\Gamma $ $=$ $G(\mathcal{A})$ be the automata group given as above. We
say that this group $\Gamma $ is a fractal group if, for any vertex $u$ of $%
\mathcal{T}_{n},$ the image of the projection $p_{u}\left( St_{\Gamma
}(u)\right) $ is group-isomorphic to $\Gamma ,$ after the identification of
the tree $\mathcal{T}_{n}$ with its subtree $\mathcal{T}_{u}$ with the root $%
u.$
\end{definition}

\strut

For instance, if $u$ is a vertex of the $2$-regular tree $\mathcal{T}_{2}$,
then we can construct a subtree $\mathcal{T}_{u},$ as follows:

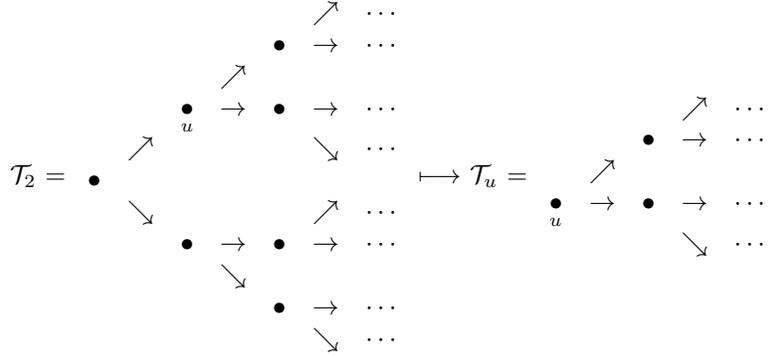
\strut

\begin{center}
$\mathcal{T}_{2}$ $=$ $
\begin{array}{lllllll}
&  &  &  &  & \nearrow & \cdots \\ 
&  &  &  & \bullet & \rightarrow & \cdots \\ 
&  &  & \nearrow &  &  &  \\ 
&  & \underset{u}{\bullet } & \rightarrow & \bullet & \rightarrow & \cdots
\\ 
& \nearrow &  &  &  & \searrow & \cdots \\ 
\bullet &  &  &  &  &  &  \\ 
& \searrow &  &  &  & \nearrow & \cdots \\ 
&  & \bullet & \rightarrow & \bullet & \rightarrow & \cdots \\ 
&  &  & \searrow &  &  &  \\ 
&  &  &  & \bullet & \rightarrow & \cdots \\ 
&  &  &  &  & \searrow & \cdots
\end{array}
$ $\longmapsto $ $\mathcal{T}_{u}$ $=$ $
\begin{array}{lllll}
&  &  & \nearrow & \cdots \\ 
&  & \bullet & \rightarrow & \cdots \\ 
& \nearrow &  &  &  \\ 
\underset{u}{\bullet } & \rightarrow & \bullet & \rightarrow & \cdots \\ 
&  &  & \searrow & \cdots
\end{array}
.$
\end{center}

\strut

As we can check, the graphs $\mathcal{T}_{2}$ and $\mathcal{T}_{u}$ are
graph-isomorphic. The above definition shows that if the automata actions $%
\mathcal{A}_{q}$'s of $\mathcal{A}$ are acting \emph{fully} on $\mathcal{T}%
_{n},$ then the automata group $G(\mathcal{A})$ is a fractal group. There
are lots of famous fractal groups, but we introduce the following example,
for our purpose.

\strut

\begin{example}
Let $\mathcal{A}$ $=$ $<X_{2n},$ $F_{n},$ $\varphi ,$ $\psi >$ be an
automaton, where $F_{n}$ is the free group with its generator set $X_{2n}$ $=
$ $\{g_{1}^{\pm 1},$ ..., $g_{n}^{\pm 1}\}.$ Then the automata group $G(%
\mathcal{A})$ is group-isomorphic to $F_{n}.$ It is easy to check that all
elements in $F_{n}$ acts fully on the $2n$-regular tree $\mathcal{T}_{2n},$
and hence $G(\mathcal{A})$ is a fractal group. For example, if $n$ $=$ $2,$
then we can get the following $0$-th and $1$-st levels of $\mathcal{T}_{4}$:

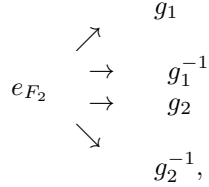
\strut 

\begin{center}
$
\begin{array}{lll}
&  & g_{1} \\ 
& \nearrow  &  \\ 
e_{F_{2}} & 
\begin{array}{l}
\rightarrow  \\ 
\rightarrow 
\end{array}
& 
\begin{array}{l}
g_{1}^{-1} \\ 
g_{2}
\end{array}
\\ 
& \searrow  &  \\ 
&  & g_{2}^{-1},
\end{array}
$
\end{center}

\strut 

where $e_{F_{2}}$ is the group-identity of $F_{2}$ $=$ $<g_{1},$ $g_{2}>.$%
\strut \strut \strut \strut \strut \strut \strut 
\end{example}

\strut \strut \strut \strut

\strut \strut \strut \strut

\subsection{Right Graph Von Neumann algebras}

\strut

\strut

\strut \strut \strut \strut

Let $G$ be a countable directed graph with its graph groupoid $\Bbb{G}.$ We
define and consider certain operators on the graph Hilbert space $H_{G}$
induced by $G$. As we defined and observed, we can regard all elements in $%
\Bbb{G}$ as reduced words in $E(\widehat{G}),$ under the admissibility with
the reduction (RR).

\strut

\begin{definition}
Let $G$ be a countable directed graph and let $\Bbb{G}$ be the corresponding
graph groupoid. Define the Hilbert space $H_{G}$ of $G$ by

\strut 

\begin{center}
$H_{G}$ $\overset{def}{=}$ $\ \underset{w\in FP_{r}(\widehat{G})}{\oplus }%
\left( \Bbb{C}\xi _{w}\right) $
\end{center}

\strut 

with its Hilbert basis $\mathcal{B}_{H_{G}}$ $=$ $\{\xi _{w}$ $:$ $w$ $\in $ 
$FP_{r}(\widehat{G})\}$ in $H_{G},$ where $FP_{r}(\widehat{G})$ is the
reduced finite path set.
\end{definition}

\strut \strut \strut

We have the following multiplication rule on $H_{G}$:

\strut

\begin{center}
$\xi _{w_{1}}\xi _{w_{2}}=\left\{ 
\begin{array}{ll}
\xi _{w_{1}w_{2}} & \text{if }w_{1}\text{ }w_{2}\text{ }\neq \text{ }%
\emptyset \\ 
\xi _{\emptyset }=0_{H_{G}} & \text{otherwise,}
\end{array}
\right. $
\end{center}

\strut

for all $\xi _{w_{1}},$ $\xi _{w_{2}}$ $\in $ $\mathcal{B}_{H_{G}}.$ Suppose 
$w_{1}$ $=$ $w$ and $w_{2}$ $=$ $w^{-1}$ in $FP_{r}(\widehat{G}).$ Then, by
the above multiplication rule, we can have an element $\xi _{w}$ $\xi
_{w^{-1}}$ $=$ $\xi _{ww^{-1}}$ in $H_{G},$ where $w$ $w^{-1}$ is a vertex
in $V(\widehat{G}).$ So, we can determine $\xi _{v}$ $\in $ $H_{G},$ for all 
$v$ $\in $ $V(\widehat{G}),$ too. Thus, for any $w$ $\in $ $\Bbb{G},$ we can
have the corresponding Hilbert space element $\xi _{w}$ in $H_{G}.$ This
multiplication rule let us define left and right multiplication operators on 
$H_{G}.$ The \emph{(left) multiplication operators} $L_{w}$'s on $H_{G}$ are
defined by

\strut

\begin{center}
$L_{w}$ $\xi _{w^{\prime }}$ $\overset{def}{=}$ $\xi _{w}$ $\xi _{w^{\prime
}}$ $=$ $\xi _{ww^{\prime }},$ for all $w,$ $w^{\prime }$ $\in $ $\Bbb{G},$
\end{center}

\strut

where $\xi _{\emptyset }$ $\overset{def}{=}$ $0_{H_{G}}.$ It is easy to see
that $L_{w_{1}}$ $L_{w_{2}}$ $=$ $L_{w_{1}w_{2}},$ for all $w_{1},$ $w_{2}$ $%
\in $ $\Bbb{G},$ and $L_{w}^{*}$ $=$ $L_{w^{-1}},$ for all $w$ $\in $ $\Bbb{G%
},$ and hence if $w$ is a reduced finite path, then $L_{w}$ is a partial
isometry and if $w$ is a vertex, then $L_{w}$ is a projection on $H_{G}.$ As
in [13], we are interested in right multiplication operators:

\strut

\begin{definition}
An operator $R_{w}$ on $H_{G}$ is defined by the right multiplication
operator with its symbol $\xi _{w}$ on $H_{G},$ for $w$ $\in $ $\Bbb{G}.$
i.e.,

\strut 

\begin{center}
$R_{w}\xi _{w^{\prime }}\overset{def}{=}\left\{ 
\begin{array}{ll}
\xi _{w^{\prime }}\xi _{w}=\xi _{w^{\prime }w} & \text{if }w^{\prime }\text{ 
}w\neq \emptyset  \\ 
\xi _{\emptyset }=0_{H_{G}} & \text{otherwise,}
\end{array}
\right. $
\end{center}

\strut 

for all $w,$ $w^{\prime }$ $\in $ $\Bbb{G}.$ The adjoint $R_{w}^{*}$ of $%
R_{w}$ is defined by $R_{w}^{*}$ $=$ $R_{w^{-1}},$ for all $w$ $\in $ $\Bbb{G%
}.$
\end{definition}

\strut

By definition, the product of two right multiplication operators $R_{w_{1}}$ 
$R_{w_{2}}$ is the multiplication operator $R_{w_{2}w_{1}}.$ i.e.,

$\strut $

\begin{center}
$R_{w_{1}}$ $R_{w_{2}}$ $=$ $R_{w_{2}w_{1}},$ \ for all $w_{1},$ $w_{2}$ $%
\in $ $\Bbb{G}.$
\end{center}

\strut

So, it is easy to check that the multiplication operators $R_{v},$ for all $%
v $ $\in $ $V(\widehat{G})$, are projections on $H_{G},$ since

$\strut $

\begin{center}
$R_{v}^{2}$ $=$ $R_{v}$ $R_{v}$ $=$ $R_{v^{2}}$ $=$ $R_{v}$ and $R_{v}^{*}$ $%
=$ $R_{v^{-1}}$ $=$ $R_{v}.$
\end{center}

\strut

Thus the right multiplication operators $R_{w},$ for all $w$ $\in $ $FP_{r}(%
\widehat{G}),$ are partial isometries; indeed,

\strut

\begin{center}
$R_{w}R_{w}^{*}R_{w}=R_{ww^{-1}w}=R_{w}$ \ 
\end{center}

and \ 

\begin{center}
$R_{w}^{*}R_{w}R_{w}^{*}=R_{w^{-1}ww^{-1}}=R_{w^{-1}}=R_{w}^{*},$
\end{center}

\strut

for all $w$ $\in $ $FP_{r}(\widehat{G}).$\strut

\strut

Now, we define certain groupoid actions $\beta $ of the graph groupoid $\Bbb{%
G}$ on $K$ $\otimes $ $H_{G},$ called \emph{right graph-representations},
where $K$ is a Hilbert space.

\strut \strut \strut \strut

\begin{definition}
Let $\Bbb{G}$ be a graph groupoid of a countable directed graph $G$ and let $%
M$ be a von Neumann algebra acting on a Hilbert space $K$ (i.e., $M$ $%
\subseteq $ $B(K)$)$.$ Define a right graph-representation (in short, a
right $G$-representation) $\beta $ $:$ $\Bbb{G}$ $\rightarrow $ $B(K$ $%
\otimes $ $H_{G})$, by a nonunital intertwined partial representation
satisfying that

\strut 

\begin{center}
$\beta _{w}(m)$ $R_{w}R_{w}^{*}=R_{w}^{*}mR_{w}=R_{w^{-1}}mR_{w},$
\end{center}

\strut 

for all $m$ $\in $ $M$ and $w$ $\in $ $\Bbb{G},$ and

\strut 

\begin{center}
$\beta _{v}(m)$ $=$ $m,$ for all $m$ $\in $ $M$ and $v$ $\in $ $V(\widehat{G}%
).$
\end{center}

\strut 

In this setting, we can regard $R_{w}$'s as $1_{K}$ $\otimes $ $R_{w}$'s on $%
K$ $\otimes $ $H_{G}$, for all $w$ $\in $ $\Bbb{G}.$ The above two rules for 
$\beta $ are called the $\beta $-relations on $M.$
\end{definition}

\strut \strut \strut \strut \strut \strut \strut \strut \strut \strut \strut
\strut \strut

For a fixed von Neumann algebra $M$ and a graph groupoid $\Bbb{G},$ we can
construct a \emph{groupoid crossed product algebra }$M$ $\times _{\beta }$ $%
\Bbb{G}.$

\strut

\begin{definition}
Let $M$ be a von Neumann algebra and $\Bbb{G},$ the graph groupoid of a
countable directed graph $G$ and let $\beta $ be a right $G$-representation.
Define the crossed product $\Bbb{M}_{G}$ $=$ $M$ $\times _{\beta }$ $\Bbb{G}$
of $M$ and $\Bbb{G},$ via $\beta ,$ by the von Neumann algebra generated by $%
M$ and $\{R_{w}$ $:$ $w$ $\in $ $\Bbb{G}\}$ in $B(K$ $\otimes $ $H_{G}),$
satisfying the $\beta $-relations on $M$. This von Neumann algebra $\Bbb{M}%
_{G}$ is called a right graph von Neumann algebra induced by $G$ over $M.$
\end{definition}

\strut \strut

\begin{remark}
In [10] and [11], we defined a graph von Neumann algebra $\Bbb{M}%
_{G}^{(left)}$ by the groupoid crossed product $M$ $\times _{\alpha }$ $\Bbb{%
G},$ where $\alpha $ $:$ $\Bbb{G}$ $\rightarrow $ $B(K$ $\otimes $ $H_{G})$
is a $G$-representation satisfying the $\alpha $-relations on $M$:

\strut 

\begin{center}
$\alpha _{w}(m)$ $L_{w}$ $L_{w}^{*}$ $=$ $L_{w}^{*}$ $m$ $L_{w},$ for all $w$
$\in $ $FP_{r}(\widehat{G}),$
\end{center}

and

\begin{center}
$\alpha _{v}(m)$ $=$ $m,$ for all $v$ $\in $ $V(\widehat{G}),$
\end{center}

\strut 

for all $m$ $\in $ $M,$ where $L_{w}$'s are the (left) multiplication
operators, for $w$ $\in $ $\Bbb{G}$. Notice that, in fact, the von Neumann
algebra $\Bbb{M}_{G}^{(left)}$ is the opposite algebra $\Bbb{M}_{G}^{op}$ of
a right graph von Neumann algebra $\Bbb{M}_{G}$ $=$ $M$ $\times _{\beta }$ $%
\Bbb{G}$, and hence they are anti-$*$-isomorphic from each other. This
guarantees that all results in [10] and [11] can be applicable in our
``right'' graph von Neumann algebraic setting.
\end{remark}

\strut \strut

Every operator $x$ in a right graph von Neumann algebra $\Bbb{M}_{G}$ $=$ $M$
$\times _{\beta }$ $\Bbb{G}$ has its expression,

\strut

\begin{center}
$x=\underset{w\in \Bbb{G}}{\sum }m_{w}$ $R_{w},$ \ for $m_{w}$ $\in $ $M.$
\end{center}

\strut \strut

\begin{lemma}
Let $m_{1}R_{w_{1}},$ ..., $m_{n}R_{w_{n}}$ be operators in a right graph
von Neumann algebra $\Bbb{M}_{G},$ for $n$ $\in $ $\Bbb{N}.$ Then

\strut \strut 

$\Pi _{k=1}^{n}\left( m_{k}R_{w_{k}}\right) $

\strut 

\begin{center}
$=\left\{ 
\begin{array}{ll}
\left(
m_{1}m_{2}^{w_{1}^{-1}}m_{3}^{(w_{1}w_{2})^{-1}}...m_{n}^{(w_{1}w_{2}...w_{n-1})^{-1}}\right) R_{w_{n}...w_{1}}
& \text{if }w_{n}...w_{1}\neq \emptyset  \\ 
&  \\ 
0_{\Bbb{M}_{G}} & \text{otherwise,}
\end{array}
\right. $
\end{center}

\strut 

where $m^{w}$ $\overset{def}{=}$ $\beta _{w}(m),$ for all $m$ $\in $ $M$ and 
$w$ $\in $ $\Bbb{G}.$ $\square $
\end{lemma}

\strut \strut

\strut \strut

\subsection{$M$-Valued Right Graph $W^{*}$-Probability Spaces}

\strut

\strut

Let $G$ be a countable directed graph with its graph groupoid $\Bbb{G}$ and
let $M$ be an arbitrary von Neumann algebra acting on a Hilbert space $K$.
In this section, we will define a $M$-diagonal right graph $W^{*}$%
-probability space $(\Bbb{M}_{G},$ $E),$ over its $M$-diagonal subalgebra $%
\Bbb{D}_{G}$, where $\Bbb{M}_{G}$ $=$ $M$ $\times _{\beta }$ $\Bbb{G}$ is a
right graph von Neumann algebra.

\strut

Let $v$ $\in $ $V(\widehat{G}).$ Then we can define a conditional
expectation $E_{v}$ $:$ $\Bbb{M}_{G}$ $\rightarrow $ $M$ $\cdot $ $R_{v},$
where $M$ $\cdot $ $R_{v}$ $=$ $\{m$ $R_{v}$ $:$ $m$ $\in $ $M\},$ by

\strut

\begin{center}
$E_{v}\left( \underset{w\in \Bbb{G}}{\sum }m_{w}R_{w}\right) =m_{v}R_{v},$
\end{center}

\strut

for all $\underset{w\in \Bbb{G}}{\sum }$ $m_{w}R_{w}$ $\in $ $\Bbb{M}_{G}.$
Notice that each $M$ $\cdot $ $R_{v}$ is a $W^{*}$-subalgebra of $\Bbb{M}%
_{G},$ which is $*$-isomorphic to $M,$ since $R_{v}$ is a projection, for
all $v$ $\in $ $V(\widehat{G}).$ A pair $(\Bbb{M}_{G},$ $E_{v})$ is a $M$%
-valued $W^{*}$-probability space, for $v$ $\in $ $V(\widehat{G}).$ We call
it a \emph{vertex-depending }(or the $v$-depending) \emph{right graph }$%
W^{*} $\emph{-probability space over }$M.$ The conditional expectation $E_{v}
$ is said to be a \emph{vertex-depending} (or the $v$-depending) \emph{%
conditional expectation} \emph{for} $v$ $\in $ $V(\widehat{G}).$ By the very
definition, the right graph von Neumann algebra $\Bbb{M}_{G}$ has $\left| V(%
\widehat{G})\right| $-many vertex-depending $W^{*}$-probability spaces over $%
M.$

\strut

\begin{definition}
By $\Bbb{D}_{G},$ denote the $W^{*}$-subalgebra $\underset{v\in V(\widehat{G}%
)}{\oplus }$ $(M$ $\cdot $ $R_{v})$ of $\Bbb{M}_{G}$ $=$ $M$ $\times _{\beta
}$ $\Bbb{G}.$ This subalgebra $\Bbb{D}_{G}$ is called the $M$-diagonal
subalgebra of $\Bbb{M}_{G}.$ Define a conditional expectation $E$ $:$ $\Bbb{M%
}_{G}$ $\rightarrow $ $\Bbb{D}_{G}$, by $\underset{v\in V(\widehat{G})}{%
\oplus }$ $E_{v},$ where $E_{v}$'s are the $v$-depending conditional
expectations, for all $v$ $\in $ $V(\widehat{G}).$ i.e.,

\strut 

\begin{center}
$E\left( \underset{w\in \Bbb{G}}{\sum }m_{w}R_{w}\right) =\underset{v\in V(%
\widehat{G})}{\sum }m_{v}R_{v},$
\end{center}

\strut 

for all $\underset{w\in \Bbb{G}}{\sum }m_{w}R_{w}$ $\in $ $\Bbb{M}_{G}.$ The
pair $(\Bbb{M}_{G},$ $E)$ is called the $M$-diagonal right graph $W^{*}$%
-probability space over $\Bbb{D}_{G}.$ We call $E$, the canonical
conditional expectation.
\end{definition}

\strut

Notice that the $M$-diagonal subalgebra \strut \strut $\Bbb{D}_{G}$ of $\Bbb{%
M}_{G}$ $=$ $M$ $\times _{\beta }$ $\Bbb{G}$ is uniquely determined for all
right $G$-representations $\beta .$

\strut

\begin{definition}
Let $G$ be a countable directed graph and let $\Bbb{G}$ be the graph
groupoid of $G.$ Define a map $\delta $ $:$ $\Bbb{G}$ $\rightarrow $ $\Bbb{G}
$ by mapping $w$ $\in $ $\Bbb{G}$ to the graphical image $\delta (w)$ of $w,$
for all $w$ $\in $ $\Bbb{G}.$ The map $\delta $ on $\Bbb{G}$ is called the
diagram map and the graphical image $\delta (w)$ of $w$ is called the
diagram of $w$, for all $w$ $\in $ $\Bbb{G}.$ We say that the elements $w_{1}
$ and $w_{2}$ in $\Bbb{G}$ are diagram-distinct, if (i) $w_{1}$ $\neq $ $%
w_{2}^{-1}$ and (ii) $\delta (w_{1})$ $\neq $ $\delta (w_{2}).$ Suppose that 
$X_{1}$ and $X_{2}$ are subsets in $\Bbb{G}.$ They are said to be
diagram-distinct if, for any pair $(w_{1},$ $w_{2})$ in $X_{1}$ $\times $ $%
X_{2},$ $w_{1}$ and $w_{2}$ are diagram-distinct.
\end{definition}

\strut

Let $l$ be a loop edge in $E(\widehat{G}),$ and let $w_{1}$ $=$ $l^{k_{1}}$
and $w_{2}$ $=$ $l^{k_{2}}$ in $\Bbb{G},$ for $k_{1}$ $\neq $ $k_{2}$ $\in $ 
$\Bbb{N}.$ Clearly, $w_{1}$ $\neq $ $w_{2}^{-1}.$ But the elements $w_{1}$
and $w_{2}$ are not diagram-distinct, because their diagrams $\delta (w_{1})$
and $\delta (w_{2})$ are identical to the diagram $\delta (l)$ $=$ $l$ of $%
l. $ Suppose either $w_{1}$ or $w_{2}$ is a non-loop finite path such that
(i) $w_{1}$ $\neq $ $w_{2}^{-1},$ and (ii) $\delta (w_{k})$'s are not loop
finite paths, for $k$ $=$ $1,$ $2.$ Then they are diagram-distinct whenever $%
w_{1}$ and $w_{2}$ are distinct in $\Bbb{G}.$

\strut \strut \strut

\begin{theorem}
The subsets $M$ $\cdot $ $R_{w_{1}}$ and $M$ $\cdot $ $R_{w_{2}}$ of a right
graph von Neumann algebra $\Bbb{M}_{G}$ $=$ $M$ $\times _{\beta }$ $\Bbb{G}$
are free over $D_{G}$ in $(\Bbb{M}_{G},$ $E)$ if and only if $w_{1}$ and $%
w_{2}$ are diagram-distinct in $\Bbb{G}$, where $M$ $\cdot $ $R_{w}$ $%
\overset{def}{=}$ $\{m$ $L_{w}$ $:$ $m$ $\in $ $M\},$ for all $w$ $\in $ $%
\Bbb{G}.$ $\square $
\end{theorem}

\strut \strut \strut \strut \strut

The readers can find the proof showing that: $M$ $\cdot $ $L_{w_{1}}$ and $M$
$\cdot $ $L_{w_{2}}$ are free over $\Bbb{D}_{G}$ in a (left) graph von
Neumann algebra $M$ $\times _{\alpha }$ $\Bbb{G}$, if and only if $w_{1}$
and $w_{2}$ are diagram-distinct, in [10] and [11], where $L_{w_{1}}$ and $%
L_{w_{2}}$ are left multiplication operators. By the slight modification of
it, we can prove the above theorem. As corollary of the previous theorem, if 
$e_{1}$ and $e_{2}$ are distinct edges in $E(G)$ ($\subset $ $E(\widehat{G})$%
), then $M$ $\cdot $ $R_{e_{1}}$ and $M$ $\cdot $ $R_{e_{2}}$ are free over $%
\Bbb{D}_{G}$ in $(\Bbb{M}_{G},$ $E),$ and vice versa. Indeed, if two edges $%
e_{1}$ and $e_{2}$ are distinct in $E(G),$ then $e_{1}$ and $e_{2}^{\pm 1}$
(resp., $e_{2}$ and $e_{1}^{\pm 1}$) are diagram-distinct in $\Bbb{G}$.

\strut

Let's consider the subset $\{e,$ $e^{-1}\}$ of $\Bbb{G},$ for $e$ $\in $ $%
E(G)$ $\subset $ $E(\widehat{G}).$ Then we can define a subset $\Bbb{G}_{e}$
of $\Bbb{G}$, by the collection of all reduced words in $\{e,$ $e^{-1}\}.$
It is easy to check that: if $e$ is a loop edge in $E(G)$, then $\Bbb{G}_{e}$
is a group, which is group-isomorphic to the infinite abelian cyclic group $%
\Bbb{Z},$ and if $e$ is a non-loop edge in $E(G)$, then $\Bbb{G}_{e}$ is a
groupoid $\{\emptyset ,$ $v_{1},$ $v_{2},$ $e,$ $e^{-1}\},$ where $e$ $=$ $%
v_{1}$ $e$ $v_{2},$ with $v_{1},$ $v_{2}$ $\in $ $V(G)$.

\strut

Notice that the graph groupoid $\Bbb{G}$ is the reduced free product $%
\underset{e\in E(G)}{*}$ $\Bbb{G}_{e}$ of $\Bbb{G}_{e}$'s (in the sense of
[10]). Moreover, each $\Bbb{G}_{e}$ can be regarded as a new graph groupoid
induced by a one-edge graph $G_{e},$ with $E(G_{e})$ $=$ $\{e\}.$ Clearly,
here, $G_{e}$ is a full-subgraph of $G,$ and $\Bbb{G}_{e}$ is a subgroupoid
of $\Bbb{G},$ for all $e$ $\in $ $E(G).$ Construct the $W^{*}$-subalgebras,

\strut

\begin{center}
$\Bbb{M}_{e}$ $\overset{def}{=}$ $vN\left( M\times _{\beta _{e}}\Bbb{G}_{e},%
\text{ }\Bbb{D}_{G}\right) $ \strut with $\beta _{e}$ $\overset{def}{=}$ $%
\beta $ $\mid _{\Bbb{G}_{e}}$ on $M,$
\end{center}

\strut

of $\Bbb{M}_{G},$ for all $e$ $\in $ $E(G).$ By the previous theorem, we can
have that:

\strut

\begin{corollary}
The $W^{*}$-subalgebras $\Bbb{M}_{e}$'s of $\Bbb{M}_{G}$ are free over $\Bbb{%
D}_{G}$, from each other, for all $e$ $\in $ $E(G),$ in $(\Bbb{M}_{G},$ $E).$
$\square $
\end{corollary}

\strut \strut \strut \strut \strut \strut

\begin{definition}
The $W^{*}$-subalgebras $\Bbb{M}_{e}$ of $\Bbb{M}_{G}$ are called $\Bbb{D}%
_{G}$-free blocks of $\Bbb{M}_{G},$ for all $e$ $\in $ $E(G).$
\end{definition}

\strut \strut \strut

Notice that $\Bbb{M}_{e}$ is (not only $*$-isomorphic but also) identically
same as $\Bbb{M}_{e^{-1}}$ in $\Bbb{M}_{G},$ for all $e$ $\in $ $E(G).$\ So,
we can get the following theorem characterizing the $\Bbb{D}_{G}$-free
structure of $(\Bbb{M}_{G},$ $E).$\strut \ The following theorem is the
direct consequence of the previous corollary, and the fact that $\{\Bbb{M}%
_{e}$ $:$ $e$ $\in $ $E(G)\}$ generates $\Bbb{M}_{G}.$

\strut

\begin{theorem}
Let $(\Bbb{M}_{G},$ $E)$ be the $M$-diagonal right graph $W^{*}$-probability
space over its $M$-diagonal subalgebra $\Bbb{D}_{G}.$ Then

\strut 

\begin{center}
$\Bbb{M}_{G}$ $\overset{*\text{-isomorphic}}{=}$ $\underset{e\in E(G)}{%
\,\,*_{\Bbb{D}_{G}}}$ $\Bbb{M}_{e},$
\end{center}

\strut 

where $\Bbb{M}_{e}$'s are the $\Bbb{D}_{G}$-free blocks of $\Bbb{M}_{G},$
for all $e$ $\in $ $E(G).$ $\square $
\end{theorem}

\strut \strut

\strut Recall that if $(A,$ $E_{B})$ is a $B$-valued $W^{*}$-probability
space over $B,$ and if $\{A_{i}\}_{i\in I}$ are $W^{*}$-subalgebras of $A,$
over their common $W^{*}$-subalgebra $B,$ then the $B$-free product algebra $%
\underset{i\in I}{\,\,*_{B}}$ $A_{i}$ can be expressed by

\strut

\begin{center}
$B\oplus \left( \oplus _{n=1}^{\infty }\left( \underset{i_{1}\neq
i_{2},i_{2}\neq i_{3},...,i_{n-1}\neq i_{n}\in I}{\oplus }%
(A_{i_{1}}^{o}\otimes _{B}...\otimes _{B}A_{i_{n}}^{o})\right) \right) ,$
\end{center}

\strut

as a Banach space, where $A_{i_{k}}^{o}$ $\overset{def}{=}$ $A_{i_{k}}$ $%
\ominus $ $B,$ for all $k$ $=$ $1,$ ..., $n.$ Since our right graph von
Neumann algebra $\Bbb{M}_{G}$ is $*$-isomorphic to the $\Bbb{D}_{G}$-free
product $\underset{e\in E(G)}{\,\,*_{\Bbb{D}_{G}}}$ $\Bbb{M}_{e},$ where $%
\Bbb{M}_{e}$'s are the $\Bbb{D}_{G}$-free blocks of $\Bbb{M}_{G},$ we can
get a similar Banach space expression. However, the $\Bbb{D}_{G}$-tensor
products in our $\Bbb{D}_{G}$-free product algebra is dependent upon the
admissibility on $\Bbb{G}.$ i.e., if two edges $e_{1}$ and $e_{2}$ are
totally distinct, in the sense that all pairs $(e_{1}^{\symbol{94}},$ $%
e_{2}^{\symbol{94}})$ are diagram-distinct, for $e_{k}^{\symbol{94}}$ $\in $ 
$\{e_{k},$ $e_{k}^{-1}\},$ for $k$ $=$ $1,$ $2,$ then $\Bbb{M}_{e_{1}}^{o}$ $%
\otimes _{\Bbb{D}_{G}}$ $\Bbb{M}_{e_{2}}^{o}$ is Banach-space isomorphic to $%
\{0_{\Bbb{D}_{G}}\},$ where $\Bbb{M}_{e_{k}}^{o}$ $=$ $\Bbb{M}_{e_{k}}$ $%
\ominus $ $\Bbb{D}_{G},$ for $k$ $=$ $1,$ $2.$

\strut

\begin{definition}
Let $G$ be a countable directed graph. Define the subset $E(\widehat{G}%
)_{r}^{*}$ of $FP_{r}(\widehat{G})$ by

\strut \strut 

\begin{center}
$E(\widehat{G})_{r}^{*}$ $\overset{def}{=}$ $E(\widehat{G})$ $\cup $ $\left(
\cup _{k=2}^{\infty }\left\{ e_{1}e_{2}...e_{k-1}e_{k}\left| 
\begin{array}{c}
e_{1}\text{ ... }e_{k}\in FP_{r}(\widehat{G}), \\ 
e_{1}\neq e_{2}^{\pm 1},\text{ }e_{2}\neq e_{3}^{\pm 1} \\ 
\text{ \ \ }...,\text{ }e_{k-1}\neq e_{k}^{\pm 1}
\end{array}
\right. \right\} \right) .$
\end{center}
\end{definition}

\strut \strut

\strut Similar to [10] and [11], we can get the following theorem.

\strut \strut

\begin{theorem}
Let $G$ be a countable directed graph and let $\Bbb{M}_{G}$ $=$ $M$ $\times
_{\beta }$ $\Bbb{G}$ be the corresponding right graph von Neumann algebra,
for a fixed von Neumann algebra $M.$ Then, as a Banach space,

\strut 

\begin{center}
$\Bbb{M}_{G}$ $=$ $\Bbb{D}_{G}\oplus \left( \underset{w\in E(\widehat{G}%
)_{r}^{*}}{\oplus }\,\Bbb{M}_{w}^{o}\right) $
\end{center}

with

\begin{center}
$\Bbb{M}_{w}^{o}$ $\overset{def}{=}$ $\Bbb{M}_{e_{1}}^{o}\otimes _{\Bbb{D}%
_{G}}$ ... $\otimes _{\Bbb{D}_{G}}$ $\Bbb{M}_{e_{n}}^{o},$
\end{center}

where

\begin{center}
$\Bbb{M}_{e}^{o}$ $\overset{def}{=}$ $\Bbb{M}_{e}$ $\ominus $ $\Bbb{D}_{G},$
for all $e$ $\in $ $E(\widehat{G}),$

\strut 
\end{center}

\strut whenever $w$ $=$ $e_{1}$ ... $e_{n}$ $\in $ $E(\widehat{G})_{r}^{*},$
with $e_{1},$ ..., $e_{n}$ $\in $ $E(\widehat{G}),$ under the identification 
$\Bbb{M}_{e}$ $=$ $M_{e^{-1}},$ for all $e$ $\in $ $E(\widehat{G}).$ $%
\square $
\end{theorem}

\strut \strut \strut

\textbf{Notation} \ To emphasize that the $\Bbb{D}_{G}$-valued free product
algebra $\underset{e\in E(G)}{\,*_{\Bbb{D}_{G}}}$ $\Bbb{M}_{e}$ has the
above Banach-space expression determined by the admissibility on $\Bbb{G},$
we will denote it by $\underset{e\in E(G)}{\,*_{\Bbb{D}_{G}}^{r}}$ $\Bbb{M}%
_{e}$. And we will call it the $\Bbb{D}_{G}$-valued \emph{reduced} free
product of the $\Bbb{D}_{G}$-free blocks $\Bbb{M}_{e}$'s. $\square $

\strut

In [11], the $\Bbb{D}_{G}$-free blocks in a (left) graph von Neumann algebra
are characterized. Similarly, we can have the following results.

\strut

\begin{theorem}
Let $\Bbb{M}_{G}$ $=$ $M$ $\times _{\beta }$ $\Bbb{G}$ be a right graph von
Neumann algebra. Then its $\Bbb{D}_{G}$-free block $\Bbb{M}_{e}$ is $*$%
-isomorphic to either $M$ $\times _{\lambda _{e}}$ $\Bbb{Z}$ (whenever $e$
is a loop-edge) or $M$ $\otimes _{\beta _{e}}$ $M_{2}(\Bbb{C})$ (whenever $e$
is a non-loop edge), where $\lambda _{e}$ is the group action of $\Bbb{Z},$
such that $\lambda _{n}(m)$ $=$ $\beta _{e^{n}}(m),$ for all $m$ $\in $ $M$
and $n$ $\in \Bbb{Z}.$ $\square $
\end{theorem}

\strut \strut

\textbf{Assumption} In the rest of this paper, we will restrict our
interests to the case where $M$ $=$ $\Bbb{C}.$ And we call a graph von
Neumann algebra $\Bbb{M}_{G}$ $=$ $\Bbb{C}$ $\times _{\beta }$ $\Bbb{G},$
``the'' right graph von Neumann algebra induced by $G.$ $\square $

\strut

Suppose $M$ $=$ $\Bbb{C}.$ Then a right graph von Neumann algebra $\Bbb{M}%
_{G}$ $=$ $\Bbb{C}$ $\times _{\beta }$ $\Bbb{G}$ is $*$-isomorphic to the
groupoid $W^{*}$-subalgebra $\overline{\Bbb{C}[\Bbb{G}]}^{w}$ of $B(H_{G}),$
for all right graph-representation $\beta ,$ where $H_{G}$ is the graph
Hilbert space. Also, in such case, the $\Bbb{D}_{G}$-free blocks $\Bbb{M}%
_{e} $'s of $\Bbb{M}_{G}$ are $*$-isomorphic to either

$\strut $

\begin{center}
$L(\Bbb{Z})$ $\overset{*\text{-isomorphic}}{=}$ $L^{\infty }(\Bbb{T})$ \quad
or\quad $M_{2}(\Bbb{C}).$
\end{center}

\strut \strut

\begin{example}
Let $G_{N}$ be a one-vertex-$N$-loop-edge graph with its vertex set $V(G_{N})
$ $=$ $\{v\}$ and $E(G_{N})$ $=$ $\{e_{j}$ $=$ $v$ $e_{j}$ $v$ $:$ $j$ $=$ $%
1,$ ..., $N\}.$ Then the right graph von Neumann algebra $\Bbb{C}_{G_{N}}$ $=
$ $\Bbb{C}$ $\times _{\beta }$ $\Bbb{G}_{N},$ where $\Bbb{G}_{N}$ is the
graph groupoid of $G_{N},$ is $*$-isomorphic to the free group factor $%
L(F_{N})$ (Also See [10]). Notice that the graph groupoid $\Bbb{G}_{N}$ of $%
G_{N}$ is a group and moreover it is group-isomorphic to the free group $%
F_{N}$ $=$ $<$ $g_{1},$ ..., $g_{N}$ $>$ with $N$-generators $\{g_{1},$ ..., 
$g_{N}\}$. Indeed, the edges $e_{1},$ ..., $e_{N}$ are generators of $\Bbb{G}%
_{N}$ as a group. So, there is a natural generator-preserving
group-isomorphism between $\Bbb{G}_{N}$ and $F_{N}.$ Also, notice that the
vertex $v$ in $\Bbb{G}_{N}$ is the group identity. And all partial
isometries $R_{e_{j}^{n}}$ correspond to the unitary representation $%
u_{g_{j}^{n}}$ of $g_{j}^{n},$ for all $n$ $\in $ $\Bbb{Z}$ $\setminus $ $%
\{0\}$ and for $j$ $=$ $1,$ ..., $N.$ Therefore, the right graph von Neumann
algebra $\Bbb{C}_{G_{N}}$ $=$ $\Bbb{C}$ $\times _{\beta }$ $\Bbb{G}_{N}$ is $%
*$-isomorphic to the classical crossed product algebra $\Bbb{C}$ $\times
_{\lambda }$ $F_{N},$ whenever $\lambda $ $=$ $\beta .$ Notice that $\Bbb{C}$
$\times _{\beta }$ $\Bbb{G}_{N}$ is $*$-isomorphic to $\overline{\Bbb{C}[%
\Bbb{G}_{N}]}^{w},$ for any $\beta ,$ by the linearity of $\beta $ on $\Bbb{C%
}.$ Also, it is well-known that a crossed product algebra $\Bbb{C}$ $\times
_{\gamma }$ $\Gamma $ is $*$-isomorphic to $\overline{\Bbb{C}[\Gamma ]}^{w},$
for any group action $\gamma ,$ for all groups $\Gamma .$ Therefore, the
right graph von Neumann algebra $\Bbb{C}_{G_{N}}$ is $*$-isomorphic to the
group von Neumann algebra $\overline{\Bbb{C}[F_{N}]}^{w}.$ More generally,
the right graph von Neumann algebra $\Bbb{M}_{G_{N}}$ $=$ $M$ $\times
_{\beta }$ $\Bbb{G}_{N}$ is $*$-isomorphic to the classical crossed product
von Neumann algebra $M$ $\times _{\lambda }$ $F_{N},$ whenever $\lambda $ $=$
$\beta .$
\end{example}

\strut \strut \strut \strut \strut

\begin{example}
Let $C_{N}$ be the one-flow circulant graph with $V(C_{N})$ $=$ $\{v_{1},$
..., $v_{N}\}$ and $E(C_{N})$ $=$ $\{e_{j}$ $=$ $v_{j}$ $e_{j}$ $v_{j+1}$ $:$
$j$ $=$ $1,$ ..., $N,$ $v_{N+1}$ $\overset{def}{=}$ $v_{1}\}.$ Then the
right graph von Neumann algebra $\Bbb{M}_{C_{N}}$ $=$ $\overline{\Bbb{C}%
\left( \Bbb{G}_{C_{N}}\right) }^{w}$ contains $W^{*}$-subalgebras which are $%
*$-isomorphic to the group von Neumann algebra $L(\Bbb{Z})$ $=$ $\overline{%
\Bbb{C}[\Bbb{Z}]}^{w}.$ Since a finite path $w$ $=$ $e_{1}$ ... $e_{N}$
induces the subset $\Bbb{S}_{w}$ of $\Bbb{G}_{C_{N}}$ consisting of all
reduced words in $S_{w}$ $=$ $\{w,$ $w^{-1}\}$ in $\Bbb{G}_{C_{N}},$ we have
the $W^{*}$-subalgebra $N$ $=$ $\overline{\Bbb{C}\left( \Bbb{S}_{w}\right) }%
^{w}$ in $\Bbb{M}_{C_{N}}.$ We can regard this $W^{*}$-subalgebra $N$ as the
right graph von Neumann algebra induced by the graph $G_{w}$ with its vertex
set $V(G_{w})$ $=$ $\{v_{1}\}$ and its edge set $E(G_{w})$ $=$ $\{w$ $=$ $%
v_{1}$ $w$ $v_{1}\}.$ Then it is an one-vertex-one-loop-edge graph. So, by
the previous example, the von Neumann algebra $N$ is $*$-isomorphic to $L(%
\Bbb{Z}),$ where $\Bbb{Z}$ is the infinite abelian cyclic group. It is
well-know that the group von Neumann algebra $L(\Bbb{Z})$ is $*$-isomorphic
to $L^{\infty }(\Bbb{T}),$ where $\Bbb{T}$ is the unit circle in $\Bbb{C}.$
Consider $w^{k},$ for all $k$ $\in $ $\Bbb{Z}$ $\setminus $ $\{0\}.$ Then $%
W^{*}$-subalgebras $N$ $=$ $\Bbb{C}$ $\times _{\beta }$ $\Bbb{G}_{w^{k}}$ of 
$\Bbb{C}_{C_{N}}$ are all $*$-isomorphic to $L(\Bbb{Z}).$

\strut 

We can easily check that $E(C_{N})_{r}^{*}$ $=$ $FP(C_{N})$ $\cup $ $%
FP(C_{N}^{-1}),$ for all $N$ $\in $ $\Bbb{N}.$ Therefore,

\strut 

$\ \ \ \ \ \ \ \ \Bbb{M}_{C_{N}}$ $=$ $\underset{k=1}{\overset{N}{\,\,\,*_{%
\Bbb{D}_{C_{N}}}}}$ $\Bbb{M}_{e_{k}}$

\strut 

where $\Bbb{M}_{e_{k}}$ $\overset{def}{=}$ $vN\left( \overline{\Bbb{C}\left( 
\Bbb{G}_{e}\right) }^{w},\text{ }\Bbb{D}_{C_{N}}\right) ,$ for all $k$ $=$ $%
1,$ ..., $N$

\strut 

\ \ $\ \ \ \ \ \ \ \ \ \ \ \ =\Bbb{D}_{C_{N}}$ $\oplus $ $\left( \underset{%
w^{*}\in E(C_{N})_{r}^{*}}{\oplus }\text{ }\Bbb{M}_{w^{*}}^{o}\right) $

\strut 

$\ \ \ \ \ \ \ \ \ \ \ \ \ \ =$ $\Bbb{D}_{C_{N}}\oplus \left( \underset{w\in
FP(C_{N})\cup FP(C_{N}^{-1})}{\oplus }\,\Bbb{M}_{w}^{o}\right) ,$

\strut 

where $\Bbb{M}_{w}^{o}$ $\overset{def}{=}$ $\Bbb{M}_{e^{(1)}}^{o}$ $\otimes
_{\Bbb{D}_{C_{N}}}$ ... $\otimes _{\Bbb{D}_{C_{N}}}$ $\Bbb{M}_{e^{(n)}}^{o},$
whenever $w$ $=$ $e^{(1)}$ ... $e^{(n)}$ in $E(C_{N})_{r}^{*},$ for $n$ $\in 
$ $\Bbb{N}.$
\end{example}

\strut \strut

\strut

\strut

\section{Labeled Graph Groupoids and Graph Automata}

\strut

\strut \strut \strut

We will give suitable weights to the elements of graph groupoids induced by
the given locally finite connected directed graphs. Then these weights give
more accurate information of the admissibility on graph groupoids. Such
admissibility conditions are explained by the labeling map and the shifting
(map) of the automata constructed by graph groupoids. Also, this process
would show how the admissibility on a graph groupoid works on the
corresponding graph Hilbert space.

\strut

Recall that a countable directed graph $G$ is \emph{locally finite}, if
every vertex of $G$ has only finitely many incident edges (equivalently, the
degree of each vertex is finite). Also, we say that $G$ is \emph{connected},
if, for any pair $(v_{1},$ $v_{2})$ of distinct vertices ($v_{1}$ $\neq $ $%
v_{2}$), there always exists at least one reduced finite path $w$ $\in $ $%
FP_{r}(\widehat{G})$ such that $w$ $=$ $v_{1}$ $w$ $v_{2}$ and $w^{-1}$ $=$ $%
v_{2}$ $w^{-1}$ $v_{1}.$

\strut

\textbf{Assumption} From now on, all given directed graphs are locally
finite and connected. $\square $

\strut

Let $G$ be a locally finite connected directed graph and let $v_{0}$ be an
any fixed vertex of $G.$ Then we can define the out-degree $\deg
_{out}(v_{0})$, the in-degree $\deg _{in}(v_{0}),$ and the degree $\deg
(v_{0})$ of $v_{0}$ as follows:

\strut

\begin{center}
$\deg _{out}(v_{0})$ $\overset{def}{=}$ $\left| \{e\in E(G):e=v_{0}\text{ }%
e\}\right| ,$
\end{center}

\strut

\begin{center}
$\deg _{in}(v_{0})$ $\overset{def}{=}$ $\left| \{e\in E(G):e=\text{ }e\text{ 
}v_{0}\}\right| ,$
\end{center}

and

\begin{center}
$\deg (v_{0})$ $\overset{def}{=}$ $\deg _{out}(v_{0})$ $+$ $\deg
_{in}(v_{0}).$
\end{center}

\strut \strut

By the locally-finiteness of $G,$ the above three cardinalities are less
than $\infty .$ We are interested in out-degrees.

\strut

\strut

\strut

\subsection{Labeled Graph Groupoids}

\strut

\strut

Let $G$ be a locally finite connected directed graph. By the locally
finiteness of $G,$ we can have that

\strut

\begin{center}
$\max \{\deg _{out}(v)$ $:$ $v$ $\in $ $V(G)\}$ $<$ $\infty .$
\end{center}

\strut \strut

Denote the maximal value of the out-degrees of all vertices by $N$. Then we
can define a set

\strut

\begin{center}
$Y$ $=$ $\{1,$ ..., $N\}$ and $Y_{0}$ $=$ $\{0\}$ $\cup $ $Y.$
\end{center}

\strut

Define \emph{lattices} $l_{0},$ $l_{1},$ ..., $l_{N}$ in the 2-dimensional
real vector space $\Bbb{R}^{2},$ by the vectors

\strut

\begin{center}
$l_{0}$ $=$ $\overrightarrow{(0,0)}$ and $l_{k}$ $=$ $\overrightarrow{(1,%
\text{ }e^{k})},$
\end{center}

\strut

for all $k$ $=$ $1,$ ..., $N.$ To distinguish the vectors and points in $%
\Bbb{R}^{2},$ we use the notation $\overrightarrow{(t_{1},\text{ }t_{2})},$
as a vector connecting $(0,$ $0)$ and $(t_{1},$ $t_{2})$ in $\Bbb{R}^{2},$
and $(t_{1},$ $t_{2}),$ as a point in $\Bbb{R}^{2},$ for $t_{1},$ $t_{2}$ $%
\in $ $\Bbb{R}.$ For the given lattices $l_{0},$ $l_{1},$ ..., $l_{N},$
determine the \emph{inverse lattices} $l_{-0},$ $l_{-1},$ ..., $l_{-N}$ by

\strut

\begin{center}
$l_{-0}$ $=$ $\overrightarrow{(0,0)}$ $=$ $l_{0},$ and $l_{-k}$ $=$ $%
\overrightarrow{(1,\text{ }-e^{k})},$
\end{center}

\strut

for all $k$ $=$ $1,$ ..., $N.$ Since $l_{-0}$ $=$ $l_{0},$ it is the
self-inverse lattice. i.e., we can understand the lattices $l_{1},$ ..., $%
l_{N}$ are the upward lattices in $\Bbb{R}^{2},$ and the inverse lattices $%
l_{-1},$ ..., $l_{-N}$ are the downward lattices in $\Bbb{R}^{2}.$

\strut

Define now the \emph{lattice paths on} $\Bbb{R}^{2},$ \emph{generated by the
lattices} $l_{\pm 1},$ ..., $l_{\pm N},$ by the diagrams in $\Bbb{R}^{2},$
constructed by the rule:

\strut

\begin{center}
$l_{i}$ $l_{j}$ $=$ the diagram connecting $l_{i}$ and $l_{j}$ by
identifying the end point $(1,$ $t_{i}e^{i})$ of $l_{i}$ and the starting
point $(0,$ $0)$ of the $l_{j},$
\end{center}

\strut

where

\begin{center}
$t_{i}$ $\overset{def}{=}$ $\left\{ 
\begin{array}{ll}
1 & \text{if }i\in \{1,\text{ }...,\text{ }N\} \\ 
-1 & \text{if }i\in \{-1,\text{ ..., }-N\},
\end{array}
\right. $
\end{center}

\strut

for all $i$ $\in $ $\{\pm 1,$ ..., $\pm N\}.$ For instance, if $l_{1}$ $=$ $%
(0,$ $e)$ and $l_{-3}$ $=$ $(0,$ $-e^{3}),$ then a lattice path $l_{1}$ $%
l_{-3}$ is the diagram in $\Bbb{R}^{2},$ starting at $(0,$ $0)$ and ending
at $(2,$ $e$ $-$ $e^{3}),$ via $(1,$ $e).$ We can do the above process
inductively, and hence we can create the lattice paths $l_{i_{1}}$ $%
l_{i_{2}} $ ... $l_{i_{n}}$ on $\Bbb{R}^{2},$ for all $(i_{1},$ ..., $i_{n})$
$\in $ $\{\pm 1,$ ..., $\pm N\}^{n},$ for $n$ $\in $ $\Bbb{N}.$ Suppose we
have the above lattice path $l_{1}$ $l_{-3},$ and let $l_{2}$ $=$ $(1,$ $%
e^{2}).$ Then, inductively, we can have the lattice path $l_{1}$ $l_{-3}$ $%
l_{2},$ by identifying the ending point $(2,$ $e$ $-$ $e^{3})$ of $l_{1}$ $%
l_{-3}$ and the starting point $(0,$ $0)$ of $l_{2}.$ Then it is a lattice
path starting at $(0,$ $0)$ and ending at $(3,$ $e$ $-$ $e^{3}$ $+$ $e^{2}),$
via $(1,$ $e) $ and $(2$, $e$ $-$ $e^{3})$ (Also, see [17] and [41]).

\strut

Let $l$ $=$ $l_{i_{1}}$ $l_{i_{2}}$ ... $l_{i_{n}}$ be a lattice path
generated by $l_{\pm 1},$ ..., $l_{\pm N}.$ Then the \emph{length }$\left|
l\right| $\emph{\ of} $l$ is defined to be $n,$ the cardinality of lattices
generating $l.$

\strut

\begin{definition}
Let $N$ $=$ $\max \{\deg _{out}(v)$ $:$ $v$ $\in $ $V(G)\}$ $\in $ $\Bbb{N},$
and let $X$ $=$ $\{l_{\pm 1},$ ..., $l_{\pm N}\}$ be the lattices defined as
above. Define the set $X_{0}$ by $X$ $\cup $ $\{l_{0}\}.$ The set $\mathcal{L%
}_{N}$ is the collection of all lattice paths generated by $X,$ containing $%
l_{0}.$ Define a subset $\mathcal{L}_{N}(k)$ of $\mathcal{L}_{N}$ by

\strut 

\begin{center}
$\mathcal{L}_{N}(k)$ $\overset{def}{=}$ $\{l$ $\in $ $\mathcal{L}_{N}$ $:$ $%
\left| l\right| $ $=$ $k\},$
\end{center}

\strut 

for all $k$ $\in $ $\Bbb{N}.$
\end{definition}

\strut

By definition, it is clear that

\strut

\begin{center}
$\mathcal{L}_{N}$ $=$ $\{l_{0}\}$ $\sqcup $ $\left( \underset{k=1}{\overset{%
\infty }{\sqcup }}\left( \mathcal{L}_{N}(k)\right) \right) ,$
\end{center}

\strut

where $\sqcup $ means the disjoint union. Suppose $l$ $=$ $l_{i_{1}}$ $%
l_{i_{2}}$ ... $l_{i_{n}}$ $\in $ $\mathcal{L}_{N}(n),$ for $(i_{1},$ ..., $%
i_{n})$ $\in $ $\{\pm 1,$ ..., $\pm N\}^{n},$ for $n$ $\in $ $\Bbb{N}.$ Also
assume that this length-$n$ lattice path $l$ starts at $(0,$ $0)$ and ends
at $(0,$ $0)$ in $\Bbb{R}^{2}.$ Then we say that the lattice path $l$ \emph{%
satisfies the axis property}. It is easy to check that:

\strut

\begin{lemma}
If $n$ is odd, then every length-$n$ lattice path does not satisfy the axis
property. $\square $
\end{lemma}

\strut

Define now the subset $\mathcal{L}_{N}^{o}(k)$ of $\mathcal{L}_{N}(k)$
satisfying the axis property.

\strut

\begin{definition}
Let $\mathcal{L}_{N}(k)$ be the set of all length-$k$ lattice paths on $\Bbb{%
R}^{2}$. Define the subset $\mathcal{L}_{N}^{o}(k)$ of $\mathcal{L}_{N}(k)$
by

\strut 

\begin{center}
$\mathcal{L}_{N}^{o}(k)$ $\overset{def}{=}$ $\{l$ $\in $ $\mathcal{L}_{N}(k)$
$:$ $l$ satisfies the axis property$\},$
\end{center}

\strut 

for all $k$ $\in $ $\Bbb{N}.$ And define the subset $\mathcal{L}_{N}^{o}$ of 
$\mathcal{L}_{N}$ by $\underset{k=1}{\overset{\infty }{\sqcup }}$ $\left( 
\mathcal{L}_{N}^{o}(k)\right) .$
\end{definition}

\strut

Now, we label the elements of a graph groupoid $\Bbb{G}$ by $X_{0}$ (or $%
\mathcal{L}_{N}$). As above, let

\strut

\begin{center}
$X$ $=$ $\{l_{1},$ ..., $l_{N}\}$ and $X_{0}$ $=$ $\{l_{0},$ $l_{1},$ ..., $%
l_{N}\},$
\end{center}

\strut where

\begin{center}
$l_{0}$ $=$ $\overrightarrow{(0,\text{ }0)}$ and $l_{k}$ $=$ $%
\overrightarrow{(0,\text{ }e^{k})},$
\end{center}

\strut

for all $j$ $=$ $1,$ ..., $N$ $=$ $\max \{\deg _{out}(v)$ $:$ $v$ $\in $ $%
V(G)\}.$

\strut

\begin{definition}
Let $G$ be a given graph and let $X$ and $X_{0}$ be the set of our lattices
as above in $\Bbb{R}^{2}$. We call both $X$ and $X_{0}$ $=$ $\{l_{0}\}$ $%
\cup $ $X,$ the (lattice) labeling sets of $G$.
\end{definition}

\strut

Let's fix $v_{0}$ $\in $ $V(G),$ and let $\deg _{out}(v_{0})$ $=$ $n$ $\leq $
$N.$ Define the subset $E_{out}^{v_{0}}$ of the edge set $E(G)$ by

\strut

\begin{center}
$E_{out}^{v_{0}}$ $\overset{def}{=}$ $\{e$ $\in $ $E(G)$ $:$ $e$ $=$ $v_{0}$ 
$e\}$ $\subseteq $ $E(G).$
\end{center}

\strut

i.e., $\deg _{out}(v_{0})$ $=$ $\left| E_{out}^{v_{0}}\right| .$ Then we
decide the weights $\{(v_{0},$ $l_{1}),$ ..., $(v_{0},$ $l_{n})\},$
contained in the set $\{v_{0}\}$ $\times $ $X,$ on $E_{out}^{v_{0}}$. i.e.,

\strut

\begin{center}
$E_{out}^{v_{0}}$ $=$ $\left\{ e_{(v_{0},\text{ }l_{1})},...,e_{(v_{0},\text{
}l_{n})}:\left| 
\begin{array}{c}
\omega _{v_{0}}(e_{(v_{0},\text{ }l_{j})})=(v_{0},\text{ }l_{j}), \\ 
\forall j=1,...,n
\end{array}
\right. \right\} ,$
\end{center}

\strut where

\begin{center}
$\omega _{v_{0}}$ $:$ $E_{out}^{v_{0}}$ $\rightarrow $ $\{v_{0}\}$ $\times $ 
$\{l_{1},$ ..., $l_{n}\}$
\end{center}

\strut

is the \emph{weighting process}, putting the weight $(v_{0},$ $l_{j})$ to an
edge $e_{(v_{0},\text{ }l_{j})}$ $\in $ $E_{out}^{v_{0}}.$

\strut

For any $v$ $\in $ $V(G),$ we can do the same process and we have the
corresponding weighting process $\omega _{v}$ and the corresponding set

$\strut $

\begin{center}
$E_{out}^{v}$ $=$ $\left\{ e_{(v,l_{1})},\text{ ..., }e_{(v,l_{k})}\left| 
\begin{array}{c}
l_{1},\text{ }...,\text{ }l_{k}\in X \\ 
\text{and }\omega _{v}(e_{(v,l_{j})})=(v,\text{ }l_{j}), \\ 
\forall \text{ }j=1,...,k,\text{ where} \\ 
k=\deg _{out}(v)\leq N
\end{array}
\right. \right\} .$
\end{center}

\strut

Notice that if $\deg _{out}(v)$ $=$ $0,$ then $E_{out}^{v}$ is empty, by the
fact that $\{v\}$ $\times $ $\varnothing $ $=$ $\varnothing ,$ where $%
\varnothing $ is the empty set. Notice also that, by the connectedness of $G$%
,

\strut

\begin{center}
$E(G)$ $=$ $\underset{v\in V(G)}{\sqcup }$ $E_{out}^{v},$
\end{center}

\strut

where ``$\sqcup $'' means the disjoint union. So, the each weighting
processes $(\omega _{v})_{v\in V(G)}$ extends to the \emph{weighting process
of }$E(G)$ by $X$:

\strut

\begin{center}
$\omega $ $:$ $E(G)$ $\rightarrow $ $V(G)$ $\times $ $X$
\end{center}

\strut by

\begin{center}
$\omega (e)$ $=$ $\omega (ve)$ $\overset{def}{=}$ $\omega _{v}(e),$
\end{center}

\strut \strut

for all $e$ $=$ $v$ $e$ $\in $ $E(G),$ with $v$ $\in $ $V(G).$ i.e., $\omega 
$ $=$ $\underset{v\in V(G)}{\cup }$ $\omega _{v}.$

\strut

Now, let $G^{-1}$ be the shadow of $G.$ We can determine the similar
weighting process on $G^{-1}.$ Assume that the graph $G$ is weighted by the
above weighting process $\omega $ for the labeling set $X.$ Define now the
sets $-X$ and $-X_{0}$ by

\strut

\begin{center}
$-X$ $\overset{def}{=}$ $\{l_{-1},$ ..., $l_{-N}\}$ and $-X_{0}$ $\overset{%
def}{=}$ $\{l_{0}\}$ $\cup $ $(-X),$
\end{center}

\strut

for all $j$ $=$ $1,$ ..., $N.$ Now, consider the set $V(G^{-1})$ $\times $ $%
(-X)$ $=$ $V(G)$ $\times $ $(-X)$. Then, we can construct the weighting
processes $(\omega _{v}^{-1})_{v\in V(G^{-1})}$ defined by

\strut

\begin{center}
$\omega _{v}^{-1}$ $:$ $E_{out}^{v}$ $\rightarrow $ $\{v\}$ $\times $ $(-X)$
\end{center}

such that

\begin{center}
$\omega _{v}^{-1}(e^{-1})$ $=$ $\omega _{v}^{-1}(e^{-1}v)$ $\overset{def}{=}$
$(v,$ $l_{-i_{j}}),$
\end{center}

\strut

whenever $w_{v}(e)$ $=$ $(v$, $l_{i_{j}})$ $\in $ $\{v\}$ $\times $ $X.$
Therefore, we can define the \emph{weighting process }$\omega ^{-1}$\emph{\
of }$E(G^{-1})$ by

\strut

\begin{center}
$\omega ^{-1}$ $:$ $E(G^{-1})$ $\rightarrow $ $V(G^{-1})$ $\times $ $(-X)$
\end{center}

by

\begin{center}
$\omega ^{-1}(e^{-1})$ $=$ $\omega ^{-1}(e^{-1}v)$ $\overset{def}{=}$ $(v,$ $%
l_{-i_{j}}),$
\end{center}

\strut whenever

\begin{center}
$\omega (e)$ $=$ $(v,$ $l_{i_{j}})$ $\in $ $\{v\}$ $\times $ $X.$
\end{center}

\strut

Define the \emph{canonical projection}

\strut

\begin{center}
$pr_{G}$ $:$ $V(G)$ $\times $ $X$ $\rightarrow $ $X$
\end{center}

by

\begin{center}
$pr_{G}\left( (v,\text{ }l_{k})\right) $ $\overset{def}{=}$ $l_{k},$ for all 
$(v,$ $l_{k})$ $\in $ $V(G)\times X.$
\end{center}

\strut

Similarly, define the canonical projection

\strut

\begin{center}
$pr_{G^{-1}}$ $:$ $V(G)$ $\times $ $(-X)$ $\rightarrow $ $-X$
\end{center}

by

\begin{center}
$pr_{G^{-1}}\left( (v,\text{ }l_{-k})\right) $ $\overset{def}{=}$ $l_{-k},$
for all $(v,$ $l_{-k})$ $\in $ $V(G^{-1})$ $\times $ $(-X).$
\end{center}

\strut \strut

\begin{definition}
Let $G$ be a locally finite connected directed graph equipped with its
weighting process $\omega $ of the edge set $E(G)$ by $V(G)$ $\times $ $X,$
where $X$ $=$ $\{l_{1},$ ..., $l_{N}\},$ and $N$ $=$ $\max \{\deg _{out}(v)$ 
$:$ $v$ $\in $ $V(G)\}.$ Sometimes, we denote the graph $G$ equipped with
the weighting process $\omega $ by the pair $(G,$ $\omega ).$ Then the
weighted graph $(G,$ $\omega )$ is called the canonical weighted graph.
Similarly, the weighted graph $(G^{-1},$ $\omega ^{-1})$ equipped with the
weighting process $\omega ^{-1}$ of $E(G^{-1})$ by $(-X)$ $\times $ $V(G)$
is called the canonical weighted shadow of $(G,$ $\omega ),$ where $G^{-1}$
is the shadow of $G.$

\strut 

Let $(G,$ $\omega )$ (resp. $(G^{-1},$ $\omega ^{-1})$) be the canonical
weighted graph (resp. the canonical weighted shadow of $(G,$ $\omega )$).
The image $pr_{G}$ $\circ $ $\omega $ of $E(G)$ (resp. $pr_{G^{-1}}$ $\circ $
$\omega ^{-1}$), contained in $X$ (resp. in $-X$), is called the canonical
labeling process of $(G,$ $\omega )$ (resp. of $(G^{-1},$ $\omega ^{-1})$).
\end{definition}

\strut

Let $(G,$ $\omega )$ be the canonical weighted graph with its canonical
weighted shadow $(G^{-1},$ $\omega ^{-1}).$ Then we can extend this
weighting processes $\omega $ and $\omega ^{-1}$ to the weighting process,
also denoted by $\omega ,$ the free semigroupoid $\Bbb{F}^{+}(\widehat{G})$
of the shadowed graph $\widehat{G}$ of $G.$ Let $\pm X_{0}^{*}$ be the
collection of all finite words in $\pm X_{0}$ $=$ $X_{0}$ $\cup $ $(-X).$ It
is easy to check that

\strut

\begin{center}
$\pm X_{0}^{*}$ $=$ $\mathcal{L}_{N},$
\end{center}

\strut

where $\mathcal{L}_{N}$ is the lattice path set generated by $\pm X$ $=$ $%
\{l_{\pm 1},$ ..., $l_{\pm N}\},$ including $l_{0}.$

\strut

Define a map

\strut

\begin{center}
$\omega $ $:$ $FP(\widehat{G})$ $\cup $ $\{\emptyset \}$ $\rightarrow $ $%
V(G)^{2}\times $ $X^{*}$
\end{center}

by

\begin{center}
$\omega (\emptyset )$ $\overset{def}{=}$ $\left( \varnothing ,\text{ }%
l_{0}\right) $
\end{center}

and

\begin{center}
$\omega \left( w\right) $ $=$ $\omega \left( v\text{ }e_{1}\text{ ... }%
e_{k}v^{\prime }\right) $ $\overset{def}{=}$ $\left( (v,\text{ }v^{\prime
}),\quad \underset{j=1}{\overset{k}{\Pi }}\widetilde{pr}\left( \omega _{\pm
}(e_{j})\right) \right) ,$
\end{center}

\strut

where $\varnothing $ is the empty set in $V(G)^{2}$ and $FP(\widehat{G})$ is
the (non-reduced) finite path set of the shadowed graph $\widehat{G},$ and
where

\strut

\begin{center}
$\widetilde{pr}$ $\overset{def}{=}pr_{G}$ $\cup $ $pr_{G^{-1}}$ and $\omega
_{\pm }$ $\overset{def}{=}$ $\omega $ $\cup $ $\omega ^{-1},$
\end{center}

whenever

\begin{center}
$w$ $=$ $v$ $w$ $v^{\prime }$ $=$ $e_{1}$ ... $e_{k},$
\end{center}

with

\begin{center}
$v,$ $v^{\prime }$ $\in $ $V(\widehat{G}),$ and $e_{1},$ ..., $e_{k}$ $\in $ 
$E(\widehat{G}),$
\end{center}

\strut

for $k$ $\in $ $\Bbb{N}.$ i.e.,

\begin{center}
$\omega _{\pm }(e)$ $=$ $\left\{ 
\begin{array}{lll}
\omega (e) &  & \text{if }e\in E(G) \\ 
\omega ^{-1}(e) &  & \text{if }e\in E(G^{-1}),
\end{array}
\right. $
\end{center}

and

\begin{center}
$\widetilde{pr}\left( g\right) $ $=$ $\left\{ 
\begin{array}{lll}
pr_{G}(g) &  & \text{if }g\in V(G)\times X \\ 
pr_{G^{-1}}(g) &  & \text{if }g\in V(G)\times (-X).
\end{array}
\right. $
\end{center}

\strut

The above map $\omega $ from $FP(\widehat{G})$ $\cup $ $\{\emptyset \}$ can
be extended to $\Bbb{F}^{+}(\widehat{G})$ into $V(G)^{2}$ $\times $ $\left(
\pm X_{0}^{*}\right) $ by putting the weights on vertices:

\strut

\begin{center}
$\omega (v)$ $\overset{def}{=}$ $\left( (v,\text{ }v),\text{ }l_{0}\right) ,$
for all $v$ $\in $ $V(G).$
\end{center}

\strut \strut \strut

\begin{definition}
The above weighting process $\omega $ of the free semigroupoid $\Bbb{F}^{+}(%
\widehat{G})$ of $\widehat{G}$ is called the \emph{canonical weighting
process of }$\Bbb{F}^{+}(\widehat{G})$ or that of $\Bbb{G}$\emph{\ induced
by the canonical weighted graph} $(G,$ $\omega ).$
\end{definition}

\strut

\begin{remark}
Set-theoretically, the graph groupoid $\Bbb{G}$ is contained in the free
semigroupoid $\Bbb{F}^{+}(\widehat{G})$ of the shadowed graph $\widehat{G}$
of $G.$ So, we can determine the weighting process $\omega $ of $\Bbb{G}$
simply by restricting $\omega $ of $\Bbb{F}^{+}(\widehat{G}).$ i.e., the
weighting process $\omega $ of $\Bbb{G}$ is understood as $\omega $ $=$ $%
\omega $ $\mid _{\Bbb{G}}.$ We will use the notation $\omega $ for the
weighting processes of $\Bbb{F}^{+}(\widehat{G})$ and $\Bbb{G},$
alternatively.
\end{remark}

\strut

Notice that the above weighting process $\omega $ of $\Bbb{F}^{+}(\widehat{G}%
)$ represents the full admissibility conditions of $\Bbb{G}.$ i.e., the
weight $\omega (w)$ of $w$ $\in $ $\Bbb{G}$ contains its initial and
terminal vertices and the connecting pattern (admissibility of edges) of $w,$
in terms of labelings in $\pm X_{0}$ $=$ $X_{0}$ $\cup $ $\left( -X\right) .$

\strut \strut

\begin{definition}
Let $(G,$ $\omega )$ be the canonical weighted graph with its canonical
weighted shadow $(G^{-1},$ $\omega ^{-1}).$ Let $\Bbb{G}$ be the graph
groupoid of $G$ and assume that $\omega $ is the canonical weighting process
of $\Bbb{G}$ induced by $(G,$ $\omega ).$ Then we call the pair $(\Bbb{G},$ $%
\omega )$ the labeled graph groupoid of $(G,$ $\omega ).$ We call the
processes $\widetilde{pr}$ $\circ $ $\omega $ on $(\Bbb{G},$ $\omega ),$ the
labeling process, where $\widetilde{pr}$ $=$ $pr_{G}$ $\cup $ $pr_{G^{-1}}$.
\end{definition}

\strut

\textbf{Assumption} From now, if we mention (locally finite connected)
directed graphs, then they are automatically assumed to be the canonical
weighted graphs. Similarly, if we mention graph groupoids, then they are
also automatically assumed to be the labeled graph groupoids. $\square $

\strut

\begin{example}
Let $G$ be the circulant graph with three vertices with

\strut 

\begin{center}
$V(G)$ $=$ $\{v_{1},$ $v_{2},$ $v_{3}\}$
\end{center}

and

\begin{center}
$E(G)$ $=$ $\{e_{j}$ $=$ $v_{j}$ $e_{j}$ $v_{j+1}$ $:$ $j$ $=$ $1,$ $2,$ $3,$
with $v_{4}$ $\overset{def}{=}$ $v_{1}\}.$
\end{center}

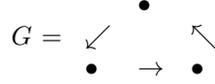
\strut \strut i.e.,

\begin{center}
$G$ $=$ $
\begin{array}{lll}
& \bullet  &  \\ 
\swarrow  &  & \nwarrow  \\ 
\bullet  & \rightarrow  & \bullet 
\end{array}
$
\end{center}

\strut 

Then, since $\deg _{out}(v_{j})$ $=$ $1,$ for all $j$ $=$ $1,$ $2,$ $3,$ we
can get the labeling sets

\strut 

\begin{center}
$X$ $=$ $\{l_{1}\}$ $=$ $\left\{ \overrightarrow{(1,\text{ }e)}\right\} $
and $X_{0}$ $=$ $\{l_{0},$ $l_{1}\},$
\end{center}

\strut 

contained in the lattice path set $\mathcal{L}_{1}.$ Relatively, we have that

\strut 

\begin{center}
$-X$ $=$ $\{l_{-1}\}$ $=$ $\left\{ \overrightarrow{(1,\text{ }-e)}\right\} ,$
\end{center}

\strut 

for the shadow $G^{-1}$ of $G.$ Then we can get the weights

\strut 

\begin{center}
$\omega (e_{j})$ $=$ $\left( (v_{j},\text{ }v_{j+1}),\text{ }l_{1}\right) ,$
for all $j$ $=$ $1,$ $2,$ $3$
\end{center}

and

\begin{center}
$\omega ^{-1}(e_{j}^{-1})$ $=$ $\left( (v_{j+1},\text{ }v_{j}),\text{ }%
l_{-1}\right) ,$ for all $j$ $=$ $1,$ $2,$ $3,$
\end{center}

\strut 

where $v_{4}$ $\overset{def}{=}$ $v_{1},$ in $V(G).$ We can construct the
labeled graph groupoid $(\Bbb{G},$ $\omega ).$ Let $w$ $=$ $e_{2}$ $e_{3}$ $%
e_{1}$ $\in $ $FP_{r}(\widehat{G}).$ Then we can get that

\strut 

\begin{center}
$\omega (w)$ $=$ $((v_{2},$ $v_{2}),$ $l_{1}l_{1}l_{1}),$ and $\omega
(w^{-1})$ $=$ $((v_{2},$ $v_{2}),$ $l_{-1}l_{-1}l_{-1}).$
\end{center}

\strut 

Take now the reduced finite path $y$ $=$ $e_{1}^{-1}e_{3}^{-1}$ $\in $ $%
FP_{r}(\widehat{G}).$ Then we can get

\strut 

\begin{center}
$\omega (y)$ $=$ $((v_{2},$ $v_{3}),$ $l_{-1}l_{-1}).$

\strut 

Suppose we have $w^{2}$ $=$ $e_{2}$ $e_{3}$ $e_{1}$ $e_{2}$ $e_{3}$ $e_{1}$
in $FP_{r}(\widehat{G}).$ Then the weight of $w^{2}$ is

\strut 

$\omega (w^{2})$ $=$ $((v_{2},$ $v_{2}),$ $l_{1}l_{1}l_{1}l_{1}l_{1}l_{1}).$
\end{center}
\end{example}

\strut

\strut

\strut

\subsection{The Operation $\theta $ and $\omega _{+}$}

\strut

\strut

In this section, we will define some operations on the labeling set $\pm
X_{0}$ $=$ $X_{0}$ $\cup $ $(-X)$ of the labeled graph groupoid $(\Bbb{G},$ $%
\omega ),$ induced by the canonical weighted graph $(G,$ $\omega ).$ Then
Define an operation

$\strut $

\begin{center}
$\theta $ $:$ $\pm X_{0}^{*}$ $=$ $\mathcal{L}_{N}$ $\rightarrow $ $\Bbb{Z}$
\end{center}

by

\strut

\begin{center}
$\theta \left( l_{i_{1}}l_{i_{2}}\text{ ... }l_{i_{n}}\right) $ $\overset{def%
}{=}$ $\sum_{j=1}^{n}$ $i_{j},$
\end{center}

\strut \strut

for all $l_{i_{1}}$ ... $l_{i_{n}}$ $\in $ $\pm X_{0}^{*}$ $=$ $\mathcal{L}%
_{N},$ for $n$ $\in $ $\Bbb{N},$ where the addition of the right-hand side
means the vector addition on $\Bbb{R}^{2}.$

\strut

By help of $\theta ,$ we can define an operation,

$\strut $

\begin{center}
$\omega _{+}$ $:$ $V(G)^{2}$ $\times $ $\left( \pm X_{0}^{*}\right) $ $%
\rightarrow $ $V(G)^{2}$ $\times $ $\Bbb{Z}$
\end{center}

by

\begin{center}
$
\begin{array}{ll}
\omega _{+}\left( (v,\text{ }v^{\prime }),\text{ }l_{i_{1}}\text{ ... }%
l_{i_{n}}\right) & \overset{def}{=}\left( (v,\text{ }v^{\prime }),\text{ }%
\theta \left( l_{i_{1}}\text{ ... }l_{i_{n}})\right) \right) \\ 
&  \\ 
& =\left( (v,\text{ }v^{\prime }),\text{ }\sum_{j=1}^{n}i_{j}\right) ,
\end{array}
$
\end{center}

\strut

for all $((v,$ $v^{\prime }),$ $l_{i_{1}}$ ... $l_{i_{n}})$ $\in $ $V(G)^{2}$
$\times $ $(\pm X_{0}^{*})$. Later, we will use the above operation $\omega
_{+}$ to verify the amalgamated moments of the labeling operators on right
graph von Neumann algebras. In conclusion, we can get that:

\strut

\begin{proposition}
Let $w_{1}$ and $w_{2}$ be the reduced finite paths in $FP_{r}(\widehat{G})$
given in the previous paragraph. Then

\strut 

\begin{center}
$\theta \left( \left( \widetilde{pr}\text{ }\circ \text{ }\omega \right)
(w_{1}w_{2})\right) $ $=$ $0$ in $\Bbb{Z}$,
\end{center}

\strut and

\begin{center}
$\omega _{+}\left( \omega (w_{1}w_{2})\right) $ $=$ $\left( (v,\text{ }v),%
\text{ }0\right) ,$
\end{center}

\strut 

if and only if the element $w_{1}$ $w_{2}$ $=$ $v$ in $\Bbb{G}$ $\setminus $ 
$\{\emptyset \}.$
\end{proposition}

\strut

\begin{proof}
($\Rightarrow $) Suppose $\omega _{+}\left( \omega (w_{1}w_{2})\right) $ $=$ 
$\left( (v,\text{ }v),\text{ }0\right) .$ By hypothesis, we have that

\strut

\begin{center}
$\omega (w_{1}w_{2})$ $=$ $\left( (v,\text{ }v),\text{ }(l_{i_{1}}\text{ ... 
}l_{i_{n}}\text{ }l_{j_{1}}\text{ ... }l_{j_{k}})\right) ,$
\end{center}

whenever

\begin{center}
$\omega (w_{1})$ $=$ $\left( (v,\text{ }v^{\prime }),\text{ }l_{i_{1}}\text{
... }l_{i_{n}}\right) $
\end{center}

and

\begin{center}
$\omega (w_{2})$ $=$ $\left( (v^{\prime },\text{ }v),\text{ }l_{j_{1}}\text{
... }l_{j_{k}}\right) .$
\end{center}

\strut

By assumption, $\theta \left( l_{i_{1}}\text{ ... }l_{i_{n}}\text{ }%
l_{j_{1}}l\text{ ... }l_{j_{k}})\right) $ $=$ $0.$ And since $w_{1}$ and $%
w_{2}$ are ``reduced'' finite paths in $\Bbb{G},$ we can verify that (i) $n$ 
$=$ $k,$ and (ii) for any $l_{i_{p}}$, there should be the unique entry $%
l_{j_{q}}$ $=$ $l_{-i_{p}},$ for $p,$ $q$ $=$ $1,$ ..., $n$ $=$ $k.$ This
shows that, if $w_{1}$ $=$ $e_{i_{1}}$ ... $e_{i_{n}}$ and $w_{2}$ $=$ $%
e_{j_{1}}$ ... $e_{j_{n}}$ in $FP_{r}(\widehat{G}),$ with $e_{i_{1}},$ ..., $%
e_{i_{n}},$ $e_{j_{1}},$ ..., $e_{j_{n}}$ $\in $ $E(\widehat{G}),$ then $%
e_{j_{q}}$ $=$ $e_{i_{p}}^{-1},$ uniquely, for $p,$ $q$ $=$ $1,$ ..., $n.$
Therefore, $w_{1}$ $w_{2}$ $=$ $v.$

\strut

($\Leftarrow $)\strut \ Suppose now $w_{1}$ $w_{2}$ $=$ $v$ $\in $ $V(G).$
Then

\strut

\begin{center}
$\omega _{+}\left( \omega (w_{1}w_{2})\right) $ $=$ $\omega _{+}\left(
\omega (v)\right) $ $=$ $\omega _{+}\left( (v,\text{ }v),\text{ }%
l_{0}\right) $ $=$ $\left( (v,\text{ }v),\text{ }0\right) .$
\end{center}

\strut
\end{proof}

\strut

The above proposition shows how the labeling on the graph groupoid $\Bbb{G}$
works. The labeling process on $\Bbb{G}$ gives a way to check how
admissibility works under (RR).

\strut

\strut

\subsection{Graph Automata}

\strut

\strut

In this section, we will construct an automaton induced by the given
canonical weighted graph $G$ $=$ $(G,$ $\omega ).$ Let $X_{0}$ $=$ $\{l_{0},$
$l_{1},$ ..., $l_{N}\}$ be the labeling set of the countable locally finite
connected directed graph $G,$ consisting of lattices in $\Bbb{R}^{2},$ where

\strut

\begin{center}
$N$ $\overset{def}{=}$ $\max \{\deg _{out}(v)$ $:$ $v$ $\in $ $V(G)\},$
\end{center}

\strut \strut

and let $\pm X_{0}$ $=$ $X_{0}$ $\cup $ $(-X)$ $=$ $\{l_{0},$ $\pm l_{1},$
..., $\pm l_{N}\}.$ Recall that the collection $\pm X_{0}^{*}$ of all words
in $\pm X_{0}$ is identical to the lattice path set $\mathcal{L}_{N}.$

\strut

Under the above setting, we can create the automaton $\mathcal{A}_{G}$ $=$ $<%
\mathcal{X}_{0},$ $E(\widehat{G}),$ $\varphi ,$ $\psi >$, where $\widehat{G}$
is the shadowed graph of $G$ and

\strut

\begin{center}
$\mathcal{X}_{0}$ $\overset{def}{=}$ $\{\emptyset _{G}\}$ $\cup $ $\left(
V(G)^{2}\times (\pm X_{0})\right) ,$
\end{center}

\strut and

\begin{center}
$\varphi \left( \left( (v_{1},\text{ }v_{2}),\text{ }l\right) ,\text{ }%
e\right) $ $\overset{def}{=}$ $\left\{ 
\begin{array}{ll}
\begin{array}{l}
\omega (e)
\end{array}
& 
\begin{array}{l}
\begin{array}{l}
\text{if }e=v_{2}\text{ }e\text{ in }E(\widehat{G}),\text{ equivalently,} \\ 
\exists \text{ }v^{\prime }\in V(G)\text{ and }l^{\prime }\in \pm X_{0}\text{
s.t.,} \\ 
\quad \omega (e)=\left( (v_{2},\text{ }v^{\prime }),\text{ }l^{\prime
}\right) \neq \emptyset _{G}
\end{array}
\end{array}
\\ 
\,\,\,\,\emptyset _{G} & \quad \text{otherwise,}
\end{array}
\right. $
\end{center}

and

\begin{center}
$\psi \left( \left( (v_{1},\text{ }v_{2}),\text{ }t\right) ,\text{ }e\right) 
$ $\overset{def}{=}$ $\left\{ 
\begin{array}{ll}
e & \text{if }e=v_{2}e\text{ in }E(\widehat{G}) \\ 
\emptyset & \quad \text{otherwise,}
\end{array}
\right. $
\end{center}

\strut

for all $e$ $\in $ $E(\widehat{G}),$ where $\emptyset _{G}$ $=$ $%
(\varnothing ,$ $l_{0})$ is the \emph{empty element of }$\mathcal{X}_{0}$
and $\emptyset $ is the empty element in $\Bbb{F}^{+}(\widehat{G}).$

\strut

\begin{definition}
We will say that the automaton $\mathcal{A}_{G}$ $=$ $<\mathcal{X}_{0},$ $E(%
\widehat{G}),$ $\varphi ,$ $\psi >$ is the graph automaton induced by the
canonical weighted graph $G$ $=$ $(G,$ $\omega ).$ Sometimes, we call $%
\varphi $ and $\psi ,$ the labeling map and the shift (or shifting map),
respectively.
\end{definition}

\strut

We can realize that the graph automaton $\mathcal{A}_{G}$ induced by $G$ is
identically same with the automaton $\mathcal{A}_{G^{-1}}$ induced by the
shadow $G^{-1}$ of $G.$ Now, fix an edge $e_{0}$ $\in $ $E(\widehat{G})$ and
a reduced finite path $w$ $=$ $v_{1}e_{1}$ ... $e_{n}v_{n}^{\prime }$ $\in $ 
$FP_{r}(\widehat{G}),$ with $e_{j}$ $=$ $v_{j}$ $e_{j}$ $v_{j}^{\prime }$ $%
\in $ $E(\widehat{G}),$ where $v_{j},$ $v_{j}^{\prime }$ $\in $ $V(\widehat{G%
}),$ for all $j$ $=$ $1,$ ..., $n,$ for some $n$ $\in $ $\Bbb{N}.$ Then we
can extend $\varphi $ and $\psi $ inductively on the free semigroupoid $\Bbb{%
F}^{+}(\widehat{G})$ of $\widehat{G}$ as follows:

\strut

$\qquad \varphi \left( \left( (v_{1},\text{ }v_{n}^{\prime }),\text{ }\omega
(e_{1})\text{ ... }\omega (e_{n})\right) ,\text{ }e_{0}\right) $

\strut

$\qquad \qquad \qquad =$ $\varphi \left( \omega (w),\text{ }e_{0}\right) $

\strut

$\qquad \qquad \qquad \overset{def}{=}$ $\left\{ 
\begin{array}{ll}
\omega (e_{0}) & \text{if }w\text{ }e_{0}\neq \emptyset \\ 
\emptyset _{G} & \text{otherwise,}
\end{array}
\right. $

\strut

and similarly,

\strut

$\qquad \psi \left( \left( (v_{1},\text{ }v_{n}^{\prime }),\text{ }\omega
(e_{1})\text{ ... }\omega (e_{n})\right) ,\text{ }e_{0}\right) $

\strut

$\qquad \qquad \qquad =$ $\psi \left( \omega (w),\text{ }e_{0}\right) $

\strut

$\qquad \qquad \qquad \overset{def}{=}$ $\left\{ 
\begin{array}{ll}
e_{0} & \text{if }w\text{ }e_{0}\neq \emptyset \\ 
\emptyset & \text{otherwise,}
\end{array}
\right. $

\strut $\strut $

for $n$ $\in $ $\Bbb{N}.$ Also, inductively, we have that, if $w^{\prime }$ $%
=$ $e_{1}^{\prime }$ ... $e_{m}^{\prime }$ $\in $ $FP_{r}(\widehat{G}),$ for 
$m$ $\in $ $\Bbb{N},$ then

\strut

$\qquad \varphi \left( \left( (v_{1},\text{ }v_{n}^{\prime }),\text{ }\omega
(e_{1})\text{ ... }\omega (e_{n})\right) ,\text{ }w^{\prime }\right) $

\strut

$\qquad \qquad \qquad \overset{def}{=}$ $\left\{ 
\begin{array}{ll}
\omega (e_{m}^{\prime }) & \text{if }w\text{ }w^{\prime }\neq \emptyset \\ 
\emptyset _{G} & \text{otherwise,}
\end{array}
\right. $

\strut

and similarly,

\strut

$\qquad \psi \left( \left( (v_{1},\text{ }v_{n}^{\prime }),\text{ }\omega
(e_{1})\text{ ... }\omega (e_{n})\right) ,\text{ }w^{\prime }\right) $

\strut

$\qquad \qquad \qquad \overset{def}{=}$ $\left\{ 
\begin{array}{ll}
w^{\prime } & \text{if }w\text{ }w^{\prime }\neq \emptyset \\ 
\emptyset & \text{otherwise.}
\end{array}
\right. $

\strut \strut \strut

\begin{remark}
As we have seen in the previous paragraphs, the (extended) definitions of $%
\varphi $ and $\psi $ from $\Bbb{F}^{+}(\widehat{G})$ into $V(G)^{2}$ $%
\times $ $(\pm X_{0}^{*})$ have the following meaning: by definition, we
have $\varphi (\omega (w),$ $w^{\prime })$ $=$ $\omega (e^{\prime }),$
whenever $w$ and $w^{\prime }$ $=$ $e^{\prime }$ ... $e^{\prime \prime }$
are admissible, where $\omega (w)$ $\in $ $V(G)^{2}$ $\times $ $(\pm
X_{0}^{*})$ and $e^{\prime },$ ..., $e^{\prime \prime }$ $\in $ $E(\widehat{G%
}),$ and hence it represents the admissibility of $w$ and $w^{\prime },$ in
terms of the weighting process $\omega .$ Equivalently,

\strut 

\begin{center}
Computing $\varphi (\omega (w),$ $w^{\prime })$ means: if $w$ and $w^{\prime
}$ are admissible, what is the weight of the connecting edge $e^{\prime }$
of $w^{\prime }$?
\end{center}

\strut 

Under the same setting, we have $\psi (\omega (w),$ $w^{\prime })$ $=$ $%
e^{\prime }$: it again represents the admissibility of $w$ and $w^{\prime },$
in terms of the edges. Thus, roughly speaking, the labeling map $\varphi $
and the shifting $\psi $ represent the admissibility conditions on $\Bbb{F}%
^{+}(\widehat{G})$ (and hence $\Bbb{G}$). It gives more detailed information
about $\Bbb{F}^{+}(\widehat{G})$ (or $\Bbb{G}$). Equivalently,

\strut 

\begin{center}
Computing $\psi (\omega (w),$ $w^{\prime })$ means: if $w$ and $w^{\prime }$
are admissible, what is the starting edge of $w^{\prime }$?
\end{center}
\end{remark}

\strut

Construct the \emph{automata actions }$\{\mathcal{A}_{e}$\emph{\ }$:$\emph{\ 
}$e$\emph{\ }$\in $\emph{\ }$E(\widehat{G})\}$\emph{\ of} $\mathcal{A}_{G},$
acting on the set $V(G)^{2}$ $\times $ $(\pm X_{0}^{*})$ (consisting of all
finite words in $\mathcal{X}_{0}$). Indeed, we can define the action $%
\mathcal{A}_{e}$ by

\strut

\begin{center}
$\mathcal{A}_{e}\left( \left( (v_{1},\text{ }v_{2}),\text{ }l_{i_{1}}\text{
... }l_{i_{n}}\right) \right) $ $\overset{def}{=}$ $\varphi \left( \left(
(v_{1},\text{ }v_{2}),\text{ }l_{i_{1}}\text{ ... }l_{i_{n}}\right) ,\text{ }%
e\right) $
\end{center}

\strut

for all $\left( (v_{1},v_{2}),\text{ }l_{i_{1}}\text{ ... }l_{i_{n}}\right) $
$\in $ $V(G)^{2}$ $\times $ $\left( \pm X_{0}^{*}\right) ,$ for all $e$ $\in 
$ $E(\widehat{G}).$ It is easy to be checked that the actions $\mathcal{A}%
_{e}$'s satisfies

\strut

\begin{center}
$\mathcal{A}_{e_{1}}$ $\mathcal{A}_{e_{2}}$ $=$ $\mathcal{A}_{e_{2}e_{1}},$
for all $e_{1},$ $e_{2}$ $\in $ $E(\widehat{G}).$
\end{center}

\strut

Thus the actions $\{\mathcal{A}_{v}$ $:$ $v$ $\in $ $V(\widehat{G})\}$ are
induced by $\mathcal{A}_{e}$'s since

\strut

\begin{center}
$\mathcal{A}_{v}$ $\overset{def}{=}$ $\mathcal{A}_{ee^{-1}}$ $=$ $\mathcal{A}%
_{e^{-1}}$ $\mathcal{A}_{e},$ whenever $e$ $=$ $v$ $e$ $\in $ $E(\widehat{G}%
).$
\end{center}

\strut

This guarantees that the automata actions $\{\mathcal{A}_{e}$ $:$ $e$ $\in $ 
$E(\widehat{G})\}$ on $V(G)^{2}$ $\times $ $\pm X_{0}^{*}$ generate the
groupoid which is groupoid-isomorphic to the graph groupoid $\Bbb{G}$ of $G.$
i.e., we can have that:

\strut

\begin{theorem}
The set $\{\mathcal{A}_{e}$ $:$ $e$ $\in $ $E(\widehat{G})\}$ of automata
actions, acting on $V(G)$ $\times $ $(\pm X_{0}^{*}),$ generates the actions 
$\{\mathcal{A}_{w}$ $:$ $w$ $\in $ $\Bbb{F}^{+}(\widehat{G})\}$ of $\mathcal{%
A}_{G},$ acting on the same set. Moreover, the groupoid generated by this
set $\{\mathcal{A}_{e}$ $:$ $e$ $\in $ $E(\widehat{G})\}$ is
groupoid-isomorphic to the graph groupoid $\Bbb{G}$ of $G.$ $\square $
\end{theorem}

\strut

\strut

\strut

\subsection{Graph-Automata Trees}

\strut

\strut

As before, the standing assumption for our graphs $G$ are as follows:
locally finite, connected, and countable. The edges of $G$ are directed, and
they are assigned weights with the use of a labeling set. In this section,
we will consider graph-automata trees where the labeled graph groupoids act
on. Let $\Bbb{G}$ be the labeled graph groupoid induced by the canonical
weighted graph $G$ $=$ $(G,$ $\omega ),$ with the labeling set $\pm X_{0}$ $%
= $ $\{l_{0},$ $\pm l_{1},$ ..., $\pm l_{N}\},$ and let the canonical
weighting process $\omega $ of the free semigroupoid $\Bbb{F}^{+}(\widehat{G}%
)$ of the shadowed graph $\widehat{G}$ of $G,$ be given as before, where

$\strut $

\begin{center}
$N$ $\overset{def}{=}$ $\max \{\deg _{out}(v)$ $:$ $v$ $\in $ $V(G)\}.$
\end{center}

\strut

The canonical weighted graph $G$ creates the corresponding graph automaton

$\strut $

\begin{center}
$\mathcal{A}_{G}$ $=$ $<\mathcal{X}_{0},$ $E(\widehat{G}),$ $\varphi ,$ $%
\psi >,$
\end{center}

where

\begin{center}
$\mathcal{X}_{0}$ $=$ $\{\emptyset _{*}\}$ $\cup $ $\left( V(G)^{2}\times
(\pm X_{0})\right) ,$
\end{center}

\strut

and $\varphi $ is the labeling map and $\psi $ is the shifting map. We
already observed that $\Bbb{G}$ acts on a $E(\widehat{G})$-set $V(G)^{2}$ $%
\times $ $(\pm X_{0}^{*}).$ For convenience, we will denote the set $%
V(G)^{2} $ $\times $ $(\pm X_{0}^{*})$ by $\mathcal{X}_{0}^{*}.$

\strut

Recall that two countable directed graphs $G_{1}$ and $G_{2}$ are \emph{%
graph-isomorphic}, if there exists a bijection

$\strut $

\begin{center}
$g$ $:$ $V(G_{1})$ $\cup $ $E(G_{1})$ $\rightarrow $ $V(G_{2})$ $\cup $ $%
E(G_{2}),$
\end{center}

\strut

such that (i) $g\left( V(G_{1})\right) $ $=$ $V(G_{2})$ and $g\left(
E(G_{1})\right) $ $=$ $E(G_{2}),$ and (ii) $g\left( e\right) $ $=$ $g(v_{1}$ 
$e$ $v_{2})$ $=$ $g(v_{1})$ $g(e)$ $g(v_{2})$ in $E(G_{2}),$ whenever $e$ $=$
$v_{1}$ $e$ $v_{2}$ $\in $ $E(G_{1}),$ with $v_{1},$ $v_{2}$ $\in $ $%
V(G_{1}).$

\strut

In [10] and [11], we showed that if two graphs $G_{1}$ and $G_{2}$ are
graph-isomorphic, then the corresponding graph groupoids $\Bbb{G}_{1}$ and $%
\Bbb{G}_{1}$ are groupoid-isomorphic. More generally, if two graphs $G_{1}$
and $G_{2}$ have the graph-isomorphic shadowed graphs $\widehat{G_{1}}$ and $%
\widehat{G_{2}}$, then the graph groupoids $\Bbb{G}_{1}$ and $\Bbb{G}_{1}$
are groupoid-isomorphic. Also, we showed that if two graph groupoids $\Bbb{G}%
_{1}$ and $\Bbb{G}_{2}$ are groupoid-isomorphic, then the (left) graph von
Neumann algebras $\overline{\Bbb{C}[\Bbb{G}_{1}]}^{w}$ and $\overline{\Bbb{C}%
[\Bbb{G}_{2}]}^{w}$ are $*$-isomorphic, as $W^{*}$-subalgebras in the
operator algebra $B(H_{G_{1}})$ $=$ $B(H_{G_{2}}).$ Therefore, if two graphs 
$G_{1}$ and $G_{2}$ have the graph-isomorphic shadowed graphs, then the
right graph von Neumann algebras $M_{G_{1}}$ and $M_{G_{2}}$ are $*$%
-isomorphic, too.

\strut

\begin{proposition}
Let $G_{1}$ and $G_{2}$ be countable directed graphs. If the shadowed graphs 
$\widehat{G_{1}}$ and $\widehat{G_{2}}$ are graph-isomorphic, then the right
graph von Neumann algebras $M_{G_{1}}$ and $M_{G_{2}}$ are $*$-isomorphic. $%
\square $
\end{proposition}

\strut

Recall also that we say a directed graph $G$ is a (directed)\emph{\ tree},
if this graph $G$ is connected and it has no loop finite paths in the free
semigroupoid $\Bbb{F}^{+}(G)$ of $G.$ Also, we say that a directed tree is 
\emph{rooted}, if we can find-and-fix a vertex $v_{0}$ of $G,$ with $\deg
_{in}(v_{0})$ $=$ $0.$ This fixed vertex $v_{0}$ is called the\emph{\ root of%
} $G.$ For instance, a graph

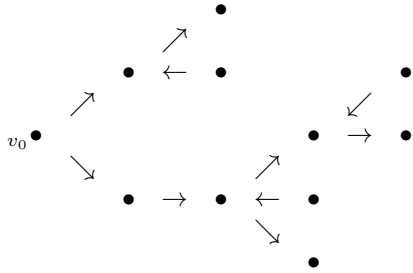
\strut

\begin{center}
$
\begin{array}{lllllllll}
&  &  &  & \bullet &  &  &  &  \\ 
&  &  & \nearrow &  &  &  &  &  \\ 
&  & \bullet & \leftarrow & \bullet &  &  &  & \bullet \\ 
& \nearrow &  &  &  &  &  & \swarrow &  \\ 
_{v_{0}}\bullet &  &  &  &  &  & \bullet & \rightarrow & \bullet \\ 
& \searrow &  &  &  & \nearrow &  &  &  \\ 
&  & \bullet & \rightarrow & \bullet & \leftarrow & \bullet &  &  \\ 
&  &  &  &  & \searrow &  &  &  \\ 
&  &  &  &  &  & \bullet &  & 
\end{array}
$
\end{center}

\strut

is a rooted tree with its (fixed) root $v_{0}$. A rooted tree $G$ is \emph{%
one-flow}, if the directions of edges oriented only one way from the root $%
v_{0}.$ For example, a graph

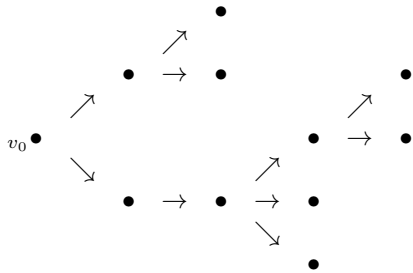
\strut

\begin{center}
$
\begin{array}{lllllllll}
&  &  &  & \bullet &  &  &  &  \\ 
&  &  & \nearrow &  &  &  &  &  \\ 
&  & \bullet & \rightarrow & \bullet &  &  &  & \bullet \\ 
& \nearrow &  &  &  &  &  & \nearrow &  \\ 
_{v_{0}}\bullet &  &  &  &  &  & \bullet & \rightarrow & \bullet \\ 
& \searrow &  &  &  & \nearrow &  &  &  \\ 
&  & \bullet & \rightarrow & \bullet & \rightarrow & \bullet &  &  \\ 
&  &  &  &  & \searrow &  &  &  \\ 
&  &  &  &  &  & \bullet &  & 
\end{array}
$
\end{center}

\strut

is a one-flow rooted tree with its root $v_{0}$. A one-flow rooted tree $G$
is said to be \emph{growing}, if $G$ is an infinitely countable directed
graph. Finally, we will say that a one-flow growing rooted tree $G$ is \emph{%
regular}, if the out-degrees of all vertices are identical. For example, a
graph

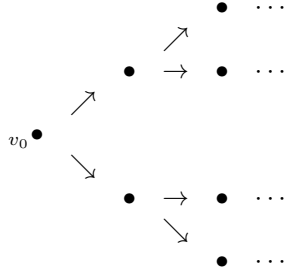
\strut

\begin{center}
$
\begin{array}{llllll}
&  &  &  &  &  \\ 
&  &  &  & \bullet & \cdots \\ 
&  &  & \nearrow &  &  \\ 
&  & \bullet & \rightarrow & \bullet & \cdots \\ 
& \nearrow &  &  &  &  \\ 
_{v_{0}}\bullet &  &  &  &  &  \\ 
& \searrow &  &  &  &  \\ 
&  & \bullet & \rightarrow & \bullet & \cdots \\ 
&  &  & \searrow &  &  \\ 
&  &  &  & \bullet & \cdots
\end{array}
$
\end{center}

\strut

is a regular one-flow growing rooted tree. In particular, if $\deg _{out}(v)$
$=$ $N,$ for all vertices $v$, then this regular one-flow growing rooted
tree is called \emph{the} $N$-\emph{regular tree}. The very above example is
the $2$-regular tree.

\strut

Let $N$ be the maximal out-degree of the graph $G,$ and let $\mathcal{T}%
_{2N} $ be the $2N$-regular tree. Then the automata actions $\{\mathcal{A}%
_{w}$ $:$ $w$ $\in $ $\Bbb{F}^{+}(\widehat{G})\}$ of the graph automaton $%
\mathcal{A}_{G}$ acts on this $2N$-regular tree $\mathcal{T}_{2N}.$ Indeed,
for the set $\{\mathcal{A}_{w}\}_{w\in \Bbb{F}^{+}(\widehat{G})}$ of
automata actions of the graph automaton $\mathcal{A}_{G},$ we can create a
one-flow growing rooted tree $\mathcal{T}_{G}$ having its ``arbitrarily''
fixed root $\omega (v_{0})$ $=$ $\left( (v_{0},v_{0}),\text{ }x_{0}\right) $ 
$\in $ $\mathcal{X}_{0}^{*}$, where

\strut

\begin{center}
$V(\mathcal{T}_{G})$ $\overset{def}{=}$ $\left\{ \varphi (X,\text{ }e)\left|
X\in \mathcal{X}_{0}^{*},\text{ }e\in E(\widehat{G})\right. \right\} $
\end{center}

and

\begin{center}
$E(\mathcal{T}_{G})$ $\overset{def}{=}$ $\left\{ \psi (X,\text{ }e)\left|
X\in \mathcal{X}_{0}^{*},\text{ }e\in E(\widehat{G})\right. \right\} ,$
\end{center}

\strut

where $\widehat{G}$ is the shadowed graph of $G.$

\strut

\begin{remark}
Notice that, by the connectedness of the shadowed graph $\widehat{G}$ of the
connected graph $G,$ we can fix any weight $\omega (v)$ of $v$ $\in $ $V(%
\widehat{G}),$ as the root of the tree $\mathcal{T}_{G}$. Suppose $\mathcal{T%
}_{v_{1}}$ and $\mathcal{T}_{v_{2}}$ are the one-flow growing trees with
their roots $((v_{1},$ $v_{1}),$ $x_{0})$ and $((v_{2},$ $v_{2}),$ $x_{0}),$
respectively. Then, $\mathcal{T}_{v_{i}}$ is embedded in $\mathcal{T}%
_{v_{j}},$ as a full-subgraph (See below), whenever $i$ $\neq $ $j$ $\in $ $%
\{1,$ $2\}.$ In general, the graphs $\mathcal{T}_{v_{1}}$ and $\mathcal{T}%
_{v_{2}}$ have no graph-isomorphic relation, but they are embedded from each
other, by the connectedness of $\widehat{G}.$ Since $v_{1}$ and $v_{2}$ are
arbitrary, we can consider only one choice $\mathcal{T}_{v},$ for $v$ $\in $ 
$V(G),$ as a candidate of the tree, where $\{\mathcal{A}_{w}\}_{w\in \Bbb{F}%
^{+}(\widehat{G})}$ act. Denote it by $\mathcal{T}_{G}$. i.e., whenever we
choose one tree $\mathcal{T}_{v},$ for $v$ $\in $ $V(G),$ then the trees $%
\mathcal{T}_{v^{\prime }}$'s are embedded in $\mathcal{T}_{v},$ for all $%
v^{\prime }$ $\in $ $V(G).$ Without loss of generality, if we write $%
\mathcal{T}_{G}$ from now, then it means a tree $\mathcal{T}_{v},$ for a
fixed vertex $v$ $\in $ $V(G).$ Remark that, the tree $\mathcal{T}_{G}$ has
its root $\omega (v),$ if and only if $\mathcal{T}_{G}$ $=$ $\mathcal{T}_{v}.
$
\end{remark}

\strut

We can easily check that

\strut

\begin{center}
$FP(\mathcal{T}_{G})$ $=$ $\left\{ \psi (X,\text{ }w)\left| X\in \mathcal{X}%
_{0}^{*},\text{ }w\in FP(\widehat{G})\right. \right\} ,$
\end{center}

and hence

\begin{center}
$
\begin{array}{ll}
\Bbb{F}^{+}(\mathcal{T}_{G})= & \{\emptyset _{G}\}\cup \left\{ \varphi (X,%
\text{ }e)\left| 
\begin{array}{c}
X\in \mathcal{X}_{0}^{*},\text{ } \\ 
e\in E(\widehat{G})
\end{array}
\right. \right\} \\ 
&  \\ 
& \quad \quad \;\cup \left\{ \psi (X,\text{ }w)\left| 
\begin{array}{c}
X\in \mathcal{X}_{0}^{*}, \\ 
w\in FP(\widehat{G})
\end{array}
\right. \right\} .
\end{array}
$
\end{center}

\strut \strut \strut

By the connectedness of the graph $G,$ and by the definition of the automata
actions, every (nonempty) finite paths in $FP(\widehat{G})$ is embedded in
the tree $\mathcal{T}_{G},$ via the automata actions. Then, we can construct
the full-subgraphs $\{\mathcal{T}_{w}$ $:$ $w$ $\in $ $FP(\widehat{G})\}$ of 
$\mathcal{T}_{G},$ where $\mathcal{T}_{w}$'s are the one-flow growing rooted
tree with their roots $\omega \left( w\right) ,$ for all $w$ $\in $ $FP(%
\widehat{G}).$

\strut

Recall that, we say that a countable directed graph $G_{1}$ is a \emph{%
full-subgraph} of a countable directed graph $G_{2},$ if

\strut

\begin{center}
$E(G_{1})$ $\subseteq $ $E(G_{2})$
\end{center}

and

\begin{center}
$V(G_{1})$ $=$ $\left\{ v\in V(G_{2})\left| 
\begin{array}{c}
e=ve\text{ or }e=ev, \\ 
\forall e\in E(G_{1})
\end{array}
\right. \right\} .$
\end{center}

\strut

Notice the difference between full-subgraphs and subgraphs: we say that $%
G_{1}$ is a \emph{subgraph} of $G_{2},$ if

\strut

\begin{center}
$V(G_{1})$ $\subseteq $ $V(G_{2})$
\end{center}

and

\begin{center}
$E(G_{1})$ $=$ $\left\{ e\in E(G_{2})\left| 
\begin{array}{c}
e=v\text{ }e\text{ }v^{\prime }, \\ 
\forall v,\text{ }v^{\prime }\in V(G_{2})
\end{array}
\right. \right\} .$
\end{center}

\strut

Every subgraph is a full-subgraph, but the converse does not hold true, in
general.\strut

\strut

\begin{definition}
Let $\mathcal{T}_{G}$ be the above one-flow growing rooted tree, where the
automata actions $\{\mathcal{A}_{w}$ $:$ $w$ $\in $\thinspace $\Bbb{F}^{+}(%
\widehat{G})\}$ of the graph automaton $\mathcal{A}_{G}$ act. This tree $%
\mathcal{T}_{G}$ is called the $\mathcal{A}_{G}$-tree, which is a
full-subgraph of the $2N$-regular tree $\mathcal{T}_{2N}.$ The
full-subgraphs $\mathcal{T}_{w}$'s for $w$ $\in $ $FP(\widehat{G})$ of $%
\mathcal{T}_{G}$ are called the $w$-parts of $\mathcal{T}_{G}.$
\end{definition}

\strut \strut

\strut The important thing is now that the $w$-parts $\mathcal{T}_{w}$'s of
the $\mathcal{A}_{G}$-graph $\mathcal{T}_{G}$ are embedded in the $2N$%
-regular tree $\mathcal{T}_{2N},$ and $\mathcal{T}_{w^{\prime }}$'s are
embedded in $\mathcal{T}_{w}$'s, whenever

\strut

\begin{center}
$\varphi \left( \omega (v),\text{ }w^{\prime }\right) $ $=$ $\left( (v,\text{
}v^{\prime }),\text{ }X^{\prime }\right) ,$
\end{center}

where

\begin{center}
$\varphi \left( \omega (v),\text{ }w\right) $ $=$ $\left( (v,\text{ }%
v^{\prime \prime }),\text{ }X\right) $ and $X^{\prime }$ $=$ $(X,$ $%
X^{\prime \prime }),$
\end{center}

\strut

for some $X^{\prime \prime }$ $\in $ $(\pm X_{0}^{*}),$ where $\omega (v)$
is the root of $\mathcal{T}_{G}$ $=$ $\mathcal{T}_{v}.$

\strut

\begin{remark}
The construction of the $\mathcal{A}_{G}$-tree $\mathcal{T}_{G}$ is nothing
but the rearrangement of the finite paths in the free semigroupoid $\Bbb{F}%
^{+}(\widehat{G})$ of the shadowed graph $\widehat{G}$ inside the $2N$%
-regular tree $\mathcal{T}_{2N},$ up to the admissibility on $\Bbb{F}^{+}(%
\widehat{G})$. Notice that, in fact, the $\mathcal{A}_{G}$-tree $\mathcal{T}%
_{G}$ contains the information about the vertices in $\Bbb{F}^{+}(\widehat{G}%
),$ too, since the vertices of $\mathcal{T}_{G}$ are contained in $\mathcal{X%
}_{0}^{*}$ $=$ $V(G)^{2}$ $\times $ $(\pm X_{0}^{*}).$ Remark that, by the
connectedness of $\widehat{G},$ the $w$-parts $\mathcal{T}_{w}$'s of $%
\mathcal{T}_{G}$ for $w$ $\in $ $FP(\widehat{G})$ are well-constructed as a
one-flow growing tree with their roots $\varphi (\omega (v),$ $w)$, where $%
\omega (v)$ is the root of $\mathcal{T}_{G}.$ Moreover, each tree $\mathcal{T%
}_{w}$ is embedded in the other trees $\mathcal{T}_{G}.$
\end{remark}

\strut \strut \strut \strut \strut

Remark that, even though $\left| \Bbb{G}\right| $ $<$ $\infty ,$ in general, 
$\left| \Bbb{F}^{+}(\widehat{G})\right| $ $=$ $\infty ,$ whenever $\left|
E(G)\right| $ $\geq $ $1.$ In [17], we observed the following special
labeled graph groupoids.

\strut

\begin{definition}
Let $G$, $\Bbb{F}^{+}(\widehat{G}),$ and $\Bbb{G}$ be given as before and
let $\mathcal{A}_{G}$ $=$ $<\mathcal{X}_{0},$ $E(\widehat{G}),$ $\varphi ,$ $%
\psi >$ be the graph automaton induced by $G,$ acting on the $2N$-regular
tree $\mathcal{T}_{2N}$. For any fixed $w$ $\in $ $FP(\widehat{G}),$ the
tree $\mathcal{T}_{w}$ is the $w$-part of the $\mathcal{A}_{G}$-tree $%
\mathcal{T}_{G},$ with its root $\varphi (\omega (v),$ $w),$ where $\omega
(v)$ is the fixed root of the $\mathcal{A}_{G}$-tree $\mathcal{T}_{G}.$ Let $%
\Bbb{G(\mathcal{T}}_{w}\Bbb{)}$ be the groupoid generated by the actions $%
\mathcal{A}_{y}$'s acting only on $\mathcal{T}_{w}$. If $\Bbb{G}(\mathcal{T}%
_{w})$'s are groupoid-isomorphic to $\Bbb{G}(\mathcal{T}_{G}),$ for all $w$ $%
\in $ $FP(\widehat{G}),$ then we say that the groupoid $\Bbb{G}(\mathcal{A}%
_{G}),$ generated by $\mathcal{A}_{G},$ is a fractaloid. Equivalently, we
say that the (labeled) graph groupoid $\Bbb{G}$ is a fractaloid.
\end{definition}

\strut

Notice that

\strut

\begin{center}
$\Bbb{G}(\mathcal{A}_{G})$ $\overset{\text{Groupoid}}{=}$ $\Bbb{G}$ $%
\overset{\text{Groupoid}}{=}$ $\Bbb{G}(\mathcal{T}_{G}).$
\end{center}

\strut

The first groupoid-isomorphic relation is shown in the previous section, and
the second groupoid-isomorphic relation holds, by the previous remark.
Readers can understand the above definition of fractaloids as the
graph-groupoid version of the fractal groups (See [1] and [17]). The
following theorem provides the graph-theoretical characterization of
fractaloids.

\strut

\begin{theorem}
(See [17]) Let $G$ be a canonical weighted graph with its labeled graph
groupoid $\Bbb{G},$ and let $\mathcal{A}_{G}$ be the graph automaton induced
by $G$ and $\mathcal{T}_{G},$ the $\mathcal{A}_{G}$-tree. Every $w$-part $%
\mathcal{T}_{w}$ of $\mathcal{T}_{G}$ is graph-isomorphic to $\mathcal{T}%
_{G},$ for all $w$ $\in $ $FP(\widehat{G}),$ if and only if $\Bbb{G}$ is a
fractaloid. $\square $
\end{theorem}

\strut \strut

Let $G$ be a canonical weighted graph with its labeled graph groupoid $\Bbb{G%
},$ and assume that the automata actions $\{\mathcal{A}_{w}$ $:$ $w$ $\in $ $%
FP(\widehat{G})\}$ of the graph automaton $\mathcal{A}_{G}$\emph{\ act fully
on the }$2N$\emph{-regular tree} $\mathcal{T}_{2N},$ in the sense that the $%
\mathcal{A}_{G}$-tree $\mathcal{T}_{G},$ which is a full-subgraph of $%
\mathcal{T}_{2N},$ is graph-isomorphic to $\mathcal{T}_{2N}.$ i.e., the
automata actions of $\mathcal{A}_{G}$ act fully on $\mathcal{T}_{2N},$ if $%
\mathcal{T}_{G}$ $\overset{\text{Graph}}{=}$ $\mathcal{T}_{2N}$.

\strut

\begin{corollary}
(See [17]) Let $G$ be a canonical weighted graph with its labeled graph
groupoid $\Bbb{G},$ and let $\mathcal{A}_{G}$ be the graph automaton induced
by $G.$ The automata actions $\{\mathcal{A}_{w}$ $:$ $w$ $\in $ $FP(\widehat{%
G})\}$ of $\mathcal{A}_{G}$ act fully on the $2N$-regular tree $\mathcal{T}%
_{2N},$ if and only if $\Bbb{G}$ is a fractaloid. $\square $
\end{corollary}

\strut \strut \strut

The above theorem provides the graph-theoretical characterization of
fractaloids: $\Bbb{G}$\strut \strut is a fractaloid if and only if

$\strut $

\begin{center}
$\mathcal{T}_{G}$ $\overset{\text{Graph}}{=}$ $\mathcal{T}_{2N}$ $\overset{%
\text{Graph}}{=}$ $\mathcal{T}_{w},$ for all $w$ $\in $ $\Bbb{F}^{+}(%
\widehat{G}).$
\end{center}

\strut

And the above corollary provides the algebraic (or automata-theoretical)
characterization of fractaloids; $\Bbb{G}$ is a fractaloid if and only if
the automata actions of the graph automaton $\mathcal{A}_{G}$ act fully on $%
\mathcal{T}_{2N}.$\strut 

\strut

\strut

\strut

\section{Labeling Operators of Graph Groupoids\strut}

\strut

\strut

In this section, we define labeling operators of labeled graph groupoids in
right graph von Neumann algebras. Let $G$ be a canonical weighted graph with
its labeled graph groupoid $\Bbb{G}.$ Then the labeling process on $\Bbb{G}$
induces the labeling operator $T_{G}$ of $\Bbb{G}$ contained in the right
graph von Neumann algebra $M_{G}$ $=$ $\overline{\Bbb{C}[\Bbb{G}]}^{w}$ of $%
G $, in $B(H_{G}).$

\strut

\strut

\subsection{The Labeling Operator of $\Bbb{G}$ in $M_{G}$}

\strut

\strut

Let $G$ $=$ $(G,$ $\omega )$ be a canonical weighted graph, where $G$ is a
countable locally finite connected directed graph, and let $\Bbb{G}$ $=$ $(%
\Bbb{G},$ $\omega )$ be the corresponding labeling graph groupoid of $G.$
Let $H_{G}$ be the graph Hilbert space in the sense of Section 2.3.
Throughout this section, put

\strut

\begin{center}
$N$ $=$ $\max \{\deg _{out}(v)$ $:$ $v$ $\in $ $V(G)\}.$
\end{center}

\strut \strut \strut

\begin{definition}
Define an operator $T_{k}$ $\in $ $B(H_{G})$ by

\strut 

\begin{center}
$
\begin{array}{ll}
T_{k}\xi _{w} & \overset{def}{=}\xi _{w}\xi _{e_{l_{k}}}=\xi _{we_{l_{k}}}
\\ 
&  \\ 
& =\left\{ 
\begin{array}{ll}
\xi _{we_{k}} & 
\begin{array}{l}
\text{if }\exists \text{ }e_{l_{k}}\in E(\widehat{G})\text{ s.t.} \\ 
w\text{ }e_{l_{k}}\neq \emptyset \text{ in }FP_{r}(\widehat{G}) \\ 
\text{and }(\widetilde{pr}\circ \omega )(e_{l_{k}})=l_{k}
\end{array}
\\ 
&  \\ 
\xi _{\emptyset }=0_{H_{G}} & \text{otherwise,}
\end{array}
\right. 
\end{array}
$
\end{center}

\strut 

for all $w$ $\in $ $\Bbb{G},$ where $\xi _{w}$ $\in $ $\mathcal{B}_{H_{G}}$ $%
\cup $ $\{\xi _{v}$ $:$ $v$ $\in $ $V(\widehat{G})\}$, and where $(%
\widetilde{pr}$ $\circ $ $\omega )(e_{l_{k}})$ $=$ $l_{k}$ $\in $ $\pm X_{0}$
is the weight of the edge $e_{l_{k}},$ where $\pm X_{0}$ $=$ $\{l_{0},$ $\pm
l_{1},$ ..., $\pm l_{N}\},$ for $k$ $=$ $\pm 1,$ ..., $\pm N.$ These
operators $T_{k}$'s are called the $k$-th labeling operators, for all $k$ $=$
$\pm 1,$ ..., $\pm N.$ Then we can define the operator $T_{G}$ $\in $ $%
B(H_{G})$ by

\strut 

\begin{center}
$T_{G}$ $\overset{def}{=}$ $\underset{k=-N}{\overset{-1}{\sum }}$ $T_{k}$ $+$
$\underset{j=1}{\overset{N}{\sum }}$ $T_{j}.$
\end{center}

\strut 

This operator $T$ is called the labeling operator of $\Bbb{G}$ on $H_{G}.$
For convenience, denote $\underset{k=-N}{\overset{-1}{\sum }}$ $T_{k}$ and $%
\underset{j=1}{\overset{N}{\sum }}$ $T_{j}$ by $T_{-}$ and $T_{+},$
respectively.
\end{definition}

\strut

The labeling operator $T_{G}$ of $\Bbb{G}$ on $H_{G}$ is a Heck-type
operator or a Ruelle-operator-like operator (See [1], [17], [33] and [34]).
The labeling operators are originally defined in [17]. The following
proposition was already proven in [17]. It shows that we can regard the $k$%
-th labeling operators $T_{k}$'s, for $k$ $=$ $\pm 1,$ ..., $\pm N,$ and the
labeling operator $T_{G},$ as elements of the right graph von Neumann
algebra $M_{G}$ of $G.$

\strut \strut

\begin{proposition}
(See [17]) Let $G$ be the given graph with its labeled graph groupoid $\Bbb{G%
},$ and let $M_{G}$ $=$ $\overline{\Bbb{C}[\Bbb{G}]}^{w}$ be the right graph
von Neumann algebra of $G,$ in $B(H_{G}).$ Then the $k$-th labeling
operators $T_{k}$ and the labeling operator $T_{G}$ of $\Bbb{G}$ are
contained in $M_{G}$, for $k$ $=$ $\pm 1,$ ..., $\pm N.$ $\square $
\end{proposition}

\strut

Indeed, we can define the elements $\tau _{k}$ $\in $ $M_{G},$ by

\strut

\begin{center}
$\tau _{k}$ $\overset{def}{=}$ $\underset{e\in E(\widehat{G}),\,\,\left( 
\widetilde{pr}\text{ }\circ \text{ }\omega \right) (e)=l_{k}}{\sum }$ $R_{e}$
$\in $ $M_{G},$ for $k$ $=$ $\pm 1,$ ..., $\pm N.$
\end{center}

\strut

Then these operator $\tau _{k}$ satisfies that

\strut \strut

\begin{center}
$
\begin{array}{ll}
\tau _{k}\xi _{w} & =\underset{e\in E(\widehat{G}),\left( \,\widetilde{pr}%
\text{ }\circ \text{ }\omega \right) (e)=l_{k}}{\sum }\left( R_{e}\xi
_{w}\right) =\underset{e\in E(\widehat{G}),\left( \,\widetilde{pr}\text{ }%
\circ \text{ }\omega \right) (e)=l_{k}}{\sum }\xi _{we} \\ 
&  \\ 
& =\left\{ 
\begin{array}{ll}
\xi _{we_{0}} & \text{if }we_{0}\neq \emptyset \\ 
&  \\ 
0_{H_{G}} & \text{otherwise.}
\end{array}
\right.
\end{array}
$
\end{center}

\strut

By the definition of labeling processes, (i) if the above $H_{G}$-value is
nonzero, then there exists a ``unique'' $e_{0}$ $\in $ $E(\widehat{G}),$
with $\omega (e_{0})$ $=$ $((v,$ $v^{\prime }),$ $l_{k}),$ whenever $w$ $=$ $%
wv$ and $e_{0}$ $=$ $v$ $e_{0}$ $v^{\prime },$ and hence $(\widetilde{pr}$ $%
\circ $ $\omega )$ $=$ $l_{k},$ and (ii) if the above $H_{G}$-value is zero,
then we can conclude $w$ $=$ $wv,$ with $\deg _{out}(v)$ $=$ $n_{v},$ and $%
\left| n_{v}\right| $ $\lneqq $ $\left| k\right| ,$ for $k$ $\in $ $\{\pm 1,$
..., $\pm N\}.$ i.e.,

\strut \strut

\begin{center}
$\left\{ 
\begin{array}{ll}
\xi _{we_{0}}=T_{k}\text{ }\xi _{w} & 
\begin{array}{l}
\text{if }\exists !\text{ }e_{0}\in E(\widehat{G})\text{ s.t. } \\ 
\quad (\widetilde{pr}\text{ }\circ \text{ }\omega )(e_{0})=x_{k}
\end{array}
\\ 
&  \\ 
\xi _{\emptyset }=0_{H_{G}} & \text{otherwise.}
\end{array}
\right. $
\end{center}

\strut

Therefore, the $k$-th labeling operators $T_{k}$ are identical to $\tau _{k}$
$\in $ $M_{G},$ and hence, by definition, the labeling operator $T_{G}$ of $%
\Bbb{G}$ is identified with the element $\tau $ $\in $ $M_{G},$ defined by

\strut

\begin{center}
$\tau $ $=$ $\left( \underset{k=-N}{\overset{-1}{\sum }}\tau _{k}\right) $ $%
+ $ $\left( \underset{j=1}{\overset{N}{\sum }}\tau _{j}\right) .$
\end{center}

\strut \strut \strut

The above proposition shows that, without loss of generality, we can regard
the $k$-th labeling operators $T_{k}$'s and the labeling operator $T_{G}$ of 
$\Bbb{G}$ as elements of the right graph von Neumann algebra $M_{G},$ for $k$
$=$ $\pm 1,$ ..., $\pm N.$ Recall that, the right graph von Neumann algebra $%
\Bbb{M}_{G}$ $=$ $\overline{\Bbb{C}[\Bbb{G}]}^{w}$ $\subseteq $ $B(H_{G})$
is $*$-isomorphic to a $D_{G}$-valued reduced free product algebra, where $%
D_{G}$ is the $\Bbb{C}$-diagonal subalgebra of $M_{G}.$ Therefore, we can
understand the operators $T_{k}$'s and $T_{G},$ as $D_{G}$-valued random
variables in $(M_{G},$ $E)$, where $E$ $:$ $M_{G}$ $\rightarrow $ $D_{G}$ is
the canonical conditional expectation. So, the labeling operators have their 
$D_{G}$-valued free distributional data.

\strut

\strut

\strut

\subsection{$D_{G}$-Valued Free Distributional Data of Labeling Operators}

\strut

\strut

As before, let $G$ $=$ $(G,$ $\omega )$ be a canonical weighted graph with
its labeled graph groupoid $\Bbb{G}$ $=$ $(\Bbb{G},$ $\omega ).$ In this
section, we will consider the $D_{G}$-valued free distributional data of the 
$k$-th labeling operators $T_{k}$'s and the labeling operator $T_{G}$ of $%
\Bbb{G}$ on the graph Hilbert space $H_{G},$ by regarding them as $D_{G}$%
-valued random variables in the right graph $W^{*}$-probability space $%
(M_{G},$ $E)$ over the $\Bbb{C}$-diagonal subalgebra $D_{G},$ where

\strut

\begin{center}
$N$ $\overset{def}{=}$ $\max \{\deg _{out}(v)$ $:$ $v$ $\in $ $V(G)\}$ $\in $
$\Bbb{N}.$
\end{center}

\strut

\textbf{Notation} Without loss of generality, we will use the notations $%
T_{k}$'s and $\tau _{k}$'s (resp., $T_{G}$ and $\tau $), alternatively, for $%
k$ $=$ $\pm 1,$ ..., $\pm N.$ $\square $

\strut \strut

Now, fix $k$ $\in $ $\{\pm 1,$ ..., $\pm N\},$ and consider the $D_{G}$%
-valued joint moments $E(T_{k}^{r_{1}}$ ... $T_{k}^{r_{n}})$ of the $k$-th
labeling operator

\strut

\begin{center}
$T_{k}$ $=$ $\underset{e\in E(\widehat{G}),\text{ }(\widetilde{pr}\text{ }%
\circ \text{ }\omega )(e)=l_{k}}{\sum }$ $R_{e}$ $\in $ $M_{G},$
\end{center}

\strut

where $(r_{1},$ ..., $r_{n})$ $\in $ $\{1,$ $*\}^{n},$ for all $n$ $\in $ $%
\Bbb{N}.$

\strut

\begin{lemma}
(See [17]) Let $T_{k}$ be the labeling operator on the graph Hilbert space $%
H_{G}.$ Then its adjoint $T_{k}^{*}$ is identified with the $(-k)$-th
labeling operator $T_{-k},$ for $k$ $\in $ $\{\pm 1,$ ..., $\pm N\}.$ $%
\square $
\end{lemma}

\strut

Indeed, for each $k,$ we can have

\strut (5.1)

\begin{center}
$
\begin{array}{ll}
T_{k}^{*} & =\underset{e\in E(\widehat{G}),\text{ }(\widetilde{pr}\text{ }%
\circ \text{ }\omega )(e)=l_{k}}{\sum }R_{e}^{*}=\underset{e\in E(\widehat{G}%
),\text{ }(\widetilde{pr}\text{ }\circ \text{ }\omega )(e)=l_{k}}{\sum }%
R_{e^{-1}} \\ 
&  \\ 
& =\underset{e^{-1}\in E(\widehat{G}),\text{ }(\widetilde{pr}\text{ }\circ 
\text{ }\omega )(e^{-1})=l_{-k}}{\sum }R_{e^{-1}}=T_{-k},
\end{array}
$
\end{center}

\strut

by the canonical weighting process, for $k$ $=$ $\pm 1,$ ..., $\pm N.$\strut
\ Thus, for any $k,$ the adjoint $T_{k}^{*}$ of the $k$-th labeling operator 
$T_{k}$ is the $(-k)$-th labeling operator $T_{-k},$ for all $k$ $=$ $\pm 1,$
..., $\pm N.$ By the above lemma, we can conclude the self-adjointness of
the labeling operator $T_{G}$ on $H_{G}.$

\strut

\begin{corollary}
(See [17]) Let $T_{G}$ be the labeling operator of $\Bbb{G}$ on $H_{G}.$
Then it is self-adjoint. $\square $
\end{corollary}

\strut \strut \strut

The above corollary guarantees that the study of $D_{G}$-valued free
distributional data of\thinspace $T_{G}$ is to study the spectral property
of $T_{G}$. Now, we will consider the $D_{G}$-freeness of $\{T_{k},$ $%
T_{-k}\}$'s in the graph $W^{*}$-probability space $(M_{G},$ $E),$ for all $%
k $ $=$ $1,$ ..., $N.$

\strut

\begin{theorem}
(See [17]) Let $\{T_{k}$ $:$ $k$ $=$ $\pm 1,$ ..., $\pm N\}$ be the $k$-th
labeling operators. Then the families $\{T_{k},$ $T_{-k}\}$'s are free over $%
D_{G}$ from each other in $(M_{G},$ $E),$ for $k$ $=$ $1,$ ..., $N.$ i.e.,
the families $\{T_{k_{1}},$ $T_{-k_{1}}\}$ and $\{T_{k_{2}},$ $T_{-k_{2}}\}$
are free over $\Bbb{D}_{G}$ in $(\Bbb{M}_{G},$ $E),$ whenever $k_{1}$ $\neq $
$k_{2}$ in $\{1,$ ..., $N\}.$ $\square $
\end{theorem}

\strut \strut \strut \strut \strut

The above theorem can be proven, by the fact: the subsets

$\strut $

\begin{center}
$\mathcal{R}_{k}$ $\overset{def}{=}$ $\{e$ $\in $ $E(\widehat{G})$ $:$ $(%
\widetilde{pr}$ $\circ $ $\omega )(e)$ $=$ $k$ or $-k\}$ $\subset $ $\Bbb{G}%
, $
\end{center}

\strut

for $k$ $=$ $1,$ ..., $N,$ are diagram-distinct from each other.

\strut Consider now the $D_{G}$-valued joint $*$-moments of $T_{\pm 1},$
..., $T_{\pm N}$. Let $(i_{1},$ ..., $i_{n})$ $\in $ $\{\pm 1,$ ..., $\pm
N\}^{n},$ for $n$ $\in $ $\Bbb{N}.$ Observe that

\strut

$\qquad E\left( T_{i_{1}}\text{ ... }T_{i_{n}}\right) $

$\strut \strut $

$\quad \qquad =$ $E\left( \left( \underset{e_{1}\in E(\widehat{G}),\text{ }(%
\widetilde{pr}\text{ }\circ \text{ }\omega )(e_{1})=l_{i_{1}}}{\sum }\text{ }%
R_{e_{1}}\right) \text{ ... }\left( \underset{e_{n}\in E(\widehat{G}),\text{ 
}(\widetilde{pr}\text{ }\circ \text{ }\omega )(e_{n})=l_{i_{n}}}{\sum }\text{
}R_{e_{n}}\right) \right) $

\strut

$\quad \qquad =$ $E\left( \underset{(e_{1},...,e_{n})\in E(\widehat{G})^{n},%
\text{ }(\widetilde{pr}\text{ }\circ \text{ }\omega )(e_{j})=l_{i_{j}}}{\sum 
}\text{ }R_{e_{n}e_{n-1}...e_{2}e_{1}}\right) $

\strut (5.2)

$\qquad \quad =$ $E\left( \underset{w\in E(\widehat{G})^{n}\text{ }\cap 
\text{ }FP(\widehat{G}),\text{ }(\widetilde{pr}\text{ }\circ \text{ }\omega
)(w)=l_{i_{1}}...l_{i_{n}}}{\sum }\text{ }R_{w}\right) ,$

\strut

where $E(\widehat{G})^{n}$ $=$ $\underset{n\text{-times}}{\underbrace{E(%
\widehat{G})\times ...\times E(\widehat{G})}},$ for $n$ $\in $ $\Bbb{N}.$

\strut

Now, recall the definition of the operations $\theta $ and $\omega _{+},$
introduced in Section 4.2.

\strut

\begin{lemma}
Let $w$ $\in $ $E(\widehat{G})^{n}$ $\cap $ $FP(\widehat{G}),$ with $\omega
(w)$ $=$ $\left( (v,\text{ }v),\text{ }l_{i_{1}}...l_{i_{n}}\right) ,$ for $n
$ $\in $ $\Bbb{N},$ for some $v$ $\in $ $V(\widehat{G}).$ Then $\omega
_{+}\left( \omega (w)\right) $ $=$ $\left( (v,\text{ }v),\text{ }0\right) $
if and only if $R_{w}$ $=$ $R_{v},$ for $v$ $\in $ $V(G).$
\end{lemma}

\strut

\begin{proof}
($\Rightarrow $) Assume that $w$ $\in $ $\Bbb{G},$ and let

\strut

$\qquad \omega _{+}\left( \omega (w)\right) $ $=$ $\omega _{+}\left( \left(
(v,\text{ }v),\text{ }l_{i_{1}}\text{ ... }l_{i_{n}}\right) \right) $ $=$ $%
\left( (v,\text{ }v),\text{ }\theta (l_{i_{1}}\text{ ... }l_{i_{n}})\right) $

\strut

$\qquad \qquad \qquad \quad =$ $\left( (v,\text{ }v),\text{ }\sum_{j=1}^{n}%
\text{ }i_{j}\right) $ $=$ $\left( (v,\text{ }v),\text{ }0\right) .$

\strut

Then, we can conclude that $w$ $=$ $v$ $w$ $v,$ for $v$ $\in $ $V(G),$ and $%
w $ $=$ $e_{i_{1}}$ ... $e_{i_{n}},$ with $e_{i_{j}}$ $\in $ $E(\widehat{G}),
$ having $\omega (e_{i_{j}})$ $=$ $l_{i_{j}},$ for all $j$ $=$ $1,$ ..., $n.$
Since $\sum_{j=1}^{n}$ $i_{j}$ $=$ $0$ in $\Bbb{Z},$ the number $n$ should
be an even number in $2\Bbb{N},$ and $\{i_{1},$ ..., $i_{n}\}$ is decomposed
by $\{j_{1},$ ..., $j_{\frac{n}{2}}\}$ and $\{q_{1},$ ..., $q_{\frac{n}{2}%
}\}.$ Futhermore, the set $\{e_{i_{1}},$ ..., $e_{i_{n}}\}$ is decomposed by 
$\mathcal{E}_{1}$ $=$ $\{e_{j_{1}},$ ..., $e_{j_{\frac{n}{2}}}\}$ and $%
\mathcal{E}_{2}$ $=$ $\{e_{q_{1}},$ ..., $e_{q_{\frac{n}{2}}}\}.$ By the
definitions of $\omega _{+}$, we can conclude that $\mathcal{E}_{2}$ $=$ $%
\mathcal{E}_{1}^{-1}$ in $E(\widehat{G}).$

\strut

Therefore, by [10] and [11], $w$ $=$ $v$ $w$ $v$ $=$ $v$ $v$ $v$ $=$ $v$ in $%
V(G).$

\strut

($\Leftarrow $) Suppose $R_{w}$ $=$ $R_{v},$ for $v$ $\in $ $V(G).$ Then

\strut

\begin{center}
$\omega _{+}\left( \omega (v)\right) $ $=$ $\omega _{+}\left( \left( (v,%
\text{ }v),\text{ }l_{0}\right) \right) $ $=$ $\left( (v,\text{ }v),\text{ }%
0\right) .$
\end{center}

\strut
\end{proof}

\strut \strut \strut \strut

By the previous observations, we obtain the following theorem providing the $%
D_{G}$-valued free distributional data of the $k$-th labeling operators $%
T_{\pm 1},$ ..., $T_{\pm N}.$

\strut

\begin{theorem}
Let $T_{k}$'s be the $k$-th labeling operators on $H_{G},$ for all $k$ $=$ $%
\pm 1,$ ..., $\pm N.$ Then the $\Bbb{D}_{G}$-valued joint $*$-moments of
them is determined by

\strut 

(5.3)

\begin{center}
$E\left( T_{i_{1}}\text{ ... }T_{i_{n}}\right) $ $=$ $\underset{w\in E(%
\widehat{G})^{n}\cap FP(\widehat{G}),\text{ }\omega _{+}\left( \omega
(w)\right) =\left( (v,\text{ }v),\text{ }0\right) }{\sum }$ $R_{w},$
\end{center}

\strut 

in $\Bbb{D}_{G},$ for all $(i_{1},$ ..., $i_{n})$ $\in $ $\{\pm 1,$ ..., $%
\pm N\},$ for $n$ $\in $ $\Bbb{N}.$ $\square $
\end{theorem}

\strut

The formula (5.3) provides the $D_{G}$-valued free distributional data of $%
\{T_{k}$ $:$ $k$ $=$ $\pm 1,$ ..., $\pm N\}$ in $(M_{G},$ $E),$ in terms of
their $D_{G}$-valued joint $*$-moments. By the Moebius inversion (e.g., [5]
and [21]: also see Section 1.1), we can get that the $D_{G}$-valued joint $*$%
-cumulants of $\{T_{k}$ $:$ $k$ $=$ $\pm 1,$ ..., $\pm N\}$:

\strut

\begin{center}
$k_{n}\left( T_{i_{1}},\text{ ..., }T_{i_{n}}\right) $ $=$ $\underset{\pi
\in NC(n)}{\sum }$ $E_{\pi }\left( T_{i_{1}},\text{ ..., }T_{i_{n}}\right) $ 
$\mu (\pi ,$ $1_{n}).$
\end{center}

\strut

Similar to [10] and [11], we obtain the following lemma.

\strut

\begin{lemma}
Let $w_{1},$ ..., $w_{n}$ $\in $ $\Bbb{G}$ and let $R_{w_{1}},$ ..., $%
R_{w_{n}}$ $\in $ $M_{G}$ be the corresponding right multiplication
operators, for $n$ $\in $ $\Bbb{N}.$ Then

\strut 

\strut (5.4)

\begin{center}
$k_{n}\left( R_{w_{1}},\text{ }R_{w_{2}},\text{ ..., }R_{w_{n}}\right) $ $=$ 
$\mu _{0}$ $E\left( R_{w_{n}...w_{2}w_{1}}\right) ,$
\end{center}

where

\begin{center}
$\mu _{0}$ $=$ $\underset{\pi \in NC(w_{1},...,w_{n})}{\sum }$ $\mu (\pi ,$ $%
1_{n}),$
\end{center}

where

\begin{center}
$NC(w_{1},$ ..., $w_{n})$ $\overset{def}{=}$ $\left\{ \pi \in NC(n)\left| 
\begin{array}{c}
E_{\pi }(R_{w_{1}},\text{ ..., }R_{w_{n}}) \\ 
=E\left( R_{w_{n}...w_{2}w_{1}}\right) \neq 0_{D_{G}}
\end{array}
\right. \right\} .$
\end{center}

$\square $
\end{lemma}

\strut

So, the above lemma will let us find the $D_{G}$-valued free distributional
data of the $k$-th labeling operators $\{T_{\pm 1},$ ..., $T_{\pm N}\},$ in
terms of the $D_{G}$-valued joint $*$-cumulants of them. Notice that, again
by Moebius inversion, the $D_{G}$-valued free distributional data of $T_{\pm
1},$ ..., $T_{\pm N},$ represented by $D_{G}$-valued $*$-moments and by $%
D_{G}$-valued $*$-cumulants, are equivalent.

\strut

Let $(i_{1},$ ..., $i_{n})$ $\in $ $\{\pm 1,$ ..., $\pm N\}^{n},$ for $n$ $%
\in $ $\Bbb{N}.$ Observe that

\strut

$\qquad k_{n}\left( T_{i_{1}},\text{ ..., }T_{i_{n}}\right) $

\strut

$\quad \qquad =$ $k_{n}\left( \left( \underset{e_{1}\in E(\widehat{G}),\text{
}(\widetilde{pr}\text{ }\circ \text{ }\omega )(e)=l_{i_{1}}}{\sum }\text{ }%
R_{e_{1}}\right) ,\text{ ..., }\left( \underset{e_{n}\in E(\widehat{G}),%
\text{ }(\widetilde{pr}\text{ }\circ \text{ }\omega )(e)=l_{i_{n}}}{\sum }%
\text{ }R_{e_{n}}\right) \right) $

\strut

$\quad \qquad =$ $\underset{(e_{1},...,e_{n})\in E(\widehat{G}),\text{ }(%
\widetilde{pr}\text{ }\circ \text{ }\omega )(e_{i_{j}})=l_{i_{j}}}{\sum }$ $%
k_{n}\left( R_{e_{1}},\text{ ..., }R_{e_{n}}\right) $

\strut

by the bimodule map property of $k_{n}(...)$ (See [21])

\strut

\strut (5.5)

$\qquad \quad =$ $\underset{(e_{1},...,e_{n})\in E(\widehat{G})^{n}\cap FP(%
\widehat{G}),\text{ }(\widetilde{pr}\text{ }\circ \text{ }\omega
)(e_{i_{j}})=l_{i_{j}}}{\sum }$ $\left( \mu _{(e_{1},...,e_{n})}\text{ }%
E(R_{e_{n}...e_{2}e_{1}})\right) $

\strut

by (5.4), where

\strut

\begin{center}
$\mu _{(e_{1},...,e_{n})}$ $=$ $\underset{\pi \in NC(e_{1},...,e_{n})}{\sum }
$ $\mu (\pi ,$ $1_{n}),$
\end{center}

where

\begin{center}
$NC(e_{1},$ ..., $e_{n})$ $=$ $\left\{ \pi \in NC(n)\left| 
\begin{array}{c}
E_{\pi }(R_{e_{1}},\text{ ..., }R_{e_{n}}) \\ 
\;=E\left( R_{e_{n}e_{n-1}...e_{2}e_{1}}\right) \\ 
\neq 0_{D_{G}}
\end{array}
\right. \right\} .$
\end{center}

\strut

Thus, the formula (5.5) is identified with the following formula (5.6)

\strut

(5.6)

\begin{center}
\strut $\underset{w=e_{1}...e_{n}\in E(\widehat{G})^{n}\cap FP(\widehat{G}),%
\text{ }\omega _{+}\left( \omega (w)\right) =\left( (v,\,v),\text{ }0\right) 
}{\sum }$ $\mu _{w}$ $R_{w},$
\end{center}

\strut

where

\strut

\strut (5.7)

\begin{center}
$\mu _{w}$ $=$ $\mu _{(e_{1},...,e_{n})},$ whenever $w$ $=$ $e_{n}$ ... $%
e_{1}$ $\in $ $\Bbb{G}.$
\end{center}

\strut

Therefore, we can get the following equivalent $D_{G}$-valued free
distributional data of the $k$-th labeling operators $T_{k}$'s, for $k$ $=$ $%
\pm 1,$ ..., $\pm N$:

\strut

\begin{theorem}
Let $T_{k}$'s be the $k$-th labeling operators, for $k$ $=$ $\pm 1,$ ..., $%
\pm N.$ Then the $\Bbb{D}_{G}$-valued joint $*$-cumulants of them are

\strut 

(5.8)$\qquad \quad $

\begin{center}
$k_{n}\left( T_{i_{1}},\text{ ..., }T_{i_{n}}\right) $ $=$ $\underset{w\in E(%
\widehat{G})^{n}\cap FP(\widehat{G}),\text{ }\omega _{+}\left( \omega
(w)\right) =\left( (v,\text{ }v),\text{ }0\right) }{\sum }$ $\left( \mu
_{w}R_{w}\right) ,$
\end{center}

\strut 

for all $(i_{1},$ ..., $i_{n})$ $\in $ $\{\pm 1,$ ..., $\pm N\}^{n},$ for $n$
$\in $ $\Bbb{N},$ where $\mu _{w}$'s are given in (5.7). $\square $
\end{theorem}

\strut

In [10] and [11], we showed that $\mu _{w}$ $=$ $0,$ whenever there exists
at least one pair $(i,$ $j)$ $\in $ $\{1,$ ..., $n\}$ such that $e_{i}$ and $%
e_{j}$ are diagram-distinct, where $w$ $=$ $e_{1}$ ... $e_{n}$ $\in $ $FP(%
\widehat{G}),$ with $e_{1},$ ..., $e_{n}$ $\in $ $E(\widehat{G}).$
Therefore, we can have the following corollary.

\strut

\begin{corollary}
Let $T_{k}$'s be the $k$-th labeling operators, for $k$ $=$ $\pm 1,$ ..., $%
\pm N.$ Then

\strut 

(1) $k_{n}(T_{i_{1}},$ ..., $T_{i_{n}})$ $=$ $0_{D_{G}},$ whenever $n$ is
odd.

\strut 

(2) $k_{n}\left( T_{i_{1}},\text{ ..., }T_{i_{n}}\right) $ $=$ $\underset{%
w\in \mathcal{W}_{n}^{(c)}(i_{1},\text{ ..., }i_{n})}{\sum }$ $(\mu
_{w}R_{w}),$ where

\strut 

\strut (5.9)

\begin{center}
$\mathcal{W}_{n}^{(c)}(i_{1},$ ..., $i_{n})$ $=$ $\left\{
w=e_{1}...e_{n}\left| 
\begin{array}{c}
w\in E(\widehat{G})^{n}\cap FP(\widehat{G}) \\ 
(\widetilde{pr}\text{ }\circ \text{ }\omega )(w)=l_{i_{1}}...l_{i_{n}} \\ 
\theta \left( l_{i_{1}}...l_{i_{n}}\right) =0,\text{ and} \\ 
\exists \text{ }e\in E(\widehat{G})\text{ s.t. }e_{1}...e_{n}\text{ is a} \\ 
\text{word only in }\{e,\text{ }e^{-1}\}.
\end{array}
\right. \right\} ,$
\end{center}

\strut 

for all $n$ $\in $ $\Bbb{N}.$ $\square $
\end{corollary}

\strut \strut \strut \strut

Now, let's consider the $D_{G}$-valued free distributional data of the
labeling operator $T_{G}$ of $\Bbb{G}$ on the graph Hilbert space $H_{G}.$
Recall that the labeling operator $T_{G}$ is self-adjoint in $M_{G}$ $%
\subseteq $ $B(H_{G}).$ Therefore, it is enough to consider its $D_{G}$%
-valued (non-joint and non-$*$-)moments $\{E(T_{G}^{n})\}_{n\in \Bbb{N}}.$
Since $T_{G}$ is self-adjoint, the data $\{E(T_{G}^{n})$ $:$ $n$ $\in $ $%
\Bbb{N}\}$ contain the full $D_{G}$-valued free distributional data of it.

\strut

Fix $n$ $\in $ $\Bbb{N}$ and consider $T_{G}^{n}$:

\strut

\begin{center}
$
\begin{array}{ll}
T_{G}^{n} & =\left( T_{-}+T_{+}\right)
^{n}=(T_{-N}+...+T_{-1}+T_{1}+...+T_{N})^{n} \\ 
&  \\ 
& =\underset{(i_{1},\text{ ..., }i_{n})\in \{\pm 1,\text{ ..., }\pm N\}^{n}}{%
\sum }T_{i_{1}}...T_{i_{n}}.
\end{array}
$
\end{center}

\strut

Thus, we can compute that

\strut

\begin{center}
$E\left( T_{G}^{n}\right) $ $=$ $\underset{(i_{1},\text{ ..., }i_{n})\in
\{\pm 1,\text{ ..., }\pm N\}^{n}}{\sum }$ $E\left( T_{i_{1}}\text{ ... }%
T_{i_{n}}\right) .$
\end{center}

\strut

Therefore, by (5.4), we can obtain the following theorem.

\strut

\begin{theorem}
Let $T_{G}$ be the labeling operator of $\Bbb{G}$ on the graph Hilbert space 
$H_{G}.$ Then

\strut 

(5.10) $\qquad \qquad \quad $

\begin{center}
$E(T_{G}^{n})$ $=$ $\underset{w\in \mathcal{W}_{n}\text{ }}{\sum }$ $R_{w},$
\end{center}

\strut 

for all $n$ $\in $ $\Bbb{N},$ where $\mathcal{W}_{n}^{(m)}$ is defined by

\strut 

\begin{center}
$\mathcal{W}_{n}^{(m)}$ $\overset{def}{=}$ $\left\{ w\in E(\widehat{G}%
)^{n}\left| 
\begin{array}{c}
w\in FP(\widehat{G}), \\ 
\omega _{+}\left( \omega (w)\right) =\left( (v,\text{ }v),\text{ }0\right) ,
\\ 
\text{for some }v\in V(G)
\end{array}
\right. \right\} ,$
\end{center}

\strut 

for all $n$ $\in $ $\Bbb{N}.$ $\square $
\end{theorem}

\strut

The above theorem provides the $D_{G}$-valued free distributional data of
the labeling operator $T_{G},$ in terms of $D_{G}$-valued moments.
Similarly, we can have the following equivalent $D_{G}$-valued free
distributional data of $T_{G},$ in terms of $D_{G}$-valued cumulants.

\strut

\begin{theorem}
Let $T_{G}$ be the labeling operator of $\Bbb{G}$ on $H_{G}.$ Then the $D_{G}
$-valued cumulants of $T_{G}$ satisfy that: if $n$ is odd, then $k_{n}(T_{G},
$ ..., $T_{G})$ $=$ $0_{D_{G}},$ and if $n$ is even, then

\strut 

(5.11) $\qquad $

\begin{center}
$k_{n}\left( \underset{n\text{-times}}{\underbrace{T_{G},\text{ ..., }T_{G}}}%
\right) $ $=$ $\underset{(i_{1},\text{ ..., }i_{n})\in \{\pm 1,\text{ ..., }%
\pm N\}^{n}}{\sum }\left( \underset{w\in \mathcal{W}_{n}^{(c)}}{\sum }\mu
_{w}R_{w}\right) ,$
\end{center}

\strut 

where $\mu _{w}$'s are defined in (5.7) and where $\mathcal{W}_{n}^{(c)}$ is
defined by

\strut 

(5.12)

\begin{center}
$\mathcal{W}_{n}^{(c)}$ $\overset{def}{=}$ $\underset{(i_{1},...,i_{n})\in
\{\pm 1,\text{ ..., }\pm N\}^{n}}{\cup }$ $\mathcal{W}%
_{n}^{(c)}(i_{1},...,i_{n}),$ for all $n$ $\in $ $\Bbb{N}$.
\end{center}

\strut 

Here, $\mathcal{W}_{n}^{(c)}(i_{1},$ ..., $i_{n})$ is introduced in (5.9). $%
\square $
\end{theorem}

\strut

Recall that if an operator $a$ of a $B$-valued $W^{*}$-probability space $%
(A, $ $E_{B}),$ where $E_{B}$ $:$ $A$ $\rightarrow $ $B$ is a conditional
expectation, is self-adjoint, then the $B$-valued free distribution $\sigma
_{a}$ of $a$ and the (operator-valued) spectral measure $\chi _{a}$ of $a$
are identical. Therefore, the formulae (5.10) and (5.11) also contain the
spectral measure theoretical data of the labeling operator $T_{G}$ (See [5]
and [23]).

\strut

We will finish this section with the following corollary, which is the
modification of (5.10). By the fact that all odd $D_{G}$-valued cumulants of 
$T_{G}$ vanish, we can conclude that all odd $D_{G}$-valued moments of $%
T_{G} $ vanish, too, by the Moebius Inversion. Indeed, we have

\strut

\begin{center}
$E\left( T_{G}^{n}\right) $ $=$ $\underset{\pi \in NC(n)}{\sum }$ $k_{n:\pi
}(T_{G},$ ..., $T_{G}),$
\end{center}

\strut

where $k_{n:\pi }(\bullet ,$ ..., $\bullet )$ are the partition-depending $%
D_{G}$-valued cumulant. Notice that, if $n$ is odd, then every noncrossing
partition $\pi $ of $NC(n)$ contains an odd block, having its odd size.
Therefore, $k_{n:\pi }(T_{G},$ ..., $T_{G})$ $=$ $0_{D_{G}},$ for all $\pi $ 
$\in $ $NC(n),$ and hence $E(T_{G}^{n})$ $=$ $0_{D_{G}},$ whenever $n$ is
odd. So, the formula (5.10) can be modified by:

\strut

\begin{corollary}
Let $T_{G}$ be the labeling operator of the given graph groupoid $\Bbb{G}$
on the graph Hilbert space $H_{G}.$ Then

\strut 

\strut (5.10)$^{\prime }$

\begin{center}
$E(T_{G}^{n})$ $=$ $\left\{ 
\begin{array}{lll}
\underset{w\in \mathcal{W}_{n}^{(m)}}{\sum }\text{ }R_{w} &  & \text{if }n%
\text{ is even} \\ 
&  &  \\ 
0_{D_{G}} &  & \text{if }n\text{ is odd,}
\end{array}
\right. $
\end{center}

\strut 

for all $n$ $\in $ $\Bbb{N}.$ In particular, $\mathcal{W}_{n}^{(m)}$ is
empty, whenever $n$ is odd. $\square $
\end{corollary}

\strut \strut

\strut

\strut

\subsection{Graph-Theoretical Characterizations of $\mathcal{W}_{2n}^{(m)}$
and $\mathcal{W}_{2n}^{(c)}$}

\strut

\strut

As we have seen in (5.10) and (5.11), the $D_{G}$-valued free distributional
data of the labeling operator $T_{G}$ on the graph Hilbert space $H_{G}$ is
characterized by the sets $\{\mathcal{W}_{2n}^{(m)}$ $:$ $n$ $\in $ $\Bbb{N}%
\}$ and $\{\mathcal{W}_{2n}^{(c)}$ $:$ $n$ $\in $ $\Bbb{N}\}$, where

\strut

\begin{center}
$\mathcal{W}_{2n}^{(m)}$ $=$ $\left\{ w=e_{1}...e_{2n}\left| 
\begin{array}{c}
w\in E(\widehat{G})^{2n}\times FP(\widehat{G})\text{,} \\ 
\omega _{+}\left( \omega (w)\right) =\left( (v,\text{ }v),\text{ }0\right) ,
\\ 
\text{for some }v\in V(G)
\end{array}
\right. \right\} $
\end{center}

and

\begin{center}
$\mathcal{W}_{2n}^{(c)}$ $=$ $\left\{ w=e_{1}...e_{2n}\left| 
\begin{array}{c}
w\in E(\widehat{G})\times FP(\widehat{G}), \\ 
\exists e\in E(\widehat{G})\text{ s.t., }w\text{ is the word} \\ 
\text{only in }\{e,\text{ }e^{-1}\},\text{ satisfying} \\ 
\omega _{+}\left( \omega (w)\right) =\left( (v,\text{ }v),\text{ }0\right) ,
\\ 
\text{for some }v\in V(G)
\end{array}
\right. \right\} ,$
\end{center}

\strut

for all $n$ $\in $ $\Bbb{N}.$ First, notice that both $\mathcal{W}%
_{2n}^{(m)} $ and $\mathcal{W}_{2n}^{(c)}$ are contained in $FP(\widehat{G})$
$\cap $ $E(\widehat{G})^{*},$ where $E(\widehat{G})^{*}$ $=$ $\cup
_{n=1}^{\infty }$ $E(\widehat{G})^{n}$ ($\supseteq $ $FP(\widehat{G})$) is
the collection of all finite words in $E(\widehat{G}),$ for all $n$ $\in $ $%
\Bbb{N},$ and $FP(\widehat{G})$ is the finite path set in the free
semigroupoid $\Bbb{F}^{+}(\widehat{G})$ of the shadowed graph $\widehat{G}$
of the given graph $G.$

\strut

Now, we define the \emph{loop-set}, $loop(\widehat{G})$ contained in the set 
$FP(\widehat{G})$ $\cap $ $E(\widehat{G})^{*},$ by

\strut

\begin{center}
$loop(\widehat{G})$ $\overset{def}{=}$ $\left\{ w\in E(\widehat{G}%
)^{*}\left| 
\begin{array}{c}
w\text{ is a loop} \\ 
\text{finite path} \\ 
\text{in }FP(\widehat{G})
\end{array}
\right. \right\} $ $\subseteq $ $FP(\widehat{G}).$
\end{center}

\strut

More precisely, we can define it by

\strut

\begin{center}
$loop(\widehat{G})$ $\overset{def}{=}$ $\underset{v\in V(G)}{\cup }$ $\left(
loop_{v}(\widehat{G})\right) ,$
\end{center}

where

\begin{center}
$loop_{v}(\widehat{G})$ $\overset{def}{=}$ $\{w$ $\in $ $FP(\widehat{G})$ $:$
$w$ $=$ $v$ $w$ $v\}.$
\end{center}

\strut

i.e., the subset $loop(\widehat{G})$ of $FP(\widehat{G})$ is decomposed by $%
loop_{v}(\widehat{G}),$ for all $v$ $\in $ $V(G).$ Clearly, if a vertex $v$
of $G$ has no incident loop finite path, then $loop_{v}(\widehat{G})$ is
empty. However, since $loop_{v}(\widehat{G})$ are contained in $\Bbb{F}^{+}(%
\widehat{G})$ (not in $\Bbb{G}$), and since our graph $G$ (and $\widehat{G}$%
) is connected, whenever there is at least one edge $e$ such that either $e$ 
$=$ $v$ $e$ or $e$ $=$ $e$ $v,$ the set $loop_{v}(\widehat{G})$ is nonempty.
For instance, if $e$ $=$ $ve$ in $E(\widehat{G}),$ then $ee^{-1}$ $\in $ $%
loop_{v}(\widehat{G}).$

\strut

By the definition of $loop(\widehat{G}),$ we can re-construct the set $%
\mathcal{W}_{2n}^{(m)}$:

\strut

\begin{center}
$\mathcal{W}_{2n}^{(m)}$ $=$ $\left\{ w\in E(\widehat{G})^{2n}\left| 
\begin{array}{c}
w\in loop(\widehat{G}),\text{ satisfying} \\ 
\omega _{+}\left( \omega (w)\right) =\left( (v,\text{ }v),\text{ }0\right) ,
\\ 
\text{for some }v\in V(G)
\end{array}
\right. \right\} ,$
\end{center}

\strut

for all $n$ $\in $ $\Bbb{N}.$

\strut

Now, recall the reduction (RR) on the free semigroupoid $\Bbb{F}^{+}(%
\widehat{G}).$ We give the reduction (RR), as the following quotient map $%
\delta $ $:$ $\Bbb{F}^{+}(\widehat{G})$ $\rightarrow $ $\Bbb{G}$ (and hence
surjective), defined by

\strut

\begin{center}
$\delta (W)$ $=$ $w$ $\in $ $\Bbb{G},$ for all $W$ $\in $ $\Bbb{F}^{+}(%
\widehat{G}),$
\end{center}

satisfying that

\begin{center}
$\delta (W$ $W^{-1})$ $=$ $v_{1}$ and $\delta (W^{-1}$ $W)$ $=$ $v_{2},$
\end{center}

\strut

for all $W$ $=$ $v_{1}$ $W$ $v_{2}$ $\in $ $\Bbb{F}^{+}(\widehat{G}),$ with $%
v_{1},$ $v_{2}$ $\in $ $V(\widehat{G})$ $=$ $V(G).$ By [62], we can check
that, indeed, the map $\delta $ $:$ $\Bbb{F}^{+}(\widehat{G})$ $\rightarrow $
$\Bbb{G}$ is the well-defined quotient map. Equivalently,

$\strut $

\begin{center}
$\Bbb{G}$ $\overset{\text{Groupoid}}{=}$ $\Bbb{F}^{+}(\widehat{G})$ $/$ $%
\delta (\Bbb{F}^{+}(\widehat{G})).$
\end{center}

\strut

We say that $\delta $ is the \emph{reduction} (\emph{quotient}) \emph{map of}
$\Bbb{F}^{+}(\widehat{G})$ \emph{onto} $\Bbb{G}.$ By help of this map $%
\delta ,$ we can re-construct the set $\mathcal{W}_{2n}^{(m)}$ by

\strut

\begin{center}
$\mathcal{W}_{2n}^{(m)}$ $=$ $\left\{ w\in E(\widehat{G})^{2n}\cap loop(%
\widehat{G})\left| \delta (w)\in V(\widehat{G})\right. \right\} ,$
\end{center}

\strut

for all $n$ $\in $ $\Bbb{N}.$ Thus the set $\mathcal{W}_{2n}^{(n)}$ is
characterized, graph-theoretically, as follows.

\strut

\begin{proposition}
For $n$ $\in $ $\Bbb{N},$ we have that

\strut 

\begin{center}
$\mathcal{W}_{2n}^{(m)}$ $=$ $\left\{ w\in E(\widehat{G})^{2n}\cap loop(%
\widehat{G})\left| \delta (w)\in V(\widehat{G})\right. \right\} .$
\end{center}

$\square $
\end{proposition}

\strut

Now, we define the subset $\mathcal{D}_{loop}(\widehat{G})$ of $loop(%
\widehat{G})$ in $FP(\widehat{G})$ $\subset $ $\Bbb{F}^{+}(\widehat{G})$ by

\strut

\begin{center}
$\mathcal{D}_{loop}(\widehat{G})$ $\overset{def}{=}$ $\left\{ w\in loop(%
\widehat{G})\left| 
\begin{array}{c}
\exists \text{ }e\in E(\widehat{G})\text{ s.t.} \\ 
w\text{ is a word} \\ 
\text{only in }\{e,\text{ }e^{-1}\}
\end{array}
\right. \right\} .$
\end{center}

\strut

Then we can re-construct $\mathcal{W}_{2n}^{(c)}$ by

\strut

\begin{center}
$\mathcal{W}_{2n}^{(c)}$ $=$ $\left\{ w\in E(\widehat{G})^{2n}\left| 
\begin{array}{c}
w\in \mathcal{D}_{loop}(\widehat{G}),\text{ satisfying} \\ 
\omega _{+}\left( \omega (w)\right) =\left( (v,\text{ }v),\text{ }0\right) ,
\\ 
\text{for some }v\in V(G)
\end{array}
\right. \right\} ,$
\end{center}

\strut

for all $n$ $\in $ $\Bbb{N}.$ Therefore, we obtain the following proposition.

\strut

\begin{proposition}
For any $n$ $\in $ $\Bbb{N},$ we have that

\strut \strut 

\begin{center}
$
\begin{array}{ll}
\mathcal{W}_{2n}^{(c)} & =\left\{ w\in E(\widehat{G})^{2n}\left| 
\begin{array}{c}
w\in \mathcal{D}_{loop}(\widehat{G}),\text{ } \\ 
\text{satisfying} \\ 
\delta (w)\in V(\widehat{G})
\end{array}
\right. \right\} .
\end{array}
$
\end{center}

$\square $
\end{proposition}

\strut

The above two propositions provide the combinatorial characterization of the
sets $\mathcal{W}_{2n}^{(m)}$ and $\mathcal{W}_{2n}^{(c)},$ for $n$ $\in $ $%
\Bbb{N}.$ It means that we can verify the $D_{G}$-valued moments and
cumulants of the labeling operator $T_{G}$ of $\Bbb{G}$ on the graph Hilbert
space $H_{G},$ simply by checking the shadowed graph $\widehat{G}$ of $G$.

\strut

\strut

\strut

\section{Examples\strut}

\strut

\strut

In this section, we give several examples of the labeling operators of the
labeled graph groupoids. In Section 4, we showed that the $D_{G}$-valued
free distributional data of the labeling operator $T_{G}$ of a (labeled)
graph groupoid $\Bbb{G}$ on $H_{G}$ is characterized by the certain subsets $%
\mathcal{W}_{n}^{(m)}$ and $\mathcal{W}_{2n}^{(c)}$ of the (non-reduced)
finite path set $FP(\widehat{G})$ of the free semigroupoid $\Bbb{F}^{+}(%
\widehat{G})$, where $\widehat{G}$ is the shadowed graph of $G,$ for all $n$ 
$\in $ $\Bbb{N}.$ The $D_{G}$-valued free distributional data of the
labeling operator $T_{G}$ of $\Bbb{G}$ is represented by the $D_{G}$-valued
moments $\{E(T_{G}^{n})\}_{n\in \Bbb{N}},$ equivalently, by the $D_{G}$%
-valued cumulants $\{k_{n}(T_{G},$ ..., $T_{G})\}_{n\in \Bbb{N}},$ where $E$ 
$:$ $M_{G}$ $\rightarrow $ $D_{G}$ is the canonical conditional expectation;

\strut

\begin{center}
$E\left( T_{G}^{2n}\right) $ $=$ $\underset{w\in \mathcal{W}_{2n}^{(m)}}{%
\sum }$ $R_{w},$ for all $n$ $\in $ $\Bbb{N},$
\end{center}

and

\begin{center}
$k_{2n}\left( \underset{2n\text{-times}}{\underbrace{T_{G},\text{ ..., }T_{G}%
}}\right) $ $=$ $\underset{w\in \mathcal{W}_{2n}^{(c)}}{\sum }$ $\left( \mu
_{w}\text{ }R_{w}\right) ,$ for all $n$ $\in $ $\Bbb{N},$
\end{center}

\strut with

\begin{center}
$\strut k_{n}(T_{G},$ ..., $T_{G})$ $=$ $0_{D_{G}}$ $=$ $E(T_{G}^{n}),$
\end{center}

\strut

for all odd $n$.

\strut

\begin{example}
(\textbf{Fractaloid Cases}) Let $G$ be a canonical weighted graph with its
labeled graph groupoid $\Bbb{G},$ and assume that $\Bbb{G}$ is a fractaloid.
i.e., the automata actions $\{\mathcal{A}_{w}$ $:$ $w$ $\in $ $\Bbb{F}^{+}(%
\widehat{G})\}$ of the $G$-automaton $\mathcal{A}_{G}$ acts fully on the $%
\mathcal{A}_{G}$-tree $\mathcal{T}_{2N},$ where $\mathcal{T}_{2N}$ is the $%
(2N)$-regular tree, where $N$ $=$ $\max \{\deg _{out}(v)$ $:$ $v$ $\in $ $%
V(G)\}$ in $\Bbb{N}.$ Now, let $T_{G}$ be the labeling operator of the
fractaloid $\Bbb{G}$ on $H_{G}.$ Consider now the subset $\mathcal{W}%
_{2n}^{(m)}$ of the finite path set $FP(\widehat{G})$ of the shadowed graph $%
\widehat{G}$ of $G,$ for all $n$ $\in $ $\Bbb{N}.$ Recall that

\strut 

\begin{center}
$\mathcal{W}_{2n}^{(m)}$ $=$ $\left\{ w\in E(\widehat{G})^{2n}\left| 
\begin{array}{c}
\omega _{+}\left( \omega (w)\right) =\left( (v,\text{ }v),\text{ }0\right) 
\\ 
\text{for some }v\in V(G)
\end{array}
\right. \right\} ,$
\end{center}

\strut 

contained in $loop(\widehat{G})$ $\subseteq $ $FP(\widehat{G}).$ Moreover,
the $D_{G}$-valued moments $E\left( T_{G}^{2n}\right) $ is determined by

\strut 

\begin{center}
$E(T_{G}^{2n})$ $=$ $\underset{w\in \mathcal{W}_{2n}^{(m)}}{\sum }$ $R_{w},$
for all $n$ $\in $ $\Bbb{N}.$
\end{center}

\strut 

We can verify that, since $\Bbb{G}$ is a fractaloid, for any $n$ $\in $ $%
\Bbb{N},$ we have

\strut 

\strut (6.1)

\begin{center}
$\underset{w\in \mathcal{W}_{2n}^{(m)}}{\sum }$ $R_{w}$ $=$ $\underset{%
(l_{i_{1}},\text{ ..., }l_{i_{2n}})\in \left( \pm X_{0}^{*}\right) ^{2n}%
\text{, }\theta \left( l_{i_{1}}...l_{i_{2n}}\right) =0}{\sum }$ $%
R_{e_{1}...e_{2n}},$
\end{center}

\strut 

where $e_{1}$ ... $e_{2n}$ $\in $ $E(\widehat{G})^{2n},$ with $(\widetilde{pr%
}$ $\circ $ $\omega )(e_{j})$ $=$ $l_{i_{j}},$ for all $j$ $=$ $1,$ ..., $2n$%
. Then, again, since $\Bbb{G}$ is a fractaloid, there are $\left|
V(G)\right| $-many words $w_{1},$ ..., $w_{\left| V(G)\right| },$ in $E(%
\widehat{G})^{2n}$ (possibly, $\left| V(G)\right| $ $=$ $\infty $)
satisfying that $w_{k}$ $=$ $e_{k_{1}}$ ... $e_{k_{2n}}$ such that $(%
\widetilde{pr}$ $\circ $ $\omega )(e_{k_{q}})$ $=$ $l_{i_{q}},$ for $q$ $=$ $%
1,$ ..., $2n.$ Since $\underset{v\in V(G)}{\sum }$ $R_{v}$ is the identity
element $1_{M_{G}}$ $=$ $1_{D_{G}}$ of the right graph von Neumann algebra $%
M_{G}$, the formula (6.1) is identical to

\strut 

(6.2)

\begin{center}
$\underset{v\in V(G)}{\sum }$ $\left( \left| \sum_{2n}^{(N)}\right| \cdot 
\text{ }R_{v}\right) $ $=$ $\left| \mathcal{L}_{N}^{o}(2n)\right| $ $\left( 
\underset{v\in V(G)}{\sum }\text{ }R_{v}\right) $ $=$ $\left| \mathcal{L}%
_{N}^{o}(2n)\right| $ $\cdot $ $1_{D_{G}},$
\end{center}

\strut 

where $\mathcal{L}_{N}^{o}(2n)$ are introduced in Section 3.1, for all $n$ $%
\in $ $\Bbb{N}.$ In [17], we showed that, the cardinalities $\left| \mathcal{%
L}_{N}^{o}(2n)\right| $ are gotten recursively by the entries of the famous
Pascal's triangle (Also, see [41]).

\strut \strut 

This example shows that if $\Bbb{G}$ is a fractaloids, then the amalgamated
(or the operator-valued) free distributional data of the labeling operator
is completely determined by the scalar-values, which are the cardinalities
of certain subsets of the lattice path set $\mathcal{L}_{N}$.
\end{example}

\strut

Remember that, in [17], we conjectured that the only ``finite'' connected
directed graphs inducing fractaloids are one-flow circulant graphs and
one-vertex-multi-loop-edge graphs and the finite iterated glued graphs of
them, in the sense of [16]. This conjecture is still open. If it is
positive, then the finite fractaloids are completely characterized, and
hence we could obtain the complete analysis of their labeling operators in
terms of their amalgamated free distributional data.

\strut

\begin{example}
Let $G$ be a directed graph with

\strut 

\begin{center}
$V(G)$ $=$ $\{v_{1},$ $v_{2},$ $v_{3}\}$
\end{center}

and

\begin{center}
$E(G)$ $=$ $\{e_{12:1},$ $e_{12:2},$ $e_{13:1},$ $e_{22:1}\},$
\end{center}

\strut 

where $e_{ij:k}$ means the $k$-th edge connecting the vertex $v_{i}$ to the
vertex $v_{j}$. Then we can have the labeling set

\strut 

\begin{center}
$X$ $=$ $\{l_{1}$ $=$ $\overrightarrow{(1,\text{ }e)},$ $l_{2}$ $=$ $%
\overrightarrow{(1,\text{ }e^{2})}\}$ $\subset $ $\mathcal{L}_{2}$
\end{center}

since

\begin{center}
$\max \{\deg _{out}(v_{k})$ $:$ $k$ $=$ $1,$ $2,$ $3\}$ $=$ $2$ $=$ $\deg
_{out}(v_{1}).$
\end{center}

\strut \strut 

Then the weighting process $\omega $ $:$ $\Bbb{F}^{+}(\widehat{G})$ $%
\rightarrow $ $V(G)^{2}$ $\times $ $\pm X_{0}^{*}$ $=$ $\mathcal{X}_{0}^{*}$
is completely determined by the indices of the edges:

\strut 

\begin{center}
$\omega (e_{ij:k})$ $=$ $\left( (v_{i},\text{ }v_{j}),\text{ }l_{k}\right) .$
\end{center}

\strut and

\begin{center}
$\omega (e_{ij:k}^{-1})$ $=$ $\left( (v_{j},\text{ }v_{i}),\text{ }%
-l_{k}\right) .$
\end{center}

\strut 

By definition, we can have the $k$-th labeling operators, for $k$ $=$ $\pm 1,
$ $\pm 2,$ as elements in the right graph von Neumann algebra $M_{G}$ as
follows:

\strut 

\begin{center}
$T_{1}$ $=$ $R_{e_{12:1}}$ $+$ $R_{e_{13:1}}$ $+$ $R_{e_{22:1}}$ and $T_{-1}$
$=$ $T_{1}^{*}$
\end{center}

and

\begin{center}
$T_{2}$ $=$ $R_{e_{12:2}},$ and $T_{-2}$ $=$ $T_{2}^{*}$ $=$ $%
R_{e_{12:2}^{-1}}.$
\end{center}

and hence

\begin{center}
$T_{G}$ $=$ $\left( T_{1}+T_{-1}\right) $ $+$ $(T_{2}$ $+$ $T_{-2}).$
\end{center}

\strut 

Then we can have that $E(T_{G}^{k})$ $=$ $0_{D_{G}},$ whenever $k$ is odd,
and

\strut 

\begin{center}
$E\left( T_{G}^{2n}\right) $ $=$ $\underset{w\in \mathcal{W}_{2n}^{(m)}}{%
\sum }$ $R_{w},$ for all $n$ $\in $ $\Bbb{N}$,
\end{center}

\strut 

by Section 5.2. In our case, we have that

\strut 

\begin{center}
$\mathcal{W}_{2}^{(m)}$ $=$ $\left\{ 
\begin{array}{c}
e_{12:1}e_{12:1}^{-1},\quad e_{12:1}^{-1}e_{12:1},\quad
e_{12:2}e_{12:2}^{-1}, \\ 
e_{12:2}^{-1}e_{12:2},\quad e_{13:1}e_{13:1}^{-1},\quad e_{13:1}^{-1}e_{13:1}
\end{array}
\right\} ,$
\end{center}

\strut 

by Section 5.3. So, we can have that

\strut 

$\quad E(T_{G}^{2})$ $=$ $R_{e_{12:1}e_{12:1}^{-1}}$ $+$ $%
R_{e_{12:1}^{-1}e_{12:1}}$ $+$ $R_{e_{12:2}e_{12:2}^{-1}}$

\strut 

$\qquad \qquad \qquad +$ $R_{e_{12:2}^{-1}e_{12:2}}$ $+$ $%
R_{e_{13:1}e_{13:1}^{-1}}$ $+$ $R_{e_{13:1}^{-1}e_{13:1}}$

\strut 

$\qquad \qquad =$ $R_{v_{1}}$ $+$ $R_{v_{2}}$ $+$ $R_{v_{1}}$ $+$ $R_{v_{2}}$
$+$ $R_{v_{1}}$ $+$ $R_{v_{3}}$

\strut 

$\qquad \qquad =$ $3$ $R_{v_{1}}$ $+$ $2$ $R_{v_{2}}$ $+$ $R_{v_{3}}.$

\strut 

The above computation $E(T_{G}^{2})$ of $T_{G}$ says that the given graph $G$
does not generate a fractaloid.
\end{example}

\strut

\strut

\strut

\strut

\strut

\strut

\section{!!!!!!!!!!!!!!Put Pictures!!!!!!!!!!!!!!!!$\Bbb{G}\mathcal{G}\omega 
\Bbb{\varphi \psi }\widetilde{pr}$!!!!!!!!!!!!!!!!!!!!!!!!!!!!}

\strut

\strut

\strut

\strut

\strut

\strut

\strut

\strut \textbf{References}

\strut

\begin{quote}
\strut

{\small [1] \ \ A. G. Myasnikov and V. Shapilrain (editors), Group Theory,
Statistics and Cryptography, Contemporary Math, 360, (2003) AMS.}

{\small [2] \ \ A. Gibbons and L. Novak, Hybrid Graph Theory and Network
Analysis, ISBN: 0-521-46117-0, (1999) Cambridge Univ. Press.}

{\small [3]\strut \ \ \ \strut B. Solel, You can see the arrows in a Quiver
Operator Algebras, (2000), preprint.}

{\small [4] \ \ C. W. Marshall, Applied Graph Theory, ISBN: 0-471-57300-0
(1971) John Wiley \& Sons}

{\small [5] \ \ \strut D.Voiculescu, K. Dykemma and A. Nica, Free Random
Variables, CRM Monograph Series Vol 1 (1992).\strut }

{\small [6] \ \ D.W. Kribs and M.T. Jury, Ideal Structure in Free
Semigroupoid Algebras from Directed Graphs, preprint.}

{\small [7] \ \ D.W. Kribs, Quantum Causal Histories and the Directed Graph
Operator Framework, arXiv:math.OA/0501087v1 (2005), Preprint.}

{\small [8] \ \ F. Balacheff, Volume Entropy, Systole and Stable Norm on
Graphs, arXiv:math.MG/0411578v1, (2004) Preprint.}

{\small [9] \ \ G. C. Bell, Growth of the Asymptotic Dimension Function for
Groups, (2005) Preprint.}

{\small [10]\ I. Cho, Graph von Neumann algebras, ACTA. Appl. Math, 95,
(2007) 95 - 135.}

{\small [11]\ I. Cho, Characterization of Free Blocks of a right graph von
Neumann algebra, Compl. An. \& Op. theo (2007) To be Appeared.}

{\small [12]\ I. Cho, Direct Producted }$W^{*}${\small -Probability Spaces
and Corresponding Free Stochastic Integration, B. of KMS, 44, No. 1, (2007),
131 - 150.}

{\small [13] I. Cho, Vertex-Compressed Algebras of a Graph von Neumann
Algebra, (2007) Submitted to ACTA. Appl. Math.}

{\small [14] I. Cho, Group-Freeness and Certain Amalgamated Freeness, J. of
KMS, 45, no. 3, (2008) 597 - 609.}

{\small [15]\ I. Cho and P. E. T. Jorgensen, }$C^{*}${\small -Algebras
Generated by Partial Isometries, JAMC, (2008) To Appear.}

{\small [16] I. Cho and P. E. T. Jorgensen, }$C^{*}${\small -Subalgebras
Generated by Partial Isometries, JMP, (2008) To Appear.}

{\small [17] I. Cho and P. E. T. Jorgensen, Applications of Automata and
Graphs: Labeling-Operators in Hilbert Space I, (2008) Submitted to ACTA
Appl. Math., Special Issues.}

{\small [18]\ I. Raeburn, Graph Algebras, CBMS no 3, AMS (2005).}

{\small [19]\ P. D. Mitchener, }$C^{*}${\small -Categories, Groupoid
Actions, Equivalent KK-Theory, and the Baum-Connes Conjecture,
arXiv:math.KT/0204291v1, (2005), Preprint.}

{\small [20] R. Scapellato and J. Lauri, Topics in Graph Automorphisms and
Reconstruction, London Math. Soc., Student Text 54, (2003) Cambridge Univ.
Press.}

{\small [21] R. Exel, A new Look at the Crossed-Product of a }$C^{*}${\small %
-algebra by a Semigroup of Endomorphisms, (2005) Preprint.}

{\small [22] R. Gliman, V. Shpilrain and A. G. Myasnikov (editors),
Computational and Statistical Group Theory, Contemporary Math, 298, (2001)
AMS.}

{\small [23] R. Speicher, Combinatorial Theory of the Free Product with
Amalgamation and Operator-Valued Free Probability Theory, AMS Mem, Vol 132 ,
Num 627 , (1998).}

{\small [24] S. H. Weintraub, Representation Theory of Finite Groups:
Algebra and Arithmetic, Grad. Studies in Math, vo. 59, (2003) AMS.}

{\small [25] V. Vega, Finite Directed Graphs and }$W^{*}${\small %
-Correspondences, (2007) Ph. D thesis, Univ. of Iowa.}

{\small [26] W. Dicks and E. Ventura, The Group Fixed by a Family of
Injective Endomorphisms of a Free Group, Contemp. Math 195, AMS.}

{\small [27] F. Radulescu, Random Matrices, Amalgamated Free Products and
Subfactors of the von Neumann Algebra of a Free Group, of Noninteger Index,
Invent. Math., 115, (1994) 347 - 389.}

{\small [28] D. A. Lind, Entropies of Automorphisms of a Topological Markov
Shift, Proc. AMS, vo 99, no 3, (1987) 589 - 595.}

{\small [29] D. A. Lind and B. Marcus, An Introduction to Symbolic Dynamics
and Coding, (1995) Cambridge Univ. Press.}

{\small [30] D. A. Lind and S. Tuncel, A Spanning Tree Invariant for Markov
Shifts, IMA Vol. Math. Appl., vo 123, (2001), 487 - 497.}

{\small [31] D. A. Lind and K. Schmidt, Symbolic and Algebraic Dynamical
Systems, Handbook of Dynamical System, Vol.\TEXTsymbol{\backslash}1A, (2002)
765 - 812.}

{\small [32] R. V. Kadison and J. R. Ringrose, Fundamentals of the Theory of
Operator Algebra, Grad. Stud. Math., vo. 15, (1997) AMS.}

{\small [33] D. E. Dutkay and P. E. T. Jorgensen, Iterated Function Systems,
Ruelle Operators and Invariant Projective Measures,
arXiv:math.DS/0501077/v3, (2005) Preprint.}

{\small [34] P. E. T. Jorgensen, Use of Operator Algebras in the Analysis of
Measures from Wavelets and Iterated Function Systems, (2005) Preprint.}

{\small [35] D. Guido, T. Isola and M. L. Lapidus, A Trace on Fractal Graphs
and the Ihara Zeta Function, arXiv:math.OA/0608060v1, (2006) Preprint.}

{\small [36] P. Potgieter, Nonstandard Analysis, Fractal Properties and
Brownian Motion, arXiv:math.FA/0701649v1, (2007) Preprint.}

{\small [37] L. Bartholdi, R. Grigorchuk, and V. Nekrashevych, From Fractal
Groups to Fractal Sets, arXiv:math.GR/0202001v4, (2002) Preprint.}

{\small [38] I. Cho, The Moments of Certain Perturbed Operators of the
Radial Operator of the Free Group Factor }$L(F_{N})${\small , JAA, 5, no. 3,
(2007) 137 - 165.}

{\small [39] I. Cho and P. E. T. Jorgensen, }$C^{*}${\small -Subalgebras
Generated by Partial Isometries in }$B(H)${\small , JMP (2008) To Appear.}

{\small [40] S. Thompson and I. Cho, Powers of Mutinomials in Commutative
Algebras, (2008) (Undergraduate Research) Submitted to PMEJ.}

{\small [41] S. Thompson, C. M. Mendoza, and A. J. Kwiatkowski, and I. Cho,
Lattice Paths Satisfying the Axis Property, (2008) (Undergraduate Research)
Preprint.}

{\small [42] I. Cho, Labeling Operators of Graph Groupoids, (2008) Preprint.}

{\small [43] R. T. Powers, Heisenberg Model and a Random Walk on the
Permutation Group, Lett. Math. Phys., 1, no. 2, (1975) 125 - 130.}

{\small [44] R. T. Powers, Resistance Inequalities for }$KMS${\small -states
of the isotropic Heisenberg Model, Comm. Math. Phys., 51, no. 2, (1976) 151
- 156.}

{\small [45] R. T. Powers, Registance Inequalities for the Isotropic
Heisenberg Ferromagnet, JMP, 17, no. 10, (1976) 1910 - 1918.}

{\small [46] E. P. Wigner, Characteristic Vectors of Bordered Matrices with
Infinite Dimensions, Ann. of Math. (2), 62, (1955) 548 - 564.}

{\small [47] D. Voiculescu, Symmetries of Some Reduced Free Product }$C^{*}$%
{\small -Algebras, Lect. Notes in Math., 1132, Springer, (1985) 556 - 588.}

{\small [48] T. Shirai, The Spectrum of Infinite Regular Line Graphs, Trans.
AMS., 352, no 1., (2000) 115 - 132.}

{\small [49] J. Kigami, R. S. Strichartz, and K. C. Walker, Constructing a
Laplacian on the Diamond Fractal, Experiment. Math., 10, no. 3, (2001) 437 -
448.}

{\small [50] I. V. Kucherenko, On the Structurization of a Class of
Reversible Cellular Automata, Diskret. Mat., 19, no. 3, (2007) 102 - 121.}

{\small [51] J. L. Schiff, Cellular Automata, Discrete View of the World,
Wiley-Interscience Series in Disc. Math .\& Optimazation, ISBN:
978-0-470-16879-0, (2008) John Wiley \& Sons Press.}

{\small [52] P. E. T. Jorgensen, and M. Song, Entropy Encoding, Hilbert
Spaces, and Kahunen-Loeve Transforms, JMP, 48, no. 10, (2007)}

{\small [53] P. E. T. Jorgensen, L. M. Schmitt, and R. F. Werner, }$q$%
{\small -Canonical Commutation Relations and Stability of the Cuntz Algebra,
Pac. J. of Math., 165, no. 1, (1994) 131 - 151. }

{\small [54] A. Gill, Introduction to the Theory of Finite-State Machines,
MR0209083 (34\TEXTsymbol{\backslash}\#8891), (1962) McGraw-Hill Book Co.}

{\small [55] J. E. Hopcroft, and J. D. Ullman, Introduction to Automata
Theory, Language, and Computation, ISBN: 0-201-02988-X, (1979)
Addision-Wesley Publication Co.}

{\small [56] M. Fannes, B. Nachtergaele, and R. F. Werner, Ground States of }%
$VSB${\small -Models on Cayley Trees, J. of Statist. Phys., 66, (1992) 939 -
973.}

{\small [57] M. Fannes, B. Nachtergaele, and R. F. Werner, Finitely
Correlated States on Quantum Spin Chains, Comm. Math. Phys., 144, no. 3,
(1992) 443 - 490.}

{\small [58] M. Fannes, B. Nachtergaele, and R. F. Werner, Entropy Estimates
for Finitely Correlated States, Ann. Inst. H. Poincare. Phys. Theor., 57, no
3, (1992) 259 - 277.}

{\small [59] M. Fannes, B. Nachtergaele, and R. F. Werner, Finitely
Correlated Pure States, J. of Funt. Anal., 120, no 2, (1994) 511 - 534.}

{\small [60] J. Renault, A Groupoid Approach to }$C^{*}${\small -Algebras,
Lect. Notes in Math., 793, ISBN: 3-540-09977-8, (1980) Springer.}

{\small [61] S. Sakai, }$C^{*}${\small -Algebras and }$W^{*}${\small %
-Algebras, MR number: MR0442701, (1971) Springer-Verlag.}

{\small [62] I. Cho, Measures on Graphs and Groupoid Measures, Compl. An.
Oper. Theor., 2, (2008) 1 - 28.}

{\small [63] B. Fuglede, and R. V. Kadison, On Determinants and a Property
of the Trace in Finite FActors, Proc. Nat. Acad. Sci., U. S. A., 37, (1951)
425 - 431.}

{\small [64] D. Aldous, and R. Lyons, Processes on Unimodular Random
Networks, Elec. J. of Probab., 12, no 54, (2007) 1454 - 1508.}

{\small [65] O. G\"{u}hhne, G. T\'{o}tth, P. Hyllus, and H. J. Briegel, Bell
Inequalities for Graphs States, Phys. Rev. Lett., 95, no 12, (2005) }

{\small [66] R. V. Kadison, and J. R. Ringrose, Fundamentals of the Theory
of Operator Algebras, vol II, Grad. Stud. in Math., 16, ISBN: 0-8218-0820-6,
(1997) AMS.}
\end{quote}

\end{document}